\documentclass[preprint]{elsarticle} 
\usepackage{stanli}
\usetikzlibrary{decorations.pathreplacing}
\makeatletter
	\def\ps@pprintTitle{%
 	\let\@oddhead\@empty
	\let\@evenhead\@empty
	\def\@oddfoot{\centerline{\thepage}}%
	\let\@evenfoot\@oddfoot}
\makeatother

\usepackage{etoolbox}
\patchcmd{\MaketitleBox}{\footnotesize\itshape\elsaddress\par\vskip36pt}{\footnotesize\itshape\elsaddress\par\parbox[b][36pt]{\linewidth}{\vfill\hfill\textnormal{\today}\hfill\null\vfill}}{}{}%
\patchcmd{\pprintMaketitle}{\footnotesize\itshape\elsaddress\par\vskip36pt}{\footnotesize\itshape\elsaddress\par\parbox[b][36pt]{\linewidth}{\vfill\hfill\textnormal{\today}\hfill\null\vfill}}{}{}%
\usepackage{tikz}
\usepackage{tikz-imagelabels}
\usepackage{icomma}
\usepackage{pgfplots} 
\usepackage{pgfgantt}
\usepackage{pdflscape}
\pgfplotsset{compat=newest} 
\pgfplotsset{plot coordinates/math parser=false} 
\newlength\fwidth
\newlength\fheight
\usepackage[
            colorlinks=true,%
            breaklinks=true,%
            linkcolor=blue,
            urlcolor=blue,%
            citecolor=blue,%
            pdftitle={Preserving Lagrangian structure in data-driven reduced-order modeling of large-scale dynamical systems},  	
            pdfkeywords={Structure-preserving model reduction; data-driven modeling; Lagrangian dynamics; scientific machine learning; operator inference},	 
            pdfauthor={Harsh Sharma and Boris Kramer},		 %
            bookmarksopen=false,%
            pdfpagemode=None]{hyperref}
            
\usepackage[margin = 1.2in]{geometry}
\usepackage[T1]{fontenc}
\usepackage[english]{babel}
\usepackage[latin1]{inputenc}
\usepackage{mathtools}
\usepackage{amsmath}
\usepackage{amsfonts}
\usepackage{amsthm}
\usepackage{amssymb}
\usepackage{array}
\usepackage{algorithm} 
\usepackage{subcaption}
\usepackage{siunitx}
\usepackage[noend]{algpseudocode}

\usepackage{color}
\usepackage{url}
\usepackage{cleveref}
\usepackage{graphicx}
\usepackage{breakcites}
\usepackage{soul} 
\usepackage{enumitem}
\usepackage[title,titletoc,toc]{appendix}
\usepackage{siunitx}

\usepackage[textsize=tiny]{todonotes}

\theoremstyle{definition}
\newtheorem{remark}{Remark}
{\noindent {\textbf{Proof}:} }%
{\hfill $\Box$ \\[1ex] }

\newcommand{\bit}{\begin{itemize}}
\newcommand{\eit}{\end{itemize}}
\newcommand{\ben}{\begin{enumerate}}
\newcommand{\een}{\end{enumerate}}

\newcommand {\real} {\mathbb{R}}

\newcommand{\trace}[1]{\ensuremath{\mathop{\mathrm{tr}}\left( #1 \right)}} 

\newcommand{\bzero}{\ensuremath{\mathbf{0}}} 
\newcommand{\Q}{\ensuremath{\mathsf{Q}}}
\newcommand{\T}{\ensuremath{\mathsf{T}}}

\newcommand{\tL}{\ensuremath{\tilde{L}}}

\newcommand{\bA}{\ensuremath{\mathbf{A}}}
\newcommand{\bB}{\ensuremath{\mathbf{B}}}
\newcommand{\bC}{\ensuremath{\mathbf{C}}}
\newcommand{\bD}{\ensuremath{\mathbf{D}}}
\newcommand{\bE}{\ensuremath{\mathbf{E}}}
\newcommand{\bF}{\ensuremath{\mathbf{F}}}

\newcommand{\bK}{\ensuremath{\mathbf{K}}}

\newcommand{\bM}{\ensuremath{\mathbf{M}}}

\newcommand{\bQ}{\ensuremath{\mathbf{Q}}}

\newcommand{\bU}{\ensuremath{\mathbf{U}}}
\newcommand{\bV}{\ensuremath{\mathbf{V}}}
\newcommand{\bW}{\ensuremath{\mathbf{W}}}

\newcommand{\bY}{\ensuremath{\mathbf{Y}}}

\renewcommand{\bf}{\ensuremath{\mathbf{f}}}

\newcommand{\bq}{\ensuremath{\mathbf{q}}}

\newcommand{\bu}{\ensuremath{\mathbf{u}}}
\newcommand{\bv}{\ensuremath{\mathbf{v}}}
\newcommand{\bw}{\ensuremath{\mathbf{w}}}

\newcommand{\by}{\ensuremath{\mathbf{y}}}
\newcommand{\bz}{\ensuremath{\mathbf{z}}}
\newcommand{\bqd}{\ensuremath{\mathbf{\dot{q}}}}

\newcommand {\bXi} {\mbox{\boldmath $\Xi$}}

\newcommand{\hK}{\hat{\bK}}
\newcommand{\hL}{\hat{L}}

\newcommand{\hQ}{\hat{\mathbf{Q}}}

\newcommand{\hU}{\ensuremath{\hat{U}}}

\newcommand{\hY}{\hat{\mathbf{Y}}}

\newcommand{\hq}{\hat{\bq}}

\newcommand{\hbq}{\ensuremath{\mathbf{\hat q}}}

\newcommand{\hbF}{\ensuremath{\mathbf{\hat F}}}
\newcommand{\hbqd}{\ensuremath{\mathbf{\dot{\hat q}}}}
\newcommand{\hbqdd}{\ensuremath{\mathbf{\ddot{\hat q}}}}
\newcommand{\hbQ}{\ensuremath{\mathbf{\hat Q}}}
\newcommand{\hbQd}{\ensuremath{\mathbf{\dot{\hat Q}}}}
\newcommand{\hbQdd}{\ensuremath{\mathbf{\ddot{\hat Q}}}}
\newcommand{\hbf}{\ensuremath{\mathbf{\hat f}}}
\newcommand{\hbM}{\ensuremath{\mathbf{\hat M}}}
\newcommand{\hbK}{\ensuremath{\mathbf{\hat K}}}
\newcommand{\hbC}{\ensuremath{\mathbf{\hat C}}}
\newcommand{\hbE}{\ensuremath{\mathbf{\hat E}}}
\newcommand{\hbB}{\ensuremath{\mathbf{\hat B}}}

\newcommand{\tbq}{\ensuremath{\mathbf{\tilde q}}}
\newcommand{\tbqd}{\ensuremath{\mathbf{\dot{\tilde q}}}}
\newcommand{\tbqdd}{\ensuremath{\mathbf{\ddot{\tilde q}}}}
\newcommand{\tbf}{\ensuremath{\mathbf{\tilde f}}}

\newcommand{\tbM}{\ensuremath{\mathbf{\tilde M}}}
\newcommand{\tbK}{\ensuremath{\mathbf{\tilde K}}}
\newcommand{\tbC}{\ensuremath{\mathbf{\tilde C}}}
\newcommand{\tbE}{\ensuremath{\mathbf{\tilde E}}}
\newcommand{\tbB}{\ensuremath{\mathbf{\tilde B}}}

\graphicspath{{./figures/}}

\begin{document}
	
	\begin{frontmatter}
		
		\title{Preserving Lagrangian structure in data-driven reduced-order modeling of large-scale dynamical systems}
		
		
		\author[ucsd]{Harsh Sharma \corref{cor1}}
		\ead{hasharma@ucsd.edu}
		
		\author[ucsd]{Boris Kramer}

		\cortext[cor1]{Corresponding author}

		\address[ucsd]{Department of Mechanical and Aerospace Engineering, University of California San Diego, CA, United States}
		
		\begin{abstract}
			This work presents a nonintrusive physics-preserving method to learn reduced-order models (ROMs) of Lagrangian systems, which includes nonlinear wave equations. Existing intrusive projection-based model reduction approaches construct structure-preserving Lagrangian ROMs by projecting the Euler-Lagrange equations of the full-order model (FOM) onto a linear subspace. This Galerkin projection step requires complete knowledge about the Lagrangian operators in the FOM and full access to manipulate the computer code. In contrast, the proposed Lagrangian operator inference approach embeds the mechanics into the operator inference framework to develop a data-driven model reduction method that preserves the underlying Lagrangian structure. The proposed approach exploits knowledge of the governing equations (but not their discretization) to define the form and parametrization of a Lagrangian ROM which can then be learned from projected snapshot data. The method does not require access to FOM operators or computer code. The numerical results demonstrate Lagrangian operator inference on an Euler-Bernoulli beam model, the sine-Gordon (nonlinear) wave equation, and a large-scale discretization of a soft robot fishtail with $779,232$ degrees of freedom. The learned Lagrangian ROMs generalize well, as they can accurately predict the physical solutions both far outside the training time interval, as well as for unseen initial conditions. 
\end{abstract}	
		
		
\end{frontmatter}

\section{Introduction} \label{sec:intro}

In today's world of ever-increasing computational power, many engineering disciplines and the physical sciences rely on numerical simulations for design and control of complex dynamical systems. Modeling and simulation of dynamical systems using the Lagrangian mechanics framework has become essential in diverse areas such as structural mechanics, aerospace engineering, biomedical engineering, high-energy physics, quantum mechanics, solid-state physics, and soft robotics. Lagrangian systems exhibit physically interpretable quantities such as momentum, energy, or vorticity; the behavior of these quantities in numerical simulation provides an important measure of accuracy of the model.

\par

During the past three decades, significant advances have been made in the field of structure-preserving numerical methods \cite{hairer2006geometric,sharma2020review} for Hamiltonian and Lagrangian systems. In fact, the class of mechanical integrators \cite{wendlandt1997mechanical} has been developed specifically for Lagrangian systems to ensure that the numerical solution captures the underlying physics accurately. However, numerical simulations of large-scale (e.g., with thousands to millions of states) Lagrangian systems using these structure-preserving methods can take days or weeks on standard computing workstations. As a result, there is a great need for computational savings in time-critical applications such as structural design optimization, real-time simulation, control, and uncertainty quantification. This need has produced several different frameworks for faster simulation, and we next review related existing approaches in reduced-order modeling and learning of structured Lagrangian systems. 

\par 
The equations of motion for Lagrangian systems are given by the Euler-Lagrange equations. Since these equations are second-order differential equations, a straightforward approach is to rewrite the Lagrangian system as a first-order system and then using existing ROM techniques to reduce their dimensionality. However, this approach destroys the Lagrangian structure, see \cite{lall2003structure} for an example. Thus, other approaches were developed for this situation. The first approach to deriving ROMs for second-order systems goes back to 1960s where the modal truncation method was extended to the second-order setting \cite{guyan1965reduction,craig1968coupling}. The dominant pole algorithms for large-scale second-order systems were developed in \cite{rommes2008computing}. Moreover, standard model reduction methods based on balanced truncation were extended to linear systems that are second order in time in \cite{meyer1996balancing} and general second-order systems in \cite{chahlaoui2006second,reis2008balanced}. Recently a variety of data-driven approaches for second-order systems have been developed, e.g. the interpolatory Loewner framework \cite{schulze2018data,pontes2022data} and vector fitting \cite{werner2022structured}. However, most of these articles focus on frequency-response data and do not derive ROMs from time-domain data.  

\par

While the field of symplectic model reduction of Hamiltonian systems  has grown considerably in recent years \cite{peng2016symplectic,afkham2017structure,gong2017structure,pagliantini2021dynamical,buchfink2021symplectic}, progress on structure-preserving model reduction of Lagrangian systems has been less rapid. The intrusive model reduction for Lagrangian systems was introduced in \cite{lall2003structure} where the authors showed that performing a Galerkin projection on the Euler-Lagrange equations preserves the Lagrangian structure. Building on this idea, the work in \cite{carlberg2015preserving} presented an efficient structure-preserving model reduction strategy for nonlinear Lagrangian systems with parameter dependence. The authors applied their method to a geometrically nonlinear parametrized truss structure with 3,000 degrees of freedom (DOFs) in the FOM. These structure-preserving model reduction approaches are intrusive in that they assume full knowledge about governing equations and their space-time discretization, and require access to full model operators in order to derive Lagrangian reduced-order models (ROMs) via intrusive projection. This limits the scope of intrusive approaches, as in many situations, the full model operators are either not accessible or the complexity of the FOM source code makes the process of obtaining the full model operators very labor-intensive. In contrast, we propose a lightweight nonintrusive method to construct a Lagrangian ROM directly from data.

\par 
 
A variety of machine learning papers have developed structure-preserving neural networks for Lagrangian systems by endowing neural networks with \textit{physics-motivated} inductive biases, e.g. Deep Lagrangian Networks (DeLaNs) \cite{lutter2019deep}, Lagrangian Neural Networks (LNNs) \cite{cranmer2020lagrangian}, Structured Mechanical Models (SMMs) \cite{gupta2020structured}. Although these structure-preserving learning methods have been applied to various learning and control tasks, a majority of these approaches are only concerned with learning Lagrangian systems when the data is coming from very low-dimensional systems, i.e. 3-4 dimensions. 

\par

In another research direction, symbolic regression has been employed to discover governing equations and conservation laws for low-dimensional systems in \cite{schmidt2009distilling,bongard2007automated}. Sparse identification of nonlinear dynamics (SINDy) \cite{brunton2016discovering} and its modifications \cite{rudy2017data,kaheman2020sindy} have also been developed for discovering conservation laws \cite{kaiser2018discovering}, Hamiltonians \cite{kaiser2021data}, and Lagrangians \cite{chu2020discovering} from data. These approaches discover governing equations from a dictionary of candidate functions that are chosen either by expert knowledge or through sparse approximation techniques. The number of candidate functions grow factorially with the state dimension, and hence, SINDy-based methods are ill-suited for large-scale dynamical systems unless additional knowledge about low-dimensionality is used.

\par
A promising method for learning interpretable low-dimensional models from high-dimensional data is the operator inference framework~\cite{peherstorfer2016data}. This nonintrusive method for data-driven model reduction is applicable to FOMs with linear or low-order polynomial nonlinear terms and can be extended to a broader class of dynamical systems with nonpolynomial nonlinear terms via lifting transformations \cite{SKHW2020_learning_ROMs_combustor,qian2019transform,QKPW2020_lift_and_learn,khodabakhshi2021non}. The approach has also been extended to a gray-box setting in \cite{benner2020operator} where knowledge about the nonpolynomial nonlinear terms in analytic form are used to learn ROMs. However, operator inference applied in its standard form to Lagrangian systems does not preserve the underlying geometric structure and therefore produces systems whose energy grows in time (see Section~\ref{sec:motivation}), rendering them unphysical. We recently developed a structure-preserving operator inference method, Hamiltonian operator inference, for canonical Hamiltonian systems in \cite{sharma2022hamiltonian}. That method is designed for Hamiltonian FOMs obtained from spatial discretizations of conservative Hamiltonian PDEs. In contrast, the present work focuses on developing a structure-preserving operator inference method for Lagrangian mechanical systems with external nonconservative forcing and spatial discretizations of Lagrangian PDEs. The preservation of the underlying Lagrangian structure in such systems leads to new constraints on the operator inference problem which can not be tackled with the Hamiltonian operator inference method.

\par 

The main goal of this work is to develop a structure-preserving model reduction method that can learn Lagrangian ROMs nonintrusively from high-dimensional data. We focus particularly on large-scale models that arise from semi-discretization of partial differential equations (PDEs), which arise for instance in structural dynamics, nonlinear wave equations, and soft robotics applications. We approach this problem by proposing the nonintrusive \textit{Lagrangian operator inference} (\textsf{L-OpInf}), a structure-preserving data-driven model reduction method for Lagrangian systems that preserves the underlying geometric structure. The method can work with high-dimensional snapshot data from a Lagrangian system. We first project this data onto a low-dimensional basis as in classical projection-based model reduction and then obtain reduced time-derivative data using an appropriate finite difference scheme. We postulate a form and parametrization of a Lagrangian ROM and then learn the reduced operators from the reduced data using a constrained operator inference procedure. This ensures that the ROMs preserve the Lagrangian structure. The main contributions of this work are:
\begin{enumerate}
\item We develop a nonintrusive physics-preserving method to learn Lagrangian ROMs of large-scale models derived from spatial discretization of Lagrangian PDEs. The proposed method exploits knowledge about the space-time continuous Lagrangian at the PDE level, specifically of the nonlinear potential energy terms, to define and parametrize a Lagrangian ROM form. The ROM can then be learned from trajectory data via a constrained linear least-squares problem;
\item We present numerical results that demonstrate the learned models' ability to provide accurate predictions outside the training time interval for the Euler-Bernoulli beam model and the nonlinear sine-Gordon equation.~Unlike the structure-preserving Hamiltonian approaches that require both trajectory and momentum data, the presented structure-preserving approach can learn accurate and stable ROMs with bounded energy error from high-dimensional trajectory data alone;
\item We learn Lagrangian ROMs for a high-dimensional soft-robotic fishtail model with dissipation and time-dependent control input to demonstrate the proposed method's versatility and robustness to unknown control inputs.
\end{enumerate}

\par 
This paper is organized as follows. Section \ref{sec:background} reviews the basics of Lagrangian FOMs and describes projection-based intrusive structure-preserving model reduction for Lagrangian mechanical systems and nonlinear wave equations. Section \ref{sec:lopinf} presents the proposed structure-preserving operator inference for deriving Lagrangian ROMs. In Section \ref{sec:numerical} we present numerical experiments where reduced Lagrangian models are inferred from data of conservative and forced Lagrangian systems. In particular, we apply the proposed structure-preserving operator inference method to a linear beam model from structural dynamics, a nonlinear space-discretized sine-Gordon equation, and a fishtail model from soft robotics. Finally, Section \ref{sec:conclusions} summarizes the contributions and suggests future research directions.
\section{Background} \label{sec:background}
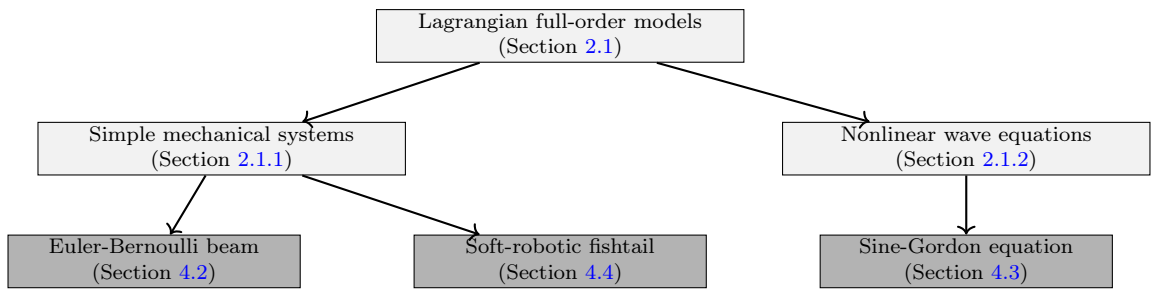
\begin{figure}[h]
\centering
\begin{tikzpicture}
[node distance = 0.5cm, auto,font=\footnotesize,
every node/.style={node distance=1.5cm},
comment/.style={rectangle, inner sep= 5pt, text width=5cm, node distance=0.25cm, font=\scriptsize\sffamily},
force2/.style={rectangle, draw, fill=black!5, inner sep=1.2pt, text width=4.8cm, text badly centered, minimum height=0.7cm, font=\footnotesize},
force/.style={rectangle, draw, fill=black!30, inner sep=1.2pt, text width=3.8cm, text badly centered, minimum height=0.7cm, font=\footnotesize}] 


\node [force2] (non) {Nonlinear wave equations\\ (Section \ref{sec:non})};
\node[force2, above of=non,left=0.5 cm of non](FOM){Lagrangian full-order models\\ (Section \ref{sec:FOM})};
\node [force2, left=5 cm of non] (simp) {Simple mechanical systems \\ (Section \ref{sec:simple})};
\node [force, below of=non] (sg) {Sine-Gordon equation \\ (Section \ref{sec:sg}) };
\node [force, left=1.5 cm of sg] (fish) {Soft-robotic fishtail \\ (Section \ref{sec:fishtail})};
\node [force, left=1.5cm of fish] (eb) {Euler-Bernoulli beam \\ (Section \ref{sec:beam})};


;


;





\path[->,thick] 
(FOM) edge (non)
(FOM) edge (simp)
(simp) edge (fish)
(simp) edge (eb)
(non) edge (sg);
\end{tikzpicture} 
\caption{Schematic overview of the types of Lagrangian FOMs considered in this work and the respective application problems (dark grey).}
\label{fig:schematic}
\end{figure}

In Section \ref{sec:FOM}, we introduce the Lagrangian FOMs considered herein (see Figure~\ref{fig:schematic}), followed by its structure-preserving time integration using variational integrators in Section \ref{sec:gni}. In Section \ref{sec:intrusive}, we derive the projection-based intrusive Lagrangian ROMs.  
\subsection{Lagrangian full-order models} \label{sec:FOM}
We consider a Lagrangian system with a finite-dimensional configuration manifold $\Q$, state space $\T\Q$ and a Lagrangian $L: \T\Q \to \real$. For a conservative system Hamilton's principle (see, e.g., \cite{arnol2013mathematical}) characterizes the solution trajectory $\bq(t)$ which passes through $\bq(t_0)=\bq_0$ at initial time $t=t_0$ to reach $\bq(t_K)=\bq_K$ at the final time $t_K>t_0$ as that which satisfies the variational principle
\begin{equation}
    \delta \mathfrak B[\bq]=\delta \int^{t_K}_{t_0} L(\bq, \dot{\bq}) \  \text{d}t = 0.
\end{equation}
Here, $\mathfrak B[\bq]$ is the action functional and the variation of the action functional $\mathfrak B$ at a trajectory $\bq$ for a variation $\delta \bq$ is defined as
\begin{equation}
    \langle \delta \mathfrak B[\bq], \delta \bq \rangle= \frac{\text d}{\text d\epsilon}\mathfrak B [\bq^\epsilon] \bigg|_{\epsilon=0},
\end{equation}
where $\bq^\epsilon$ is any of the one-parameter families that defines the variation $\delta \bq$.
The corresponding Euler-Lagrange equations are 
\begin{equation}
    \frac{\partial L(\bq,\dot{\bq})}{\partial \bq} -\frac{\text d}{\text dt} \left( \frac{\partial L(\bq,\dot{\bq})}{\partial \dot{\bq}} \right) = \bzero.
    \label{eq:el}
\end{equation}
Conservative dynamical systems governed by the Euler-Lagrange equations exhibit important qualitative properties. For autonomous Lagrangian systems the total  energy $E(\bq,\dot{\bq})=\frac{\partial L}{\partial \dot{\bq}}\cdot \dot{\bq} - L(\bq,\dot{\bq})$ is conserved along solution trajectories. In addition, by Noether's theorem \cite{noether1971invariant}, there exists an invariant of the motion corresponding to each symmetry of the system Lagrangian. The flow map for these Lagrangian systems also preserves the symplectic Lagrangian form.
\par 

For Lagrangian systems with external nonconservative forcing $\bf(\bq,\dot{\bq},t)$, the Lagrange-d'Alembert principle (see, e.g., \cite{arnol2013mathematical}) seeks $\bq(t)$ satisfying
\begin{equation}
        \delta \mathfrak B(\bq)=\delta \int^{t_K}_{t_0} L(\bq, \dot{\bq}) \  \text{d}t + \int^{t_K}_{t_0}\bf(\bq,\dot{\bq},t) \cdot \delta \bq  \ \text{d}t = 0,
        \label{eq:lda}
\end{equation}
which leads to the forced Euler-Lagrange equations 
\begin{equation}
         \frac{\partial L(\bq,\dot{\bq})}{\partial \bq} -\frac{\text d}{\text dt} \left( \frac{\partial L(\bq,\dot{\bq})}{\partial \dot{\bq}} \right) + \bf(\bq,\dot{\bq},t)  = \bzero.
         \label{eq:fel}
\end{equation}

For this work we focus on two types of Lagrangian FOMs. First, we consider \textit{simple mechanical systems} under Rayleigh viscous damping and external time-dependent forcing in Section~\ref{sec:simple}. Second, we consider spatial discretizations of nonlinear Lagrangian PDEs in Section~\ref{sec:non}.
\subsubsection{Mechanical systems with nonconservative forcing} 
\label{sec:simple}
Simple mechanical systems are defined in terms of three components: a configuration manifold $\Q$, a Riemannian metric $g(\dot{\bq},\dot{\bq})$ where $\dot{\bq}$ belongs to the tangent bundle of $\Q$, i.e. $\T\Q$, and a scalar potential function $U(\bq)$ defined on the configuration manifold $\Q$. See \cite{bullo2019geometric} for more details. In particular, we consider simple mechanical systems with configuration manifold $\Q=\real^n$ where $n$ denotes the large number of degrees of freedom in the FOM. The Lagrangian $L(\bq,\dot{\bq})=T(\dot{\bq})-U(\bq)$ for these systems represents the difference between the kinetic energy $T(\dot{\bq})$ and potential energy $U(\bq)$. The notion of Riemannian metric plays a key role in describing the kinetic energy of simple mechanical systems. In geometric terms, a Riemannian metric $g(\dot{\bq}_1,\dot{\bq}_2)$ is a positive-definite, symmetric covariant 2-tensor field on $\Q$, where $\dot{\bq}_1,\dot{\bq}_2 \in \T\Q$. The kinetic energy $T(\dot{\bq})$ of simple mechanical systems can be written in terms of the Riemannian metric as 
\begin{equation}
        T(\dot{\bq})=\frac{1}{2}g(\dot{\bq},\dot{\bq})= \frac{1}{2}\dot{\bq}^\top\bM\dot{\bq},
\end{equation}
where $\bM \in \real^{n \times n}$ is the symmetric positive-definite mass matrix, i.e. $\bM=\bM^\top \succ 0$. The Lagrangian for simple mechanical systems can be expressed as
\begin{equation}
        L(\bq,\dot{\bq})= T(\dot{\bq}) - U(\bq)=\frac{1}{2}\dot{\bq}^\top\bM\dot{\bq} - \frac{1}{2}\bq^\top \bK \bq,
        \label{eq:Lagrangian}
\end{equation}
where $\bK$ is the stiffness matrix. As a result of structure-preserving space discretization via finite element methods, the stiffness matrix is often symmetric and positive-definite, i.e. $\bK=\bK^\top\succ 0$. 
If the system is subjected to a nonconservative external forcing, i.e., $\bf \neq \bzero$, that force is often modeled via a dissipative force and an external time-dependent input as
\begin{equation}
        \bf(\bq,\dot{\bq},t)=-\bC \dot{\bq} + \bB \bu(t),
        \label{eq:fn}
\end{equation}
where $\bC \in \real^{n \times n}$ is the damping matrix, $\bB\in \real^{n \times m}$ is the input matrix, and $\bu(t) \in \real^m$ is the vector of time-dependent control inputs. The dissipative behavior for these mechanical systems is modeled using Rayleigh damping, i.e., the damping matrix $\bC$ is proportional to the mass and stiffness matrix: 
\begin{equation}
        \bC=\alpha\bM + \beta \bK, \quad \quad \alpha,\beta \in \real^+.
\end{equation}
Substituting expressions for the system Lagrangian $L(\bq,\dot{\bq})$ from equation \eqref{eq:Lagrangian} and the nonconservative external forcing $\bf(\bq,\dot{\bq},t)$ from equation \eqref{eq:fn} into the forced Euler-Lagrange equations \eqref{eq:fel} yields the resulting governing equations of motion 
\begin{equation}
    \label{eq:fom}
        \bM \ddot{\bq}(t) + \bC \dot{\bq}(t) + \bK \bq(t) = \bB \bu(t).
\end{equation}
The corresponding output equation is given by
\begin{equation}
           \by(t)=\bE\bq(t),
           \label{eq:output}     
\end{equation}
where $\by(t) \in \real^p$ is the output vector and $\bE \in \real^{p \times n}$ is the output matrix.

\par

Although the nonconservative external forcing $\bf(\bq,\dot{\bq},t)$ from equation \eqref{eq:fn} violates the symplectic structure and, in general, breaks the symmetries of the Lagrangian, the variational approach reveals how the external forcing affects the time evolution of these quantities. This plays a crucial role in developing time integrators for Lagrangian FOMs that track the change in energy or some conserved quantity accurately.

\subsubsection{Nonlinear wave equations} 
\label{sec:non}
The governing PDEs for many space-time continuous physical systems can be derived from a space-time continuous Lagrangian using the Euler-Lagrange equations. We consider nonlinear FOMs obtained via structure-preserving spatial discretization of infinite-dimensional nonlinear Lagrangian PDEs. For illustration purposes, we focus on one-dimensional nonlinear wave equations of the type
\begin{equation}
\frac{\partial ^2 q}{\partial t^2}-\frac{\partial ^2 q}{\partial x^2} + \frac{\text{d} U_{\text{nl}}(q)}{\text{d} q}=0,
\label{eq:wave}
\end{equation}
where $x$ is the spatial variable and $U_{\text{nl}}(q)$ is the nonlinear component of the potential energy. The corresponding space-time continuous Lagrangian is
\begin{equation}
\mathcal{L}(x,q,q_x,q_t)=\frac{1}{2}\left(  \left( \frac{\partial q}{\partial t}\right)^2 - \left( \frac{\partial q}{\partial x}\right)^2\right) - U_{\text{nl}}(q).
\label{eq:unl}
\end{equation}
This specific form of the space-time continuous Lagrangian covers a number of important nonlinear wave equations found in engineering and science applications, see~\cite{strauss1990nonlinear,marsden1998multisymplectic,cheviakov2020invariant}. Space-discretized nonlinear Lagrangian FOMs are usually derived from the nonlinear wave equation \eqref{eq:wave} by symmetric finite differences or pseudo-spectral methods~\cite{li2017spectral}. Those methods discretize the spatial dimension in an accurate and structure-preserving way, so that the resulting space-discretized FOM is a finite-dimensional Lagrangian system. One way to do this is to discretize the Lagrangian density directly and then derive Euler-Lagrange equations for the space-discretized Lagrangian. Direct discretization of the space-time continuous Lagrangian~\eqref{eq:unl} with $n$ equally spaced points leads to 
\begin{equation*}
L(\bq,\dot{\bq})= \frac{1}{2} \sum_{i=1}^n\left( \left( \frac{\partial q_i}{\partial t}\right)^2 -  \left(\sum_{k=1}^n D_{ik}q_k\right)^2 \right) - \sum_{i=1}^n U_{\text{nl}}(q_i),
\end{equation*}
where $\bq=[q_1,q_2, \cdots, q_n]^\top$ with $q_i:=q(t,x_i)$, and the derivative of $q$ with respect to $x$ is approximated by an appropriate differentiation matrix $\bD=(D_{ij})_{i,j=1}^n$, i.e., $\frac{\partial q}{\partial x}(x_i)\approx  \sum_{k=1}^n D_{ik}q_k$. For $\Delta x \to 0$ with $n\Delta x=\ell$, the term $L(\bq,\dot{\bq})\Delta x $ converges to $\int_0^\ell\mathcal L(x,q,q_x,q_t) \ \text{d}x$. The Euler-Lagrange equations for the space-discretized Lagrangian are
\begin{equation}
\ddot{\bq}=\bK\bq +  \frac{\text{d} U_{\text{nl}}(\bq)}{\text{d} \bq}
\label{eq:non_FOM}
\end{equation}
where the linear Lagrangian FOM operator $\bK=\bK^\top$ is always symmetric regardless of the spatial approximation. The nonlinear FOM described by~\eqref{eq:non_FOM} conserves the total energy 
\begin{equation}
E(\bq,\dot{\bq})=\frac{1}{2}\dot{\bq}^\top\dot{\bq} + \frac{1}{2}\bq^\top\bK\bq +  \sum_{i=1}^n U_{\text{nl}}(q_i).
\label{eq:FOM_energy}
\end{equation}
The solution trajectories of~\eqref{eq:non_FOM} also preserve the Lagrangian symplectic form. Additionally, if the system Lagrangian possesses symmetries then the solution trajectories exhibit additional invariants of motion.

\subsection{Structure-preserving time integration of Lagrangian FOMs} 
\label{sec:gni}
Variational integrators provide a systematic way of deriving structure-preserving numerical integrators for Lagrangian systems. We closely follow \cite{marsden2001discrete} to give a brief review of the construction of variational integrators for Lagrangian systems with nonconservative external forcing. These time integrators are based on a discrete version of Lagrange-d'Alembert principle where the basic idea is to first construct discrete approximations of both the action integral and virtual work terms in \eqref{eq:lda} and then use concepts from discrete mechanics to derive variational integrators. 
\par

For a fixed time step $\Delta t$, the discrete trajectory $\{ \bq_{k}\}^{k=K}_{k=0}$ is defined by the  configuration of the Lagrangian system at the discrete time values $\{t_k=k\Delta t\}^{k=K}_{k=0}$. We introduce the discrete Lagrangian function $L_{\text d}(\bq_k,\bq_{k+1})$ along with discrete forcing terms $\bf^+$ and $\bf^-$, which approximate the action integral and virtual work terms between $t_k$ and $t_{k+1}$ in the following sense
\begin{align*}
        L_{\text d}(\bq_k,\bq_{k+1})& \approx \int_{t_k}^{t_{k+1}} L(\bq,\dot{\bq}) \text{d}t, \\
        \bf^+_{\text d}(\bq_k,\bq_{k+1}) \cdot \delta \bq_{k+1} + \bf^-_{\text d}(\bq_k,\bq_{k+1}) \cdot \delta \bq_{k} & \approx \int^{t_{k+1}}_{t_{k}}\bf(\bq,\dot{\bq},t) \cdot \delta \bq  \ \text{d}t.
\end{align*}
Using these discrete approximations of integral terms in \eqref{eq:lda}, we consider a discrete version of the Lagrange-d'Alembert principle that seeks $\{ \bq_{k}\}^{k=K}_{k=0}$ that satisfy
\begin{equation}
        \delta \sum_{k=0}^{K-1} L_{\text d}(\bq_k,\bq_{k+1}) + \sum_{k=0}^{K-1}[\bf_{\text d}^+(\bq_k,\bq_{k+1})\cdot\delta \bq_{k+1} +  \bf_{\text d}^-(\bq_k,\bq_{k+1})\cdot\delta \bq_{k}] = 0,
\end{equation}
which yields the following discrete forced Euler-Lagrange equations
\begin{equation}
        \frac{\partial L_{\text d}(\bq_{k-1},\bq_{k})}{\partial \bq_k} + \frac{\partial L_{\text d}(\bq_k,\bq_{k+1})}{\partial \bq_{k}} + \bf_{\text d}^+(\bq_{k-1},\bq_k) + \bf_{\text d}^-(\bq_{k},\bq_{k+1}) = \bzero   \quad k=1,..., K-1.
\end{equation}

\par

For conservative Lagrangian systems with $\bf=\bzero$ in \eqref{eq:fn}, these variational integrators are automatically symplectic and exhibit bounded energy error for exponentially long times. For Lagrangian systems with nonconservative external forcing, these integrators have been shown to track the change in energy accurately \cite{kane2000variational,sharma2018energy}. In addition to their excellent energy behavior, variational integrators also conserve invariants of the dynamics associated with the symmetries of the FOM Lagrangian via a discrete version of Noether's theorem \cite{marsden2001discrete}.

\par
The conservative Lagrangian FOMs and ROMs for the linear Euler-Bernoulli beam example in Section \ref{sec:beam} and the nonlinear sine-Gordon equation in Section~\ref{sec:sg} are numerically integrated using a variational integrator based on the midpoint rule. For the nonconservative soft-robotic fishtail example in Section \ref{sec:fishtail}, we use the Newmark integrator with $\gamma=0.5$ for all FOMs and ROMs. The Newmark integrator is a variational integrator for Lagrangian systems which preserves the symplectic structure and exhibits bounded energy error for both FOM and ROM simulations. Details about the geometric properties of this variational integrator can be found in \cite{kane2000variational}.
\subsection{Projection-based intrusive structure-preserving model reduction of Lagrangian FOMs} \label{sec:intrusive}
The proposed data-driven Lagrangian-preserving model reduction method is strongly motivated by the intrusive projection-based structure-preserving model reduction of Lagrangian FOMs~\cite{lall2003structure,carlberg2015preserving} which we briefly review in this section.
Galerkin projection-based methods first replace the FOM configuration space $\Q$ by a reduced-order configuration space $\Q_r$, define an intrusive reduced Lagrangian $\tL$ on the reduced state space $\T\Q_r$, and then derive governing equations for the Lagrangian ROM using a set of reduced-dimension coordinates. The state is approximated via $\bq=\bV_r\tbq$ where $\bV_r = [\bv_1, \cdots, \bv_r] \in \real^{n \times r}$ is an orthogonal basis matrix whose columns span an $r$-dimensional reduced-order configuration space $\Q_r$. Proper orthogonal decomposition (POD) \cite{holmes2012turbulence} computes the basis matrix $\bV_r$ from the snapshot matrix $\bQ=[\bq_1, \cdots, \bq_K]$ via singular value decomposition (SVD). The intrusive reduced Lagrangian $\tL_r$ is defined as
\begin{equation}
        \tL_r(\tbq,\tbqd):=L(\bV_r\tbq,\bV_r\tbqd)=T(\bV_r\tbqd) - U(\bV_r\tbq).
        \label{eq:Ltilde_gen}
\end{equation}
The intrusive reduced nonconservative forcing $\tbf$ is
\begin{equation}
        \tbf(\tbq,\tbqd,t):=\bV_r^\top\bf(\bV_r\tbq,\bV_r\tbqd,t).
        \label{eq:ftilde_gen}
\end{equation}
Following the variational derivation for forced Lagrangian systems outlined in Section \ref{sec:FOM}, the forced Euler-Lagrange equations in $r$ dimensions are
\begin{equation}
        \frac{\partial \tL_r(\tbq,\tbqd)}{\partial \tbq} -\frac{\text d}{\text dt} \left( \frac{\partial \tL_r(\tbq,\tbqd)}{\partial \tbqd} \right) + \tbf(\tbq,\tbqd,t)  = \bzero.
    \label{eq:ELr_gen}
\end{equation}
Since the reduced equations are derived from a Lagrangian system, the resulting reduced equations of motion preserve the underlying Lagrangian structure. In this work, we are mainly interested in two types of Lagrangian systems.
\begin{enumerate}
\item {\textbf{Mechanical systems with external nonconservative forcing.}}
For the Lagrangian FOMs of the form \eqref{eq:fom}, the intrusive reduced Lagrangian $\tL_r$ is
\begin{equation}
        \tL_r(\tbq,\tbqd)= \frac{1}{2}\tbqd^\top(\bV_r^\top\bM\bV_r)\tbqd - \frac{1}{2}\tbq^\top (\bV_r^\top\bK\bV_r) \tbq,
        \label{eq:Ltilde_simple}
\end{equation}
and the intrusive reduced forcing $\tbf$ is
\begin{equation}
        \tbf(\tbq,\tbqd,t)=\bV_r^\top\bf(\bV_r\tbq,\bV_r\tbqd, t)=-(\bV_r^\top\bC\bV_r)\tbqd + (\bV_r^\top\bB) \bu(t),
        \label{eq:frint_simple}
\end{equation}
where the last equality results from the specific form of the external forcing, see \eqref{eq:fn}. 
Substituting the intrusive reduced Lagrangian $\tL_r$ from equation \eqref{eq:Ltilde_simple} and the intrusive reduced forcing $\tbf$ from \eqref{eq:frint_simple} into \eqref{eq:ELr_gen} results in the reduced equations of motion 
\begin{equation}
                \tbM \tbqdd(t) + \tbC \tbqd(t) + \tbK \tbq(t) = \tbB \bu(t), 
                \label{eq:rom_int}
\end{equation}
where $\tbM:=\bV_r^\top\bM\bV_r$, $\tbC:=\bV_r^\top\bC\bV_r$, $\tbK:=\bV_r^\top\bK\bV_r$, and $\tbB:=\bV_r^\top\bB$ are the intrusive reduced operators. It is clear from these expressions for the intrusive reduced operators $\tbM, \tbC,$ and $\tbK$ that they are all symmetric and positive-definite matrices. The corresponding output equation is 
\begin{equation}
        \by(t)=\tbE \tbq (t),
\end{equation}
where $\tbE=\bE\bV_r$ is the reduced output operator. 
\item {\textbf{Nonlinear wave equations.}} For the nonlinear Lagrangian FOMs of the form~\eqref{eq:non_FOM}, the intrusive reduced Lagrangian $\tL$ is
\begin{equation}
        \tL_r(\tbq,\tbqd)= \frac{1}{2}\tbqd^\top\tbqd- \frac{1}{2}\tbq^\top (\bV_r^\top\bK\bV_r) \tbq - U_{\text{nl}}(\bV_r\tbq).
        \label{eq:Ltilde_non}
\end{equation}
Substituting the intrusive reduced Lagrangian $\tL_r$ from equation \eqref{eq:Ltilde_non} into \eqref{eq:ELr_gen} gives us the reduced equations of motion 
\begin{equation}
    \hbqdd(t) = \tbK\hbq(t) + \frac{ \text{d} U_{\text{nl}} (\bV_r\tbq)}{\text{d} \tbq},
    \label{eq:rom_int_non}
\end{equation}
where $\tbK:=\bV_r^\top\bK\bV_r$ is the intrusive reduced operator that retains the symmetric property of the structure-preserving spatial discretizations.
\end{enumerate}

It was shown in \cite{lall2003structure} that Galerkin projection carried out on the forced Euler-Lagrange equations \eqref{eq:fom} lead to the same reduced equations, and hence, preserve the underlying Lagrangian structure. 
\section{Lagrangian operator inference} \label{sec:lopinf}
In this section, we propose \textsf{L-OpInf}, a Lagrangian operator inference framework to learn reduced-order operators of mechanical systems of the form \eqref{eq:fom} and nonlinear wave equations of the form~\eqref{eq:non_FOM}. In Section~\ref{sec:motivation} we motivate the need for \textsf{L-OpInf} by demonstrating on a model of an Euler-Bernoulli beam how the standard operator inference from \cite{peherstorfer2016data} does not preserve the underlying Lagrangian structure leading to unbounded energy growth. Based on the observations from this motivating example, we present \textsf{L-OpInf} in Section~\ref{sec:method}. We summarize the computational procedure of \textsf{L-OpInf} and some practical considerations in Section~\ref{sec:comp}.
\begin{figure}[h]
\centering
\scalebox{0.85}{


  \begin{tikzpicture}
    \draw[thick,->] (0,-2.1) -- (0.6,-2.1) node[left=-0.5]{$x$};
  \draw[thick,->] (0,-2.1) -- (0,-1.5) node[above=0]{$w(x,t)$};
    \point{support1}{-0.2}{-0.1};
    \point{support2}{15.2}{-0.1};
    \point{begin}{0}{0};
    \point{middle}{7.5}{2.2};
    \point{middle2}{7.5}{-1.5};
    \point{end}{15}{0};
    \support{2}{begin}[0];
    \support{2}{end}[0];

     \draw[->] (middle) to [out=180,in=100] (support1);
    \draw[->] (middle) to [out=0,in=60] (support2);
    \draw[->] (middle2) -- (3.5,.28);
    \draw[->] (middle2) -- (8.5,-.6);
    \draw[->] (middle2) -- (9.4,-1.1);
    \draw[<->] (0,-2.5) -- (15,-2.5);
    
    \node[text width=2cm] at (8.5,-2.8) {\Large $\ell$};
    \node[text width=5.3cm,fill=white] at (8.5,-1.8) {Transverse vibrations $w(x,t)$};
    \node[text width=5.3cm,fill=white] at (8.5,2.2) {Simply-supported at both ends};
        \draw[scale=1, domain=0:15, smooth, variable=\x, black,line width=3mm] plot ({\x}, {4*(\x/15) - 4*(\x/15)^2});
        \draw[scale=1, domain=0:15, smooth, variable=\x, black,dashed,line width=0.5mm] plot ({\x}, {-4*(\x/15) + 4*(\x/15)^2});
    \draw[scale=1, domain=0:15, smooth, variable=\x, black,dashed,line width=0.5mm] plot ({\x}, {2*(\x/15) - 2*(\x/15)^2});
     \draw[scale=1, domain=0:15, smooth, variable=\x, black,dashed,line width=0.5mm] plot ({\x}, {-2*(\x/15) + 2*(\x/15)^2});
     
  \end{tikzpicture}
}
 \caption{Euler-Bernoulli beam: A schematic showing transverse vibrations in response to a nonzero initial condition. The beam is simply supported at both ends which allows for rotation but not for vertical displacement.}
 \label{fig:LB_schematic}
\end{figure}
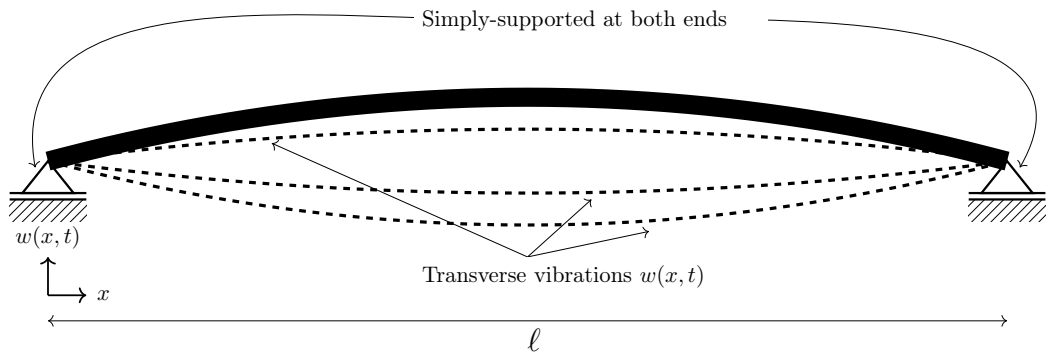
%
%
%
%
%
\subsection{Motivational example}
\label{sec:motivation}
Consider the numerical example of transverse vibrations of an Euler-Bernoulli beam, which we revisit with much more detail in Section~\ref{sec:beam}. A schematic of a simply-supported Euler-Bernoulli beam vibrating in response to some nonzero initial condition is shown in Figure \ref{fig:LB_schematic}.  We apply standard operator inference to the second-order system arising from the beam FOM model to demonstrate how violating the underlying Lagrangian structure leads to unstable ROMs.

\par

The governing PDE is discretized in space using finite elements  which leads to a discretized state vector $\bq \in \real^{400}$ which contains the transverse deflection $w(x,t)$ (see Figure \ref{fig:LB_schematic}) at different grid points along its length. The resulting FOM is integrated using a variational integrator based on the midpoint rule for $T=0.18$~\si{\s} using a fixed time step of $\Delta t=10^{-5}$~\si{\s}. Details about the beam FOM implementation can be found in \ref{sec:beam_fom}. Based on the snapshot data $\bQ=[\bq_1, \cdots,\bq_K]$ from $t=0$~\si{\s} to $t=0.03$~\si{\s}, we compute the POD basis $\bV_{2r}$ (see equation \eqref{eq:beam_basis} later). We then project the FOM snapshot data onto the basis matrix $\bV_{2r}$ to obtain the reduced snapshot data $\hbQ \in \real^{2r \times K}$. We also obtain the  reduced second-order time-derivative data $\hbQdd \in \real^{2r \times K}$ using an eighth-order finite difference approximation, see equation \eqref{eq:second} below. 

\par

\begin{figure}[h]
\captionsetup[subfigure]{oneside,margin={1.8cm,0 cm}}
\begin{subfigure}{.38\textwidth}
       \setlength\fheight{6 cm}
        \setlength\fwidth{\textwidth}
%
%
\begin{tikzpicture}

\begin{axis}[%
width=0.951\fheight,
height=0.636\fheight,
at={(0\fheight,0\fheight)},
scale only axis,
xmin=0,
xmax=40,
xlabel style={font=\color{white!15!black}},
xlabel={\small Reduced dimension $2r$},
ymode=log,
ymin=0.001,
ymax=0.1,
yminorticks=true,
ylabel style={font=\color{white!15!black}},
ylabel={\small Relative state error},
axis background/.style={fill=white},
xmajorgrids,
ymajorgrids,
legend style={at={(0.0,1.3)}, anchor=south west, legend cell align=left, align=left, draw=white!15!black,font=\small}
]
\addplot [color=red, dotted, line width=3.0pt, mark size=4.0pt, mark=triangle, mark options={solid, red}]
  table[row sep=crcr]{%
4	0.0596072966559552\\
8	0.024371353252337\\
12	0.0143783294157558\\
16	0.0100455850959177\\
20	0.00779259510965184\\
24	0.00644391860250537\\
28	0.00557756355977969\\
32	0.00501711927174082\\
36	0.00472650076150242\\
40	0.00458018343758944\\
};
\end{axis}

\begin{axis}[%
width=1.227\fheight,
height=0.658\fheight,
at={(-0.16\fheight,-0.072\fheight)},
scale only axis,
xmin=0,
xmax=1,
ymin=0,
ymax=1,
axis line style={draw=none},
ticks=none,
axis x line*=bottom,
axis y line*=left
]
\end{axis}
\end{tikzpicture}%
\caption{State error (training data)}
\label{fig:motivation_state}
    \end{subfigure}
    \hspace{1.4cm}
    \begin{subfigure}{.38\textwidth}
           \setlength\fheight{6 cm}
           \setlength\fwidth{\textwidth}
\input{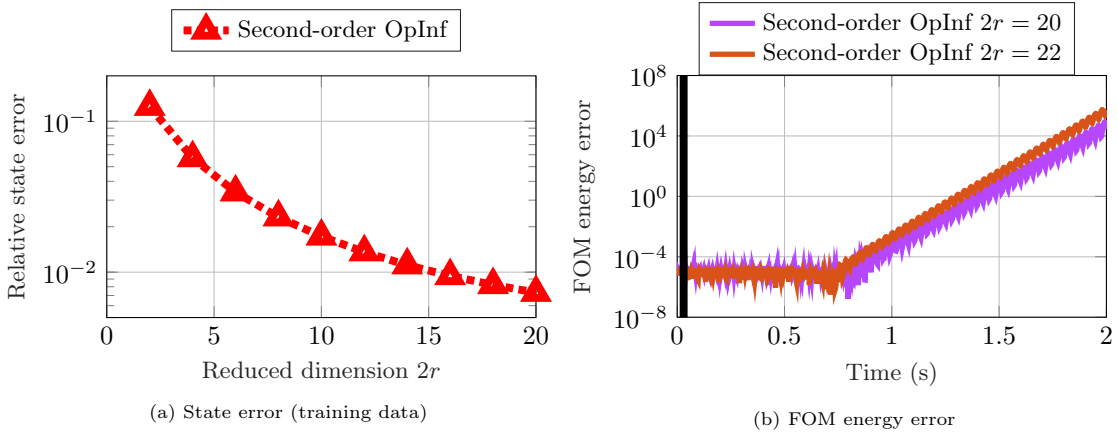}
\caption{FOM energy error}
\label{fig:motivation_energy}
    \end{subfigure}
\caption{Euler-Bernoulli beam: Even though plot (a) shows low state approximation error for second-order operator inference in the training data regime, the corresponding FOM energy error \eqref{eq:fom_energy_error} behavior in plot (b) reveals that the standard second-order operator inference violates the underlying Lagrangian structure which leads to unbounded growth of the energy error outside the training data. The vertical black line in plot (b) indicates end of training time interval $[0, 0.03]$ s, which is a small fraction of the full simulation time shown in (b). While this training data is enough to learn the model accurately (see (a)), it still does not produce a conservative Lagrangian system.}
 \label{fig:motivation}
\end{figure}

Based on the conservative nature of the considered FOM, we postulate a model form $\hbqdd + \hK \hbq = \bzero$ for learning the second-order ROM. Thus, for the Euler-Bernoulli beam equation, the second-order operator inference problem solves the following least-squares problem for the reduced stiffness matrix $\hbK \in \real^{2r \times 2r}$:
    \begin{equation}\label{eq:std_opinf}
     \min_{\hbK}   || \hbQdd + \hbK \hbQ ||_F.
    \end{equation}
Figure \ref{fig:motivation_state} shows the relative state error over the training data for learned second-order ROMs of different dimensions. We see that the state approximation error decreases monotonically from $2r=2$ to $2r=40$, i.e., the learned ROMs  approximate the state solution accurately in the training data regime. However, the reduced operators do not preserve the underlying Lagrangian structure. The energy error plots in Figure \ref{fig:motivation_energy} show that the FOM energy grows unbounded in the testing data for every standard second-order OpInf ROM from $2r=20$ to $2r=40$. These plots with unbounded energy error growth indicate unphysical solutions that eventually blow up when making predictions outside the training data, demonstrating that neglecting the physical structure can have detrimental effects on the predictive capability of data-driven ROMs.

We note that the FOM energy error growth in Figure \ref{fig:motivation_energy} is not due to a lack of training data. We have trained the models with a larger training interval $[0, 0.06]$~\si{\s}, and the learned ROMs still exhibit FOM energy blowup (not shown here) because standard second-order operator inference neglects the underlying geometric structure. A regularized operator inference method has been presented in~\cite{mcquarrie2021data} to promote stability of long-time integration. It should be noted, however, that this regularized operator inference method has no built-in constraints on the ROM operators and therefore, the resulting ROM would still violate the underlying geometric structure and does not produce a Lagrangian system model. Hence, the regularization techniques presented in~\cite{mcquarrie2021data} can not resolve the unbounded energy error in Figure~\ref{fig:motivation_energy}.

\subsection{Lagrangian operator inference}
\label{sec:method}

In this work, our goal is to learn a Lagrangian ROM for FOMs of the form \eqref{eq:fom} and \eqref{eq:non_FOM}, where the trajectories and outputs can be computed, but where the FOM operators are unavailable. We consider the situation that we have knowledge about the governing equations at the PDE level and we have simulated data from a convergent FOM time integrator. The proposed physics-preserving inference-based approach learns a Lagrangian ROM from data of a large-scale mechanical system with nonconservative external forcing or a nonlinear wave equation so that the learned ROM:
\begin{enumerate}
    \item is a Lagrangian system;
    \item retains the physical interpretation of state variables;
    \item preserves the Riemannian metric for large-scale mechanical systems;
    \item respects the symmetric property of system matrices.
\end{enumerate}

\subsubsection{Lagrangian operator inference for mechanical systems with nonconservative forcing}
\label{sec:method_simple}
Next, we introduce a novel \textsf{L-OpInf} framework for Lagrangian FOMs described in~\ref{sec:simple}. Given known control inputs $\bu(t)$ and initial conditions $(\bq(0),\bqd(0))$, let $\bq_1, \cdots, \bq_K$ be the solutions to the Lagrangian FOM \eqref{eq:fom} at $t_1, \cdots, t_K$ computed using a variational integrator. Let $\by_1, \cdots, \by_K$ be the corresponding outputs at those time instances. We collect this data in the snapshot matrices 
\begin{equation}
    \bQ=[\bq_1, \cdots, \bq_K] \in \real^{n \times K}, \quad \quad \bY=[\by_1, \cdots, \by_K] \in \real^{p \times K}.
    \label{eq:snapshot}
\end{equation}
We also define the input snapshot matrix 
\begin{equation}
    \bU=[\bu(t_1),\cdots,\bu(t_K)] \in \real^{m \times K}.
      \label{eq:snapshot_input}
\end{equation}
\par

Given these snapshot matrices of Lagrangian FOM state and output trajectories along with the input snapshot matrix, our goal is to learn a Lagrangian ROM directly from the data. To learn the reduced operators, we first prepare the data, i.e., we project FOM trajectories onto low-dimensional subspaces of the high-dimensional state spaces and compute the reduced time-derivative data. We then propose to fit operators to the projected trajectories in a structure-preserving way.

\par
As a projection bases, we use the POD basis, which can be computed via the SVD of the snapshot data matrix $\bQ$. We compute
\begin{equation}
    \bQ=\bV\bXi \bW^\top
    \label{eq:svd}
\end{equation}
where $\bV \in \real^{n \times n}$, $\bXi \in \real^{n \times n}$, and $\bW \in \real^{K \times n}$, and we assume that the singular values $\xi_1 \geq \xi_2 \geq \cdots $ in $\bXi$ are ordered from largest to smallest.
Then POD basis matrix is $\bV_r \in \real^{n \times r}$, the leading $r$ columns of $\bV$. The projected state $\hbq(t)\in \real^r$ of the FOM state $\bq(t) \in \real^n$ is $\hbq(t)=\bV_r^\top\bq(t)$. Using this relation, we obtain reduced snapshot data via the projections onto the POD basis matrix $\bV_r$ as
\begin{equation}
    \hbQ=\bV_r^\top\bQ=[\hbq_1, \cdots, \hbq_K] \in \real^{r \times K}.
    \label{eq:qhat}
\end{equation}
We also compute $\hbqd$ and $\hbqdd$ from the reduced trajectory data $\hbq$ using a finite difference scheme, e.g., via the eighth-order central finite difference scheme 
\begin{align}
  \label{eq:first}
    \dot{\hq}_k & \approx \frac{4(\hq_{k+1}-\hq_{k-1})}{5\Delta t} -\frac{(\hq_{k+2}-\hq_{k-2})}{5\Delta t} + \frac{4(\hq_{k+3}-\hq_{k-3})}{105\Delta t} - \frac{(\hq_{k+4}-\hq_{k-4})}{280\Delta t}, \\
    \ddot{\hq}_k & \approx -\frac{205\hq_{k}}{72\Delta t^2} + \frac{8(\hq_{k+1}+\hq_{k-1})}{5\Delta t^2} -\frac{(\hq_{k+2}+\hq_{k-2})}{5\Delta t^2} + \frac{8(\hq_{k+3}+\hq_{k-3})}{315\Delta t^2} - \frac{(\hq_{k+4}+\hq_{k-4})}{560\Delta t^2}.
    \label{eq:second}
\end{align}
We derive ROMs for the Euler-Bernoulli beam example using reduced time-derivative data obtained via finite difference schemes of increasing accuracy and observe that the learned ROM accuracy improves marginally with higher accuracy in second-order time derivatives. Based on the results from this study (not shown here), we use the eighth-order central finite difference scheme to compute the reduced time-derivative data. These time-derivative approximations are used to build the snapshot matrices of the reduced first-order and second-order time-derivative data
\begin{equation}
    \hbQd=[\hbqd_1, \cdots, \hbqd_K] \in \real^{r \times K}, \quad \quad \hbQdd=[\hbqdd_1, \cdots, \hbqdd_K] \in \real^{r \times K}.
    \label{eq:qhatdot}
\end{equation}
We postulate the form of the reduced Lagrangian (motivated by \eqref{eq:Ltilde_simple}) as
  \begin{equation}
        \hL_r(\hbq,\hbqd)=\frac{1}{2}\hbqd^\top\hbqd - \frac{1}{2}\hbq^\top \hbK \hbq,
        \label{eq:Lhat}
    \end{equation}
where $\hbK \in \real^{r \times r}$ is the symmetric positive-definite reduced stiffness matrix that is learned from data. We postulate the form of the nonintrusive reduced forcing $\hbf$ based on \eqref{eq:fn} and \eqref{eq:frint_simple} as
    \begin{equation}
        \hbf(\hbq,\hbqd,t)=- \hbC\hbqd + \hbB \bu(t),
    \end{equation}
where $\hbC \in \real^{r \times r}$ is the symmetric positive-definite reduced damping matrix that is learned from data and $\hbB \in \real^{r \times m}$ is the reduced input matrix that is learned from data. Based on the assumed model form for $\hL_r(\hbq,\hbqd)$ and nonintrusive reduced forcing $\hbf$, we derive the governing equations for the reduced system via the forced Euler-Lagrange equations \eqref{eq:fel} and obtain
\begin{equation}
    \hbqdd(t) + \hbC \hbqd(t) + \hbK \hbq(t) = \hbB \bu(t),
    \label{eq:rom_learn}
\end{equation}
along with the reduced output equation 
\begin{equation}
    \by(t)=\hbE\hbq(t), 
      \label{eq:rom_output}
\end{equation}
where $\hbE \in \real^{p \times r}$ is the reduced output operator that is learned from data.

\par

The Lagrangian ROM form leads us to propose the following optimization problems  to compute $\hbC \in \real^{r \times r},\hbK \in \real^{r \times r}$, and $\hbB \in \real^{r \times m}$:
\begin{equation}
\min _{\substack{\hbK=\hbK^\top \succ 0, \hbC=\hbC^\top \succ 0,\\\hbB}}
\lVert \hbQdd + \hbC\hbQd + \hbK\hbQ- \hbB \bU \rVert_{F}.
\label{eq:lopinf}
\end{equation}
The symmetric positive-definite constraints on $\hbK$ and $\hbC$ ensure that the learned reduced operators respect the underlying geometric structure of FOMs, and hence, the ROMs learned via \textsf{L-OpInf} are Lagrangian systems. The constrained optimization problem~\eqref{eq:lopinf} has a unique minimizer if and only if $\hbQ$, $\hbQd$, and $\mathbf{U}$ have full column rank; see~\cite{allwright1988positive} for more details about linear least-squares problems with symmetric positive definite constraints.

For the inference of the reduced output operator $\hbE \in \real^{p \times r}$ we solve the least-squares problem
\begin{equation}
\min _{\hbE}
\lVert \bY - \hbE\hbQ \rVert_{F}.
\label{eq:lopinf_output}
\end{equation}

Although the \textsf{L-OpInf} framework for mechanical systems is strongly motivated by analogies to the intrusive projection-based approach described in Section~\ref{sec:intrusive}, the nonintrusive ROM model form \eqref{eq:rom_learn} differs from \eqref{eq:rom_int} in the following sense. The proposed nonintrusive Lagrangian \eqref{eq:Lhat} can be interpreted as a special case of the more general nonintrusive Lagrangian $\hL_r(\hbq,\hbqd)=\frac{1}{2}\hbqd^\top\hbM\hbqd - \frac{1}{2}\hbq^\top \hbK \hbq$ with $\hbM=\mathbf{I}_r$. While this specific choice for the reduced mass matrix restricts the search of reduced Lagrangian operators to a subset of the full solution space, it plays a key role in making the constrained optimization problem for inferring the reduced Lagrangian operators tractable. We emphasize that even though we search for models with reduced mass matrix $\hbM=\mathbf{I}_r$, the nonintrusive ROM model form~\eqref{eq:rom_learn} is not equivalent to premultiplying~\eqref{eq:rom_int} by $\hbM^{-1}$. The specific choice of $\hbM=\mathbf{I}_r$ ensures that the ROM model form retains all the relevant geometric properties and the Lagrangian structure whereas premultiplying~\eqref{eq:rom_int} by $\hbM^{-1}$ loses the symmetric property of the reduced stiffness matrix $\hbK$ and reduced damping matrix $\hbC$ and therefore violates the Lagrangian structure.
\subsubsection{Lagrangian Operator Inference for Nonlinear Wave Equations}
\label{sec:method_non}
We use knowledge about the nonlinear potential energy $U_{\text{nl}}$ at the PDE level to develop a gray-box approach to derive nonintrusive Lagrangian ROMs for spatial discretizations of nonlinear wave equations. For the nonlinear Lagrangian FOMs discussed in Section~\ref{sec:non}, the nonlinear potential energy term $U_{\text{nl}}(q)$ in~\eqref{eq:unl} is assumed to be given explicitly, whereas the quadratic terms in~\eqref{eq:unl} and details about their structure-preserving spatial discretization are unavailable. We define the nonlinear forcing 
\begin{equation}
  \bf_{\text{nl}}(\bq)=\left[ \frac{\text{d} U_{\text{nl}}}{\text{d}  q}(q_1), \cdots,  \frac{\text{d}  U_{\text{nl}}}{\text{d}  q}(q_n)  \right]^\top \in \real^n.
\label{eq:f_non}
\end{equation}
We build the nonlinear forcing snapshot data matrix 
\begin{equation}
\bF_{\text{nl}}=\left[ \bf_{\text{nl}}(\bq_1), \cdots, \bf_{\text{nl}}(\bq_K) \right] \in \real^{n \times K},
\label{eq:snapshot_non}
\end{equation}
where we compute $\bf_{\text{nl}}$ in~\eqref{eq:f_non} at different time instances using the FOM snapshot data. We compute the POD basis $\bV_r \in \real^{n \times r}$ via the SVD of the augmented snapshot data matrix $\bQ_{\text{aug}}=[\bQ, \bF_{\text{nl}}]\in \real^{n \times 2K}$. We obtain projections of the snapshot data $\bQ$ and $\bF_{\text{nl}}$ as
\begin{equation}
\hbQ=\bV_r^\top\bQ\in \real^{r\times K}, \quad \quad \hbF_{\text{nl}}=\bV_r^\top\bF_{\text{nl}}\in \real^{r\times K}.
\label{eq:snapshot_proj}
\end{equation}
For the nonlinear Lagrangian FOMs discussed in Section~\ref{sec:non}, we postulate the form of the reduced Lagrangian 
    \begin{equation}
        \hL_r(\hbq,\hbqd)=\frac{1}{2}\hbqd^\top\hbqd - \frac{1}{2}\hbq^\top \hbK \hbq  - \hU_{\text{nl}}(\hbq),
   \label{eq:Lhat_non}
    \end{equation}
where $\hbK \in \real^{r \times r}$ is the symmetric reduced operator that is learned from data and $\hU_{\text{nl}}(\hbq):=U_{\text{nl}}(\bV_r\hbq)$ is the reduced nonlinear potential energy. Based on the assumed model form for $\hL_r(\hbq,\hbqd)$ in~\eqref{eq:Lhat_non}, we derive the governing equations for the nonlinear reduced system via the Euler-Lagrange equations \eqref{eq:el} and obtain
\begin{equation}
    \hbqdd(t) = \hbK\hbq(t) +  \frac{ \text{d} \hU_{\text{nl}} (\hbq)}{\text{d} \hbq}.
    \label{eq:rom_learn_non}
\end{equation}
We solve the following constrained optimization problem 
\begin{equation}
\min _{\substack{\hbK=\hbK^\top}}
\lVert \hbQdd - \hbF_{\text{nl}} -\hbK\hbQ \rVert_{F}.
\label{eq:lopinf_non}
\end{equation}
to infer the symmetric ROM operator $\hbK$. We impose a symmetry constraint on $\hbK$ to ensure that the learned ROM operator retains the symmetric property of the linear FOM operator introduced during the structure-preserving spatial discretization, see Section~\ref{sec:non}. The symmetric linear least-squares problem~\eqref{eq:lopinf_non} has a unique solution if and only if the reduced snapshot data matrix $\hbQ$ has full column rank.
%
\begin{remark}
The space-discretized FOM for nonlinear wave equations \eqref{eq:non_FOM} can also be written in the Hamiltonian form where the governing equations are a set of $2n$ coupled first-order ODEs. The Hamiltonian operator inference method~\cite{sharma2022hamiltonian} for learning Hamiltonian ROMs requires both trajectory and momentum data whereas the proposed Lagrangian operator inference approach has the advantage that it learns Lagrangian ROMs purely from trajectory data.  
\end{remark}

%
\subsection{Computational procedure}
\label{sec:comp}
Algorithm \ref{alg:lopinf_simple} and Algorithm \ref{alg:lopinf_non} summarize \textsf{L-OpInf} for mechanical systems with nonconservative external forcing and \textsf{L-OpInf} for nonlinear wave equations as discussed in Section \ref{sec:method_simple} and \ref{sec:method_non}, respectively. The constrained optimization problems in~\eqref{eq:lopinf} and~\eqref{eq:lopinf_non} are solved using the CVX optimization package \cite{cvx}, a MATLAB-based software for constrained optimization problems which allows constraints and objectives to be specified using standard MATLAB expression syntax. The CVX package supports four solvers with different capabilities and various levels of performance. The constrained optimization problems arising in this work are solved using SDPT3~\cite{toh2012implementation}. The optimization algorithm implemented in SDPT3 is a primal-dual interior point algorithm that uses the path-following paradigm.
\subsubsection{Computational cost}
The computational cost of Algorithm \ref{alg:lopinf_simple} and Algorithm \ref{alg:lopinf_non} is typically dominated by the final step which requires solution of constrained linear least-squares problem. The size of this constrained operator inference problem mainly depends on the reduced dimension $r$ and the number of snapshots $K$. Figure~\ref{fig:cvx} compares the MATLAB wall clock time (averaged over 20 runs) of solving the constrained optimization problem~\eqref{eq:lopinf} using MATLAB 2020b on a quad-core Intel i7 processor with 2.3 GHz and 32 GB RAM. We observe that the computational cost increases exponentially with increase in the reduced dimension $2r$. Moreover, for fixed reduced dimension $2r$, the computational cost also increases with an increase in the number of training snapshots, $K$.

\begin{figure}[h]
\centering
  \setlength\fheight{9 cm}
        \setlength\fwidth{\textwidth}
%
%
\begin{tikzpicture}

\begin{axis}[%
width=0.976\fheight,
height=0.59\fheight,
at={(0\fheight,0\fheight)},
scale only axis,
xmin=0,
xmax=20,
xlabel style={font=\color{white!15!black}},
xlabel={\small Reduced Dimension $2r$},
ymin=0,
ymax=15,
ylabel style={font=\color{white!15!black}},
ylabel={\small MATLAB wall clock time [\si{\s}]},
axis background/.style={fill=white},
xmajorgrids,
ymajorgrids,
legend style={at={(0.043,0.505)}, anchor=south west, legend cell align=left, align=left, draw=white!15!black,font=\small}
]
\addplot [color=red, dotted, line width=2.0pt, mark size=4.0pt, mark=o, mark options={solid, red}]
  table[row sep=crcr]{%
   2.0000    0.3170\\
    4.0000    0.3860\\
    6.0000    0.5503\\
    8.0000    0.9575\\
   10.0000    1.5813\\
   12.0000    1.8716\\
   14.0000    2.7019\\
   16.0000    3.7594\\
   18.0000    4.2467\\
   20.0000    5.2112\\
};
\addlegendentry{Training snapshots $K=3,000$}


\addplot [color=black, dashed, line width=2.0pt, mark size=4.0pt, mark=triangle, mark options={solid, black}]
  table[row sep=crcr]{%
    2.0000    0.3158\\
    4.0000    0.5103\\
    6.0000    0.9597\\
    8.0000    1.5997\\
   10.0000    2.4459\\
   12.0000    3.3271\\
   14.0000    4.3692\\
   16.0000    6.1348\\
   18.0000    7.5483\\
   20.0000    9.3486\\
};
\addlegendentry{Training snapshots $K=6,000$}

\addplot [color=blue, dashdotted, line width=2.0pt, mark size=4pt, mark=diamond, mark options={solid, blue}]
  table[row sep=crcr]{%
  2.0000    0.3265\\
    4.0000    0.6677\\
    6.0000    1.2730\\
    8.0000    2.1025\\
   10.0000    3.2839\\
   12.0000    4.6342\\
   14.0000    6.2025\\
   16.0000    8.8516\\
   18.0000   10.7667\\
   20.0000   14.4189\\
};
\addlegendentry{Training snapshots $K=10,000$}

\end{axis}

\begin{axis}[%
width=1.259\fheight,
height=0.743\fheight,
at={(-0.164\fheight,-0.097\fheight)},
scale only axis,
xmin=0,
xmax=1,
ymin=0,
ymax=1,
axis line style={draw=none},
ticks=none,
axis x line*=bottom,
axis y line*=left
]
\end{axis}
\end{tikzpicture}%
 \caption{Euler-Bernoulli beam: Computational cost of solving the constrained optimization problem~\eqref{eq:lopinf} for ROMs of reduced dimension $2r$ and different sizes of the training data, $K$.}
 \label{fig:cvx}
\end{figure}

\subsubsection{Practical considerations for solving the constrained optimization problems}
For a nonlinear wave equation with FOM equations described by \eqref{eq:non_FOM} the proposed inference method requires solving the symmetric linear least-squares problem \eqref{eq:lopinf_non}. Since CVX supports declaration of symmetric matrix variables, we declare $\hbK$ as a symmetric matrix and then solve for $\hbK$. For a mechanical system with FOM equations described by \eqref{eq:fom} the proposed inference method requires solving the constrained linear least-squares problem \eqref{eq:lopinf}. To ensure that the reduced operators satisfy the hard constraints required for preserving the Lagrangian structure we utilize the semidefinite programming (SDP) mode provided by CVX. Similar to the nonlinear wave equation problem, we first declare both $\hbK$ and $\hbC$ as symmetric matrices. In floating point arithmetic, the difference between a positive-definite and positive-semidefinite matrices become blurred due to rounding issues, so we impose the following constraints in SDP mode 
    \begin{equation}
        \hbK - \epsilon_{\texttt{tol}} \cdot \mathbf{I}_r \succeq 0, \quad \quad \hbC - \epsilon_{\texttt{tol}} \cdot \mathbf{I}_r \succeq 0,
        \label{eq:cvx}
    \end{equation}
to ensure that the smallest allowed eigenvalue for both $\hbK$ and $\hbC$ is $\epsilon_{\texttt{tol}}$.

We make two additional remarks on practical aspects of solving the constrained optimization problem in \textsf{L-OpInf}. First, the CVX solver performance may depend on scaling of the input data. Input data that spans many orders of magnitude is often seen in multiphysics models where the FOM state vector contains different physical quantities with different scale of magnitudes~\cite{SKHW2020_learning_ROMs_combustor,JMcQK,qian2022reduced}. Thus, it is important to exploit any available information about the FOM system to improve the numerical scaling of the snapshot data $\bQ$. Second, in practice, CVX is not guaranteed to yield reduced operators that satisfy the positive-definite constraint \eqref{eq:cvx}, and in some cases, CVX gives solution status as `\texttt{Inaccurate/Solved}'. This indicates that the CVX solver was unable to find a solution within the default numerical tolerance. However, the CVX solver still returns inferred reduced operators that satisfy a relaxed tolerance value instead, and these reduced operators may still be useful as ROMs. However, extra care is needed in these situations, e.g., monitoring the eigenvalues of the inferred reduced operators $\hbK$ and $\hbC$ can help in detecting if the hard constraints are satisfied or not. In the numerical examples below we always ensure that the hard constraints are satisfied by testing the validity of the inferred ROM operators before using them for simulating the ROMs.
\subsubsection{Alternate constrained optimization problems and their challenges}
For a more general reduced Lagrangian model form $\hL_r(\hbq,\hbqd)=\frac{1}{2}\hbqd^\top\hbM\hbqd - \frac{1}{2}\hbq^\top \hbK \hbq$, the Lagrangian ROM dynamics are    
\begin{equation}\label{eq:rom_general_form}
        \hbM \hbqdd + \hbC\hbqd + \hbK \hbq = \hbB \bu,
\end{equation}
where $\hbM, \hbC,$ and $\hbK$ are symmetric and positive-definite matrices. The corresponding constrained operator inference problem 
\begin{equation}
\min _{\substack{\hbM=\hbM^\top \succ 0, \hbK=\hbK^\top \succ 0,\\ \hbC=\hbC^\top \succ 0, \bB}} \lVert \hbM \hbQdd + \hbC\hbQd + \hbK\hbQ -\bB\mathbf{U} \rVert_F,
\label{eq:constrained_general}
\end{equation}
is, however, challenging to solve. In all of the numerical experiments, we observe that the CVX solver fails to find a reduced mass matrix $\hbM$ that satisfies the positive definite constraint. For low-dimensional ROMs with $r<8$, the CVX solver often gives solution status as `\texttt{Infeasible}'  or `\texttt{Inaccurate/Infeasible}'. For ROMs of size $r>10$, we observe that the CVX solver fails to make sufficient progress towards a solution, even to within the ``relaxed'' tolerance setting. To overcome this problem, we simplify the optimization problem by restricting the solution space of~\eqref{eq:constrained_general} to reduced mass matrix $\hbM=\mathbf{I}_r$ which still leads to a structure-preserving Lagrangian ROM that is inferred by fitting a Lagrangian ROM of the form~\eqref{eq:rom_learn} to the projections of the FOM data from~\eqref{eq:fom}.  

One way to circumvent the numerical challenges associated with solving~\eqref{eq:constrained_general} is to premultiply the ROM form in~\eqref{eq:rom_general_form} by $\hbM^{-1}$ and then solve for $\hbM^{-1}\hbK$ and $\hbM^{-1}\hbC$. However, this premultiplication by $\hbM^{-1}$ destroys the Lagrangian structure as we do not have a unique way of recovering $\hbM^{-1}$, $\hbK$, and $\hbC$ from $\hbM^{-1}\hbK$ and $\hbM^{-1}\hbC$. Another approach to tackling~\eqref{eq:constrained_general} is to solve it iteratively by breaking it into two constrained optimization problems. The idea is to start with an initial guess for $\hbM=\hbM_0$ and solving for $\hbC$ and $\hbK$ followed by solving for $\hbM$. We observe that this iterative approach converges only when we start with an initial guess $\hbM_0$ close to the intrusive reduced mass matrix. However, such an approach would require access to FOM operators which is not possible in the nonintrusive setting considered herein.
 \begin{algorithm}
\caption{Lagrangian operator inference (\textsf{L-OpInf}) for mechanical systems with nonconservative forcing}
\begin{algorithmic}[1]
\Require Snapshot data $\bQ \in \real^{n \times K}$, $\bU \in \real^{m \times K}$, $\bY \in \real^{p \times K}$ arranged as in \eqref{eq:snapshot}-\eqref{eq:snapshot_input} and reduced dimension $r$. 
\Ensure Reduced operators $\hbK,\hbC,\hbB$ for Lagrangian ROM~\eqref{eq:rom_learn}, and $\hbE$ \eqref{eq:rom_output}.
    \State Use knowledge of $L(\bq,\dot{\bq})$ to identify correct model form for the reduced Lagrangian $\hL_r$ in \eqref{eq:Lhat}.
    \State Build basis matrix $\bV_r \in \real^{n \times r}$ from SVD of $\bQ$ \eqref{eq:svd}.
    \State Project to obtain reduced state data $\hbQ \in \real^{r \times K}$ \eqref{eq:qhat}.
    \State Compute  reduced time-derivative data $\hbQd,\hbQdd \in \real^{r \times K}$ as in \eqref{eq:qhatdot} from the projected data $\hbQ$, i.e., by using the eighth-order central finite-difference scheme \eqref{eq:first}-\eqref{eq:second}.
    \State Solve constrained linear least-squares problems \eqref{eq:lopinf}-\eqref{eq:lopinf_output} to nonintrusively infer reduced operators $\hbK,\hbC,\hbB$, and $\hbE$.
\end{algorithmic}
\label{alg:lopinf_simple}
\end{algorithm}
%
%
%
 \begin{algorithm}
\caption{Lagrangian operator inference (\textsf{L-OpInf}) for nonlinear wave equations}
\begin{algorithmic}[1]
\Require Snapshot data $\bQ \in \real^{n \times K}$ arranged as in \eqref{eq:snapshot} and reduced dimension $r$. 
\Ensure Reduced operator $\hbK$ for nonlinear Lagrangian ROM~\eqref{eq:rom_learn_non}.
    \State Use knowledge of $U_{\text{nl}}(q)$ to identify correct model form for the reduced Lagrangian $\hL_r$ in \eqref{eq:Lhat_non}.
    \State Build nonlinear forcing snapshot data $\bF_{\text{nl}}\in \real^{n \times K}$~\eqref{eq:snapshot_non}.
    \State Build basis matrix $\bV_r \in \real^{n \times r}$ from SVD of $\bQ_{\text{aug}}=[\bQ, \bF_{\text{nl}}]\in \real^{n \times 2K}$.
    \State Project to obtain reduced state data $\hbQ \in \real^{r \times K}$  and reduced nonlinear forcing data $\hbF_{\text{nl}} \in \real^{r \times K}$ \eqref{eq:snapshot_proj}.
    \State Compute  reduced time-derivative data $\hbQdd \in \real^{r \times K}$ as in \eqref{eq:qhatdot} from the projected data $\hbQ$, i.e., by using the eighth-order central finite-difference scheme \eqref{eq:second}.
    \State Solve symmetric linear least-squares problem \eqref{eq:lopinf_non} to nonintrusively infer reduced operator $\hbK$.
\end{algorithmic}
\label{alg:lopinf_non}
\end{algorithm}

%
%
%
%

\section{Numerical results} \label{sec:numerical}
In this section, we study the numerical performance of \textsf{L-OpInf} for three Lagrangian systems. The reported error measures are detailed in Section~\ref{sec:implementation}. In Section~\ref{sec:beam} we revisit the Euler-Bernoulli beam model from structural dynamics and demonstrate that \textsf{L-OpInf} produces accurate Lagrangian ROMs for conservative mechanical systems in high dimensions. In Section~\ref{sec:sg} we consider the sine-Gordon equation to show the effectiveness of \textsf{L-OpInf} for nonlinear wave equations with nonpolynomial nonliearities. In Section~\ref{sec:fishtail} we consider a large-scale model of a soft-robotic fishtail with $n = 779,232$ DOFs to investigate the numerical performance of \textsf{L-OpInf} for mechanical systems with dissipation and external forcing.

\subsection{Error measures for accuracy and structure preservation}
\label{sec:implementation}
The state and output error plots reported in this section compute the 
\begin{equation}
\text{Relative state error}=\frac{\lVert \bQ - \bV_r\hQ \rVert_F}{\lVert \bQ \rVert_F}, \quad \quad \text{Relative output error}= \frac{\lVert \bY - \hY \rVert_F}{\lVert \bY \rVert_F},
\label{eq:errors}
\end{equation}
where $\bQ$ and $\bY$ are obtained from the Lagrangian FOM (e.g., of the form~\eqref{eq:fom} or~\eqref{eq:non_FOM}), and $\hQ$ and $\hY$ are either obtained from the nonintrusive Lagrangian ROM or the intrusive Lagrangian ROM. When reporting approximation errors in the training phase, we only consider trajectories in the training time interval $[0,T_{\text{train}}]$. For test data plots, we consider trajectories starting from the end of the training time interval to the end of total simulation time $T$, i.e. testing time interval $[T_{\text{train}},T]$.

The energy error plots reported in Section~\ref{sec:beam} and Section~\ref{sec:sg} compute the 
\begin{equation}
   \text{FOM energy error}= |E(\bV_r\hq(t),\bV_r\dot{\hq}(t)) - E(\bV_r\hq(0),\bV_r\dot{\hq}(0)) |,
    \label{eq:fom_energy_error}
\end{equation}
where $E(\bV_r\hq(t),\bV_r\dot{\hq}(t))$ is the FOM energy approximation (see equation \eqref{eq:FOM_energy}) either obtained from the nonintrusive Lagrangian ROM or the intrusive Lagrangian ROM. For the nonconservative soft-robotic fishtail example with varying energy, we compare the time evolution of FOM energy $E(\bq(t),\dot{\bq}(t))$ with $E(\bV_r\hq(t),\bV_r\dot{\hq}(t))$ to understand how well the ROMs track the energy. 

\subsection{Transverse vibrations of an Euler-Bernoulli beam}
\label{sec:beam}
The Euler-Bernoulli beam theory, also known as the classical beam theory, is most commonly used for calculating load-carrying and deflection characteristics of beams in structural and mechanical engineering. This linear beam theory is based on the observation that bending effect plays a key role in modeling transverse vibrations of beams. This theory is a  simplification of the linear elasticity theory and ignores the effects of shear deformation and rotary inertia. 
\subsubsection{PDE formulation}
\label{sec:beam_pde}
We consider the transverse vibrations of an unforced 1-D linear beam. The general dynamic equation for a conservative Euler-Bernoulli beam is given by 
\begin{equation*}
        \frac{\partial^2}{\partial x^2}\left(EI \frac{\partial^2 w(x,t)}{\partial x^2}\right)=-\mu \frac{\partial^2 w(x,t)}{\partial t^2},
\end{equation*}
where $x \in [0,\ell]$ is the spatial variable, $w(x,t)$ is the transverse deflection, the product $EI$ is the flexural rigidity, and $\mu$ is the mass per unit length. For a homogeneous beam with constant flexural rigidity, the governing PDE for modeling the transverse vibrations is
\begin{equation}
        EI \frac{\partial^4 w(x,t)}{\partial x^4}=-\mu \frac{\partial^2 w(x,t)}{\partial t^2}.
\end{equation}
The beam, as shown in Figure~\ref{fig:LB_schematic}, is simply-supported at both ends, i.e.,
\begin{equation}
w(0,t)=0, \quad w(\ell,t)=0,\quad \frac{\partial ^2 w}{\partial x^2}(0,t)=0, \quad \frac{\partial ^2 w}{\partial x^2}(\ell,t)=0.
\end{equation} 
In this study, we consider the following stationary initial condition
\begin{equation*}
    w(x,0)=0.001\left( x^4 - 2\ell x^3 + \ell^3x \right), \quad \quad \frac{\partial w}{\partial t}(x,0)=0.
    \label{eq:beam_ic}
\end{equation*}

\subsubsection{FOM and ROM implementation}
\label{sec:beam_fom}
We consider a steel beam of length $\ell=1$~\si{\metre} with mass per unit length $\mu=6.28 \times 10^{-3}$~\si{\kilogram \per\metre} and flexural rigidity $EI=9.81 \times 10^{-3}$~\si{\newton\metre\squared}. Using Hermite shape functions for beam finite elements, the governing PDE is spatially discretized to yield a Lagrangian FOM of the form
\begin{equation}
   \underbrace{ \begin{bmatrix}
    \bM_{ww} & \bM_{w\theta} \\ \bM_{w\theta} & \bM_{\theta\theta}
    \end{bmatrix}}_{\bM} \begin{bmatrix}
    \mathbf{\ddot{w}} \\ \pmb{\ddot{\theta}} 
    \end{bmatrix}+ \underbrace{\begin{bmatrix}
    \bK_{ww} & \bK_{w\theta} \\ \bK_{w\theta} & \bK_{\theta\theta}
    \end{bmatrix}}_{\bK}  \begin{bmatrix}
    \bw \\ \pmb{\theta} 
    \end{bmatrix} = \bzero,  
    \label{eq:beamfom}
\end{equation}
where $\bw \in \real^n$ is the deflection perpendicular to the beam length and $\pmb{\theta} \in \real^n$ is the rotation in the deformation plane. The state vector $\bq$ for this beam FOM is partitioned as $\bq=[\bw^\top,\pmb{\theta}^\top]^\top$ where both $\bw$ and $\pmb{\theta}$ carry specific physical meaning. The FOM mass and stiffness matrices in \eqref{eq:beamfom} also possess a block structure that reflects the partitioning of the state vector $\bq \in \real^{2n}$. 

\par 

We choose $n=200$ equally spaced grid points leading to a discretized state $\bq \in \real^{400}$. The FOM is numerically integrated for total time $T=0.18$~\si{\s} using a variational integrator based on the midpoint rule with $\Delta t=10^{-5}$~\si{\s}. The resulting time-marching equations require solving a linear system of $2n=400$ equations at every time step.

\par 

For the projection step in Step 3 of Algorithm \ref{alg:lopinf_simple}, we choose a projection matrix with a block diagonal structure, i.e. we approximate
\begin{equation}
    \begin{bmatrix}
    \bw \\ \pmb{\theta} 
    \end{bmatrix} \approx \underbrace{ \begin{bmatrix}
    \bV_{w} & \bzero \\ \bzero & \bV_{\theta}
    \end{bmatrix}}_{\bV_{2r}} \begin{bmatrix}
    \hat{\bw} \\ \pmb{\hat{\theta}} 
    \end{bmatrix}.
    \label{eq:beam_basis}
\end{equation}
This specific block diagonal structure of the basis matrix $\bV_{2r} \in \real^{2n \times 2r}$ retains the physical meaning of the $\bw$ and $\pmb \theta$ variables in the reduced dimensions. The same basis matrix is also used for the intrusive projection-based ROMs which we show for comparison of our results.

\subsubsection{Results} \label{sec:beam_results}

\begin{figure}
\captionsetup[subfigure]{oneside,margin={1.8cm,0 cm}}
\begin{subfigure}{.38\textwidth}
       \setlength\fheight{6 cm}
        \setlength\fwidth{\textwidth}
%
%
\begin{tikzpicture}

\begin{axis}[%
width=0.951\fheight,
height=0.536\fheight,
at={(0\fheight,0\fheight)},
scale only axis,
xmin=2,
xmax=20,
xlabel style={font=\color{white!15!black}},
xlabel={\small Reduced dimension $2r$},
ymode=log,
ymin=0.01,
ymax=2,
yminorticks=true,
ylabel style={font=\color{white!15!black}},
ylabel={\small Relative state error},
axis background/.style={fill=white},
xmajorgrids,
ymajorgrids,
legend style={at={(0.75,1.15)}, anchor=south west, legend cell align=left, align=left, draw=white!15!black,font=\small}
]
\addplot [color=green, dotted, line width=3.0pt, mark size=4.0pt, mark=x, mark options={solid, green}]
  table[row sep=crcr]{%
2	0.134645093182097\\
4	0.0596371586385031\\
6	0.0356015158560238\\
8	0.0245118185355717\\
10	0.0184610050037561\\
12	0.014860777683229\\
14	0.0129921318684437\\
16	0.0124702481277441\\
18	0.0140806145054294\\
20	0.0143597184812163\\
};
\addlegendentry{L-OpInf}

\addplot [color=blue, dotted, line width=3.0pt, mark size=4.0pt, mark=o, mark options={solid, blue}]
  table[row sep=crcr]{%
2	1.33558174918106\\
4	0.220940714893537\\
6	0.117335316976932\\
8	0.0871921758296532\\
10	0.063730922134043\\
12	0.0518052422133428\\
14	0.0421416816408257\\
16	0.0365380057213715\\
18	0.0305183396623845\\
20	0.0278349600715258\\
};
\addlegendentry{Intrusive Lagrangian ROM}

\end{axis}
\end{tikzpicture}%
\caption{Training regime $[0, 0.03]$~\si{\s}}
\label{fig:LB_state_train}
    \end{subfigure}
    \hspace{1.4cm}
    \begin{subfigure}{.38\textwidth}
           \setlength\fheight{6 cm}
           \setlength\fwidth{\textwidth}
\raisebox{-56mm}{
%
%
\begin{tikzpicture}

\begin{axis}[%
width=0.951\fheight,
height=0.536\fheight,
at={(0\fheight,0\fheight)},
scale only axis,
xmin=2,
xmax=20,
xlabel style={font=\color{white!15!black}},
xlabel={\small Reduced dimension $2r$},
ymode=log,
ymin=0.01,
ymax=2,
ylabel style={font=\color{white!15!black}},
ylabel={\small Relative state error},
axis background/.style={fill=white},
xmajorgrids,
ymajorgrids,
legend style={at={(0.511,0.699)}, anchor=south west, legend cell align=left, align=left, draw=white!15!black}
]
\addplot [color=green, dotted, line width=3.0pt, mark size=4.0pt, mark=x, mark options={solid, green}]
  table[row sep=crcr]{%
2	0.200479867331029\\
4	0.0862607036437419\\
6	0.0522467317795059\\
8	0.0379741020920877\\
10	0.0344379434904288\\
12	0.0354281742955599\\
14	0.0408201405381367\\
16	0.0419864529105539\\
18	0.0441369992782738\\
20	0.0436320176604149\\
};

\addplot [color=blue, dotted, line width=3.0pt, mark size=4.0pt, mark=o, mark options={solid, blue}]
  table[row sep=crcr]{%
2	1.45378467039996\\
4	0.359682317083506\\
6	0.313741003636607\\
8	0.201357127597659\\
10	0.159497333119719\\
12	0.123896703546996\\
14	0.103088037904053\\
16	0.0859280520969877\\
18	0.0741190369473179\\
20	0.0658214922004747\\
};

\end{axis}

\begin{axis}[%
width=1.227\fheight,
height=0.658\fheight,
at={(-0.16\fheight,-0.072\fheight)},
scale only axis,
xmin=0,
xmax=1,
ymin=0,
ymax=1,
axis line style={draw=none},
ticks=none,
axis x line*=bottom,
axis y line*=left
]
\end{axis}
\end{tikzpicture}%
}
\caption{Testing regime $[0.03, 0.18]$~\si{\s}}
\label{fig:LB_state_test}
    \end{subfigure}
\caption{Euler-Bernoulli beam: \textsf{L-OpInf} ROMs achieve lower state error than intrusive Lagrangian ROMs in both training and test interval, yet in the testing interval the state errors level off after $r\geq 12$.}
 \label{fig:LB_state}
\end{figure}
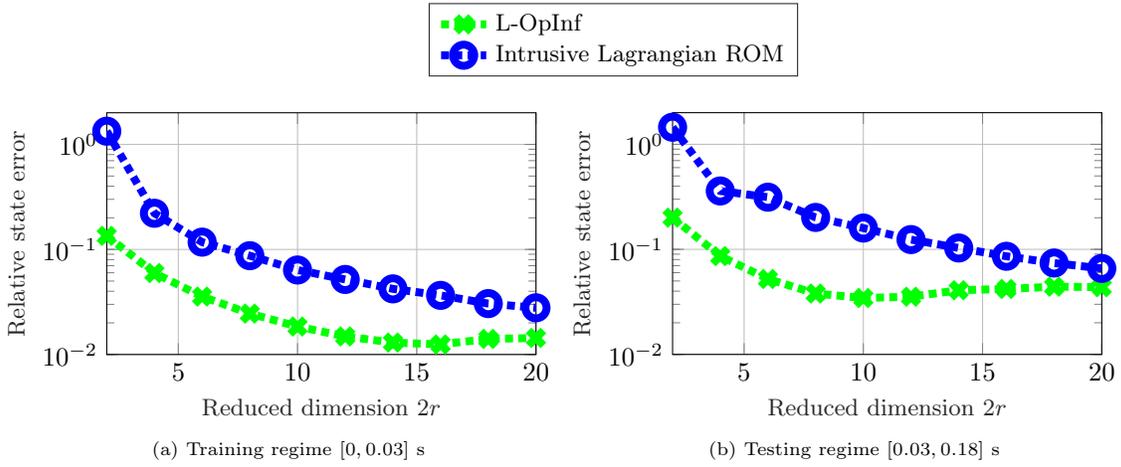

Figure~\ref{fig:LB_state} shows a comparison of the numerical performance of the intrusive and nonintrusive Lagrangian ROMs. The state error plots in Figure~\ref{fig:LB_state_train} over the training time interval $[0,0.03]$~\si{\s} show that the proposed nonintrusive approach performs better than the intrusive Lagrangian ROM. For the testing interval $[0.03, 0.18]$~\si{\s}, the nonintrusive Lagrangian ROMs perform better than the intrusive Lagrangian ROMs for all $2r$-dimensional reduced models in Figure~\ref{fig:LB_state_test}. We also observe that in the testing regime, the state approximation error for learned Lagrangian ROMs does not decrease as favorably with increase in reduced dimension for $2r>12$.

\begin{figure}[h]
\centering
\setlength\fheight{8.5 cm}
\setlength\fwidth{\textwidth}
\input{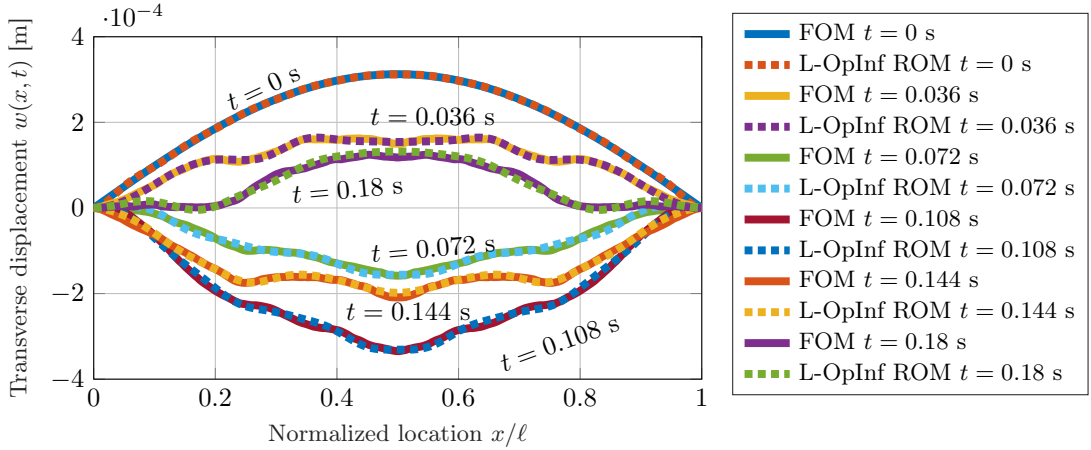}
 \caption{Euler-Bernoulli beam: Solutions using \textsf{L-OpInf} ROM of size $2r=20$ at different $t$ values (dashed lines) compared with FOM solutions at the same times (solid lines). The nonintrusive Lagrangian ROM captures the correct beam shape even at $t=0.18$~\si{\s}, which is 500\% outside training interval in a purely predictive setting.}
 \label{fig:LB_shape}
\end{figure}

Figure~\ref{fig:LB_shape} compares the \textsf{L-OpInf} ROM solution and the FOM solution at different $t$ values. Even though the reduced operators are learned from data in the training interval $[0, 0.03]$~\si{\s}, the nonintrusive Lagrangian ROM captures the correct beam shape at $t=0.18$~\si{\s} which is $500\%$ past the training time interval. The ROM solutions show very good qualitative and quantitative agreement with the FOM solutions. 

\par
%
%

\begin{figure}
\captionsetup[subfigure]{oneside,margin={1.8cm,0 cm}}
\begin{subfigure}{.37\textwidth}
       \setlength\fheight{5.8 cm}
\input{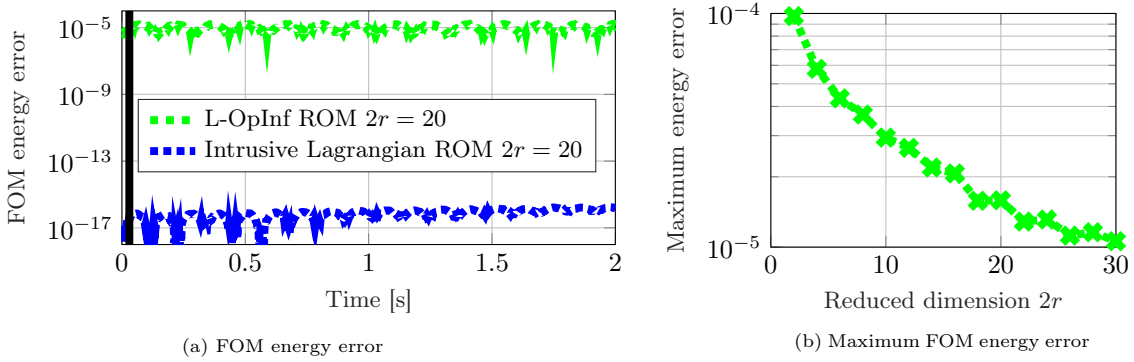}
\caption{FOM energy error}
 \label{fig:LB_Efom}
    \end{subfigure}
    \hspace{2.8cm}
\begin{subfigure}{.37\textwidth}
    \setlength\fheight{5.8 cm}
        \setlength\fwidth{\textwidth}
%
%
\begin{tikzpicture}

\begin{axis}[%
width=0.791\fheight,
height=0.536\fheight,
at={(0\fheight,0\fheight)},
scale only axis,
xmin=0,
xmax=30,
xlabel style={font=\color{white!15!black}},
xlabel={\small Reduced dimension $2r$},
ymode=log,
ymin=1e-05,
ymax=1e-04,
yminorticks=true,
ylabel style={font=\color{white!15!black}},
ylabel={\small Maximum energy error},
axis background/.style={fill=white},
xmajorgrids,
ymajorgrids,
yminorgrids
]
\addplot [color=green, dotted, line width=3.0pt, mark size=4.0pt, mark=x, mark options={solid, green}, forget plot]
  table[row sep=crcr]{%
2	9.74388665576276e-05\\
4	5.79441757399237e-05\\
6	4.33250837919472e-05\\
8	3.69325773148572e-05\\
10	2.94268622764858e-05\\
12	2.66016884213425e-05\\
14	2.1917509644766e-05\\
16	2.06800237510175e-05\\
18	1.57949640353719e-05\\
20	1.58631792730347e-05\\
22	1.28661003636763e-05\\
24	1.31332273975366e-05\\
26	1.12087077018334e-05\\
28	1.15783575711705e-05\\
30	1.06012197173616e-05\\
};
\end{axis}
\end{tikzpicture}%
\caption{Maximum FOM energy error}
\label{fig:LB_Emax}
    \end{subfigure}
\caption{Euler-Bernoulli beam: (a) The \textsf{L-OpInf} ROM exhibits bounded FOM energy error \eqref{eq:fom_energy_error} of approximately $10^{-5}$ far outside the training interval. The intrusive Lagrangian ROM shows exact energy preservation due to its specific choice of reduced Lagrangian. The black line in plot (a) indicates end of training time interval at $t=0.03$~\si{\s}. (b) The maximum FOM energy error for the \textsf{L-OpInf} ROM decreases with increase in reduced dimension $2r$.}
 \label{fig:LB_energy}
\end{figure}

In Figure~\ref{fig:LB_Efom}, we compare the FOM energy error for the \textsf{L-OpInf} ROM and the intrusive Lagrangian ROM, both simulated for $T=2$~\si{\s} and of dimension $2r=20$, to demonstrate the stability of nonintrusive Lagrangian ROMs far outside the training regime. Due to its specific choice of reduced Lagrangian, the intrusive Lagrangian ROM conserves the energy with the same accuracy as the FOM simulation. The nonintrusive ROM of dimension $2r=20$ exhibits bounded energy error due to its Lagrangian nature. The bounded energy error at $t=2$~\si{\s} ($6,500\%$ past the training interval) suggests that the nonintrusive Lagrangian ROMs simulate a perturbation of the intrusive Lagrangian ROM exactly. Thus, the FOM energy error for nonintrusive Lagrangian ROM of dimension $2r=20$ remains bounded well beyond the training data. In fact, every nonintrusive ROM learned via \textsf{L-OpInf} from $20\leq2r\leq 40$ demonstrates bounded energy error in the testing data. This shows a true strength of the proposed \textsf{L-OpInf}, namely that if the Lagrangian structure is respected in every aspect of discretization and the learning method, then long-term stable predictions are possible. 
This is in stark contrast to the standard second-order operator inference approach in Section~\ref{sec:motivation} which learned unstable ROMs from dimension $2r=20$ to $2r=40$. Compared to the the standard second-order operator inference results in Figure~\ref{fig:motivation}, \textsf{L-OpInf} shows bounded energy error behavior while also approximating the FOM state with similar accuracy. This emphasizes that \textsf{L-OpInf} learns the underlying Lagrangian dynamics rather than mere interpolations between training data snapshots.

\par

Figure~\ref{fig:LB_Emax} shows the maximum FOM energy error as a function of reduced dimension for \textsf{L-OpInf} ROMs for the entire length of simulation. The maximum FOM energy error decreases with increasing reduced dimension $2r$ which is in agreement with the state error results in Figure~\ref{fig:LB_state}. The errors level off after $2r=26$ at approximately $10^{-5}$, which is accurate enough for most structural dynamics applications.

We note that for $2r>16$, the relative state error levels off in Figure~\ref{fig:LB_state_train}. This stagnation occurs because the projected trajectories correspond to non-Markovian dynamics in the reduced setting even though the Lagrangian FOM dynamics are Markovian. The state error leveling-off for operator inference has been resolved by a re-projection sampling scheme that works for fully discrete systems with explicit time-marching schemes~\cite{peherstorfer2020sampling}. 
If the data is processed with that scheme, the learned models recover the intrusive ROMs preasymptotically under certain conditions. However, the large-scale Lagrangian dynamical systems considered herein require fully implicit time integrators to preserve the underlying geometric structure. Thus, re-projection in its current form cannot be used. Extending this algorithm to fully implicit and structure-preserving integrators remains an open problem.
%
\begin{remark}
From Figure~\ref{fig:motivation_state}, we observe that the nonintrusive ROMs obtained via the unconstrained second-order operator inference method demonstrate a higher accuracy in the training data compared to the nonintrusive structure-preserving \textsf{L-OpInf} ROMs in Figure~\ref{fig:LB_state_train}. This is to be expected, as the unconstrained inference problem from~\eqref{eq:std_opinf} solves for the reduced operator $\hat{\bK}$ without any constraints whereas the \textsf{L-OpInf} ROM is obtained by solving the constrained operator inference problem~\eqref{eq:lopinf} to ensure that the ROM is Lagrangian. The unconstrained second-order operator inference method learns reduced operators that overfit the data, and as a result, yields ROMs with lower state errors than the \textsf{L-OpInf} ROMs in the training data but they violate the underlying Lagrangian structure which leads to an unbounded energy error growth in Figure~\ref{fig:motivation_energy}. In contrast, the \textsf{L-OpInf} ROMs, due to their Lagrangian nature, yield bounded energy error in Figure~\ref{fig:LB_Efom}.
\end{remark}

In Figure~\ref{fig:LB_center}, we compare the transverse displacement of the beam at $x/l=0.5$ for the nonintrusive \textsf{L-OpInf} ROM and the nonintrusive second-order ROM, both simulated for $T=3$~\si{\s} and of dimension $2r=20$, to demonstrate the predictive capability of \textsf{L-OpInf} ROMs. Figure~\ref{fig:LB_center} shows that the \textsf{L-OpInf} method provides accurate and stable predictions in the testing data whereas the second-order operator inference method yields inaccurate solutions that eventually blow up in the predictive regime. The unphysical predictions in Figure~\ref{fig:LB_center} and the unbounded energy error growth in Figure~\ref{fig:motivation_energy} for the second-order operator inference method highlight that preserving physical properties is crucial for accurate and stable predictions and only looking at accuracy in the training data could be misleading.

\begin{figure}
\centering
\setlength\fheight{8.5 cm}
\setlength\fwidth{\textwidth}
\input{Figures/added_figure/center_comp.tex}
 \caption{Euler-Bernoulli beam: The \textsf{L-OpInf} method provides accurate and stable predictions far outside the training data regime whereas the second-order OpInf ROM learned via the unconstrained inference problem from~\eqref{eq:std_opinf} yields solutions that become unstable in the testing data. The black line indicates end of training time interval at $t=0.03$~\si{\s}.}
 \label{fig:LB_center}
\end{figure}
%
\begin{figure}
\captionsetup[subfigure]{oneside,margin={1.8cm,0 cm}}
\begin{subfigure}{.38\textwidth}
       \setlength\fheight{6 cm}
        \setlength\fwidth{\textwidth}
%
%
\begin{tikzpicture}

\begin{axis}[%
width=0.951\fheight,
height=0.536\fheight,
at={(0\fheight,0\fheight)},
scale only axis,
xmin=2,
xmax=20,
xlabel style={font=\color{white!15!black}},
xlabel={\small Reduced dimension $2r$},
ymode=log,
ymin=0.01,
ymax=2,
yminorticks=true,
ylabel style={font=\color{white!15!black}},
ylabel={\small Relative state error},
axis background/.style={fill=white},
xmajorgrids,
ymajorgrids,
legend style={at={(0.75,1.15)}, anchor=south west, legend cell align=left, align=left, draw=white!15!black,font=\small}
]
\addplot [color=green, dotted, line width=3.0pt, mark size=4.0pt, mark=x, mark options={solid, green}]
  table[row sep=crcr]{%
2	0.135264612198233\\
4	0.0605161486093417\\
6	0.0369435688050082\\
8	0.0263729228740884\\
10	0.0207936686796834\\
12	0.0176194755513933\\
14	0.0185625000739138\\
16	0.0184474777758963\\
18	0.0150738138514901\\
20	0.0159988741808267\\
};
\addlegendentry{L-OpInf}
\addplot [color=blue, dotted, line width=3.0pt, mark size=4.0pt, mark=o, mark options={solid, blue}]
  table[row sep=crcr]{%
2	1.33532569526819\\
4	0.221777809087806\\
6	0.117964296489332\\
8	0.0886842916868088\\
10	0.0648603714032791\\
12	0.0531469640965796\\
14	0.0523102899189476\\
16	0.0428091271376596\\
18	0.0377420150661271\\
20	0.0322597042794849\\
};
\addlegendentry{Intrusive Lagrangian ROM}
\end{axis}
\end{tikzpicture}%
\caption{Training data}
\label{fig:LB_unknown_train}
    \end{subfigure}
    \hspace{1.4cm}
    \begin{subfigure}{.38\textwidth}
           \setlength\fheight{6 cm}
           \setlength\fwidth{\textwidth}
\raisebox{-56mm}{
%
%
\begin{tikzpicture}

\begin{axis}[%
width=0.951\fheight,
height=0.536\fheight,
at={(0\fheight,0\fheight)},
scale only axis,
xmin=2,
xmax=20,
xlabel style={font=\color{white!15!black}},
xlabel={\small Reduced dimension $2r$},
ymode=log,
ymin=0.01,
ymax=2,
yminorticks=true,
ylabel style={font=\color{white!15!black}},
ylabel={\small Relative state error},
axis background/.style={fill=white},
xmajorgrids,
ymajorgrids,
legend style={at={(0.511,0.699)}, anchor=south west, legend cell align=left, align=left, draw=white!15!black}
]
\addplot [color=green, dotted, line width=3.0pt, mark size=4.0pt, mark=x, mark options={solid, green}]
  table[row sep=crcr]{%
2	0.136698375611215\\
4	0.0632605859008834\\
6	0.0411183827387285\\
8	0.0317576088524211\\
10	0.0270964343271557\\
12	0.0246087628445978\\
14	0.0264426205340073\\
16	0.0312284839184275\\
18	0.0230141323331054\\
20	0.0222671953105724\\
};
\addplot [color=blue, dotted, line width=3.0pt, mark size=4.0pt, mark=o, mark options={solid, blue}]
  table[row sep=crcr]{%
2	1.33520703666044\\
4	0.222957955450289\\
6	0.11980345756632\\
8	0.0906914391143856\\
10	0.0671954799215056\\
12	0.056392467951478\\
14	0.0587310639428867\\
16	0.0494330407464298\\
18	0.0448201013985711\\
20	0.0400628637262197\\
};
\end{axis}

\begin{axis}[%
width=1.227\fheight,
height=0.658\fheight,
at={(-0.16\fheight,-0.072\fheight)},
scale only axis,
xmin=0,
xmax=1,
ymin=0,
ymax=1,
axis line style={draw=none},
ticks=none,
axis x line*=bottom,
axis y line*=left
]
\end{axis}
\end{tikzpicture}%
}
\caption{ Test data}
\label{fig:LB_unknown_test}
    \end{subfigure}
\caption{Euler-Bernoulli beam: The \textsf{L-OpInf} ROMs achieves higher accuracy than the intrusive Lagrangian ROMs for both training and test initial conditions. }
 \label{fig:LB_unknown}
\end{figure}
%
%

To further highlight the generalizability of the method, we consider a different prediction scenario where we train \textsf{L-OpInf} ROMs using multiple initial conditions and then study their accuracy for unseen initial conditions. In this study, we consider a parametric initial condition of the form
\begin{equation*}
    w(x,0)=0.001\left( ax^5 -2\ell x^4 + \left( \frac{10a}{3}-4 \right)\ell^2 x^3 + \left(6- \frac{10a}{3} -a\right) \ell^4x \right), \quad \quad \frac{\partial w}{\partial t}(x,0)=0,
\end{equation*}
where $a\in \real$ is a scalar parameter. We build a training dataset by simulating the FOM until $T=0.03$~\si{\s} for $a=0.2$, $a=0.3$, $a=0.4$, and $a=0.5$. We then derive \textsf{L-OpInf} ROMs of different sizes from this training dataset. We consider two test initial conditions based on $a=0.1$ and $a=0.6$ to evaluate how the \textsf{L-OpInf} ROMs generalize for initial conditions that are not included---and are even outside---the training dataset. The comparison of the relative state error~\eqref{eq:errors} between \textsf{L-OpInf} ROMs and intrusive Lagrangian ROMs is shown
in Figure~\ref{fig:LB_unknown_train} and Figure~\ref{fig:LB_unknown_test}  for the training and the test initial conditions, respectively. The comparison in Figure~\ref{fig:LB_unknown} shows that the data-driven \textsf{L-OpInf} ROMs yield lower relative state error than the intrusive Lagrangian ROMs for both training and test initial conditions. These results show that the \textsf{L-OpInf} ROMs are robust to perturbations in the initial conditions used for building the training dataset.
\subsection{Sine-Gordon equation}
\label{sec:sg}
The sine-Gordon equation is a nonlinear hyperbolic PDE with a nonpolynomial nonlinearity. Its name is a wordplay on its similarity with the well-known Klein-Gordon wave equation. The sine-Gordon equation is a universal model for combining the wave dispersion and the nonlinearity which is a periodic function of the field variable. The sine-Gordon equation is used for modeling nonlinear phenomena in a wide variety of physical applications such as the  self-induced transparency in nonlinear optics~\cite{mccall1969self}, propagation of fluxons in Josephson junctions between superconductors~\cite{josephson1962possible}, relativistic field theory~\cite{samuel1978grand}, hydrodynamics~\cite{coullet1986resonance}, and charge-density-wave conductors~\cite{rice1976weakly}.  
\subsubsection{PDE formulation}
We consider the one-dimensional sine-Gordon equation
\begin{equation}
\frac{\partial^2 q}{\partial t^2} = \frac{\partial^2 q}{\partial x^2} - \sin (q),
 \label{eq:sg}
\end{equation}
where $t$ is the nondimensional time unit and the field variable $q(x,t)$ has the meaning of phase in the respective physical setting. This equation can be formulated as a Lagrangian PDE with the following space-time continuous Lagrangian  
\begin{equation*}  
\mathcal{L}(x,q,q_x,q_t)=\frac{1}{2}\left(  \left( \frac{\partial q}{\partial t}\right)^2 - \left( \frac{\partial q}{\partial x}\right)^2\right) - (1-\cos(q)).
\end{equation*}
In this study, we consider periodic boundary conditions with the following initial conditions 
\begin{equation*}
q(x,0)=0, \quad \quad \frac{\partial q}{\partial t}(x,0)=\frac{4}{\cosh(x)}.
\end{equation*}

\subsubsection{FOM and ROM implementation}
We study the sine-Gordon equation over $x\in [-L/2,L/2]$ with $L=40$. The nonlinear PDE is spatially discretized using $n=2,000$ equally spaced grid points leading to a discretized state $\bq \in \real^{2,000}$. We discretize the space-time continuous Lagrangian which yields the following space-discretized Lagrangian 
\begin{equation*}
L(\bq,\dot{\bq})= \frac{1}{2}\dot{\bq}^\top\dot{\bq} + \frac{1}{2}\bq^\top\bK_{\text{fd}}\bq  - \sum_{i=1}^n \left(1-\cos(q_i)\right),
\end{equation*}
where $\bK_{\text{fd}}$ denotes the symmetric finite difference approximation for the spatial derivative $\partial_{xx}$. The resulting Lagrangian FOM is represented by the following second-order nonlinear ODE system
\begin{equation*} 
    \ddot{\bq}=\bK_{\text{fd}}\bq - \begin{bmatrix} \sin(q_1) \\\vdots \\ \sin (q_n)\end{bmatrix}.
\end{equation*}
The FOM is numerically integrated until time $T=25$ using a variational integrator based on the midpoint rule with $\Delta t=0.005$. The resulting time-marching equations for the FOM require solving a system of $n=2,000$ coupled nonlinear equations at every time step. The intrusive Lagrangian ROMs are derived by projecting the FOM onto $\bV_r$ ( see equation~\eqref{eq:rom_int_non}) whereas the nonintrusive Lagrangian ROMs are inferred from the FOM simulation data by solving the constrained optimization problem in equation~\eqref{eq:lopinf_non}. We numerically integrate the ROMs of size $r$ with a variational integrator based on the midpoint rule with $\Delta t=0.005$.
\subsubsection{Results}
\begin{figure}
\captionsetup[subfigure]{oneside,margin={1.8cm,0 cm}}
\begin{subfigure}{.38\textwidth}
       \setlength\fheight{6 cm}
        \setlength\fwidth{\textwidth}
%
%
\begin{tikzpicture}

\begin{axis}[%
width=0.951\fheight,
height=0.536\fheight,
at={(0\fheight,0\fheight)},
scale only axis,
xmin=0,
xmax=14,
xlabel style={font=\color{white!15!black}},
xlabel={\small Reduced dimension $r$},
ymode=log,
ymin=1e-08,
ymax=1,
yminorticks=true,
ylabel style={font=\color{white!15!black}},
ylabel={\small Relative state error},
axis background/.style={fill=white},
xmajorgrids,
ymajorgrids,
yminorgrids,
legend style={at={(0.75,1.15)}, anchor=south west, legend cell align=left, align=left, draw=white!15!black,font=\small}
]
\addplot [color=green, dotted, line width=3.0pt, mark size=4.0pt, mark=x, mark options={solid, green}]
  table[row sep=crcr]{%
1	0.124732207305284\\
2	0.0124403797281656\\
3	0.00197334203905294\\
4	0.0003422294040307\\
5	6.33260865233879e-05\\
6	1.23631788758146e-05\\
7	2.7728007630461e-06\\
8	1.03864334878972e-06\\
9	7.50968972459573e-07\\
10	6.43271297435287e-07\\
11	6.19779142377094e-07\\
12	6.13718606367844e-07\\
13	6.12886133515571e-07\\
14	6.12642699298105e-07\\
};
\addlegendentry{L-OpInf}


\addplot [color=blue, dotted, line width=3.0pt, mark size=4.0pt, mark=o, mark options={solid, blue}]
  table[row sep=crcr]{%
1	0.120080304454082\\
2	0.0644959763036192\\
3	0.0116892088213839\\
4	0.00242462013265923\\
5	0.000589840738903601\\
6	0.000136643888237454\\
7	3.20997794664524e-05\\
8	7.49188670719831e-06\\
9	3.48786029868137e-06\\
10	2.49571172727627e-06\\
11	1.23277322556857e-06\\
12	6.36944828118952e-07\\
13	5.23218393398983e-07\\
14	1.61755936191815e-07\\
};
\addlegendentry{Intrusive Lagrangian ROM}


\end{axis}

\begin{axis}[%
width=1.259\fheight,
height=0.743\fheight,
at={(-0.164\fheight,-0.097\fheight)},
scale only axis,
xmin=0,
xmax=1,
ymin=0,
ymax=1,
axis line style={draw=none},
ticks=none,
axis x line*=bottom,
axis y line*=left
]
\end{axis}
\end{tikzpicture}%
\caption{Training regime $[0, 5]$~\si{\s}}
\label{fig:SG_state_train}
    \end{subfigure}
    \hspace{1.4cm}
    \begin{subfigure}{.38\textwidth}
           \setlength\fheight{6 cm}
           \setlength\fwidth{\textwidth}
\raisebox{-56mm}{
%
%
\begin{tikzpicture}

\begin{axis}[%
width=0.951\fheight,
height=0.536\fheight,
at={(0\fheight,0\fheight)},
scale only axis,
xmin=0,
xmax=14,
xlabel style={font=\color{white!15!black}},
xlabel={\small Reduced dimension $r$},
ymode=log,
ymin=1e-03,
ymax=10,
yminorticks=true,
ylabel style={font=\color{white!15!black}},
ylabel={\small Relative state error},
axis background/.style={fill=white},
xmajorgrids,
ymajorgrids,
yminorgrids,
legend style={at={(0.153,0.162)}, anchor=south west, legend cell align=left, align=left, draw=white!15!black}
]

\addplot [color=green, dotted, line width=3.0pt, mark size=4.0pt, mark=x, mark options={solid, green}]
  table[row sep=crcr]{%
1	1.08345508061647\\
2	1.27724906091825\\
3	1.02321961304577\\
4	0.700923053887164\\
5	0.364932736322165\\
6	0.216949862191018\\
7	0.13965666253219\\
8	0.0850434273652803\\
9	0.0832672149111198\\
10	0.0710233246487764\\
11	0.036152907805433\\
12	0.0274366269226685\\
13	0.0187218429761372\\
14	0.0112595570643362\\
};


\addplot [color=blue, dotted, line width=3.0pt, mark size=4.0pt, mark=o, mark options={solid, blue}]
  table[row sep=crcr]{%
1	1.09307390463823\\
2	1.26041654594852\\
3	0.736291786810487\\
4	0.252653412804712\\
5	0.193092560378604\\
6	0.131572592601752\\
7	0.0933696216484394\\
8	0.0647868110270725\\
9	0.0438284925983995\\
10	0.0256095868714142\\
11	0.0119075975046953\\
12	0.00773001064260091\\
13	0.00742045392838723\\
14	0.004288534759978\\
};

\end{axis}

\begin{axis}[%
width=1.259\fheight,
height=0.743\fheight,
at={(-0.164\fheight,-0.097\fheight)},
scale only axis,
xmin=0,
xmax=1,
ymin=0,
ymax=1,
axis line style={draw=none},
ticks=none,
axis x line*=bottom,
axis y line*=left
]
\end{axis}
\end{tikzpicture}%
}
\caption{Testing regime $[5, 25]$~\si{\s}}
\label{fig:SG_state_test}
    \end{subfigure}
\caption{Sine-Gordon equation: (a) In the training time interval, \textsf{L-OpInf} ROMs achieve lower state error than intrusive Lagrangian ROMs up tp $r=12$. For $r>10$, we observe a leveling off of the state error in the training data regime for the \textsf{L-OpInf} ROM. (b) For test data, the intrusive Lagrangian ROMs achieve lower state error than the \textsf{L-OpInf} ROMs.  }
 \label{fig:SG_state}
\end{figure}
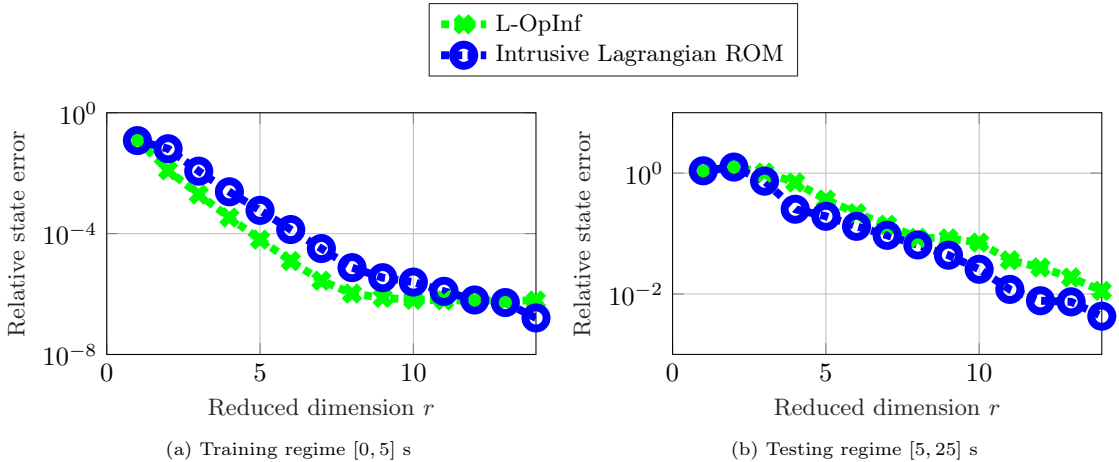
%
\begin{figure}
\captionsetup[subfigure]{oneside,margin={1cm,0 cm}}
\begin{subfigure}{.25\textwidth}
       \setlength\fheight{3 cm}
        \setlength\fwidth{\textwidth}
%
%
\definecolor{mycolor1}{rgb}{0.00000,0.44700,0.74100}%
\definecolor{mycolor2}{rgb}{0.85000,0.32500,0.09800}%
\definecolor{mycolor3}{rgb}{0.92900,0.69400,0.12500}%
\begin{tikzpicture}

\begin{axis}[%
width=1.251\fheight,
height=1\fheight,
at={(0\fheight,0\fheight)},
scale only axis,
xmin=-20,
xmax=20,
xlabel style={font=\color{white!15!black}},
xlabel={$x$},
ymin=-1,
ymax=7,
ylabel style={font=\color{white!15!black}},
ylabel={$q(x,t)$},
axis background/.style={fill=white},
xmajorgrids,
ymajorgrids,
legend style={at={(1.125,1.2)}, anchor=south west, legend cell align=left, align=left, draw=white!15!black,font=\footnotesize}
]
\addplot [color=red, line width=1.5pt]
  table[row sep=crcr]{%
-20	1.68308517345395e-07\\
-10.78	0.000831746819489609\\
-9.62	0.00265321617960979\\
-8.92	0.00534291662499342\\
-8.4	0.00898691455999412\\
-8	0.0134068476821376\\
-7.66	0.0188357861989168\\
-7.38	0.0249219260487834\\
-7.14	0.0316814995332422\\
-6.92	0.0394767283476085\\
-6.72	0.0482155765030967\\
-6.54	0.0577223645960636\\
-6.38	0.0677347557838566\\
-6.22	0.0794827039958328\\
-6.08	0.091420686609144\\
-5.94	0.105149881046831\\
-5.82	0.118545567584508\\
-5.7	0.133645414836295\\
-5.58	0.15066543883594\\
-5.48	0.166490076574345\\
-5.38	0.183973007311483\\
-5.28	0.203286756050257\\
-5.18	0.224621259055532\\
-5.1	0.243283490279747\\
-5.02	0.263488613794799\\
-4.94	0.285361896465208\\
-4.86	0.309038092429763\\
-4.78	0.334661990762882\\
-4.7	0.362388953734001\\
-4.62	0.392385428765365\\
-4.56	0.416478781874126\\
-4.5	0.442028094326627\\
-4.44	0.469116692881677\\
-4.38	0.497831790240046\\
-4.32	0.528264491668033\\
-4.26	0.560509759268815\\
-4.2	0.594666324435394\\
-4.14	0.630836537640917\\
-4.08	0.669126143211148\\
-4.02	0.709643965109002\\
-3.96	0.752501488089443\\
-3.9	0.797812316925853\\
-3.84	0.845691494869662\\
-3.78	0.896254661226244\\
-3.72	0.94961702708995\\
-3.66	1.00589214811685\\
-3.6	1.06519047400458\\
-3.54	1.12761765643855\\
-3.48	1.1932726010365\\
-3.42	1.26224525469705\\
-3.36	1.33461412815368\\
-3.3	1.41044356483523\\
-3.24	1.48978078160486\\
-3.18	1.57265272466417\\
-3.12	1.65906280465736\\
-3.06	1.74898759813316\\
-3	1.84237362689925\\
-2.94	1.93913435064407\\
-2.88	2.03914752919134\\
-2.82	2.14225312598325\\
-2.76	2.24825193068849\\
-2.68	2.39366569522787\\
-2.6	2.5431201782352\\
-2.5	2.73440937085582\\
-2.38	2.9681708596922\\
-2.16	3.39790703584921\\
-2.06	3.58903850975955\\
-1.98	3.73823508418251\\
-1.9	3.88320944921896\\
-1.84	3.98871799499738\\
-1.78	4.09116036736912\\
-1.72	4.1903041645781\\
-1.66	4.28595541318199\\
-1.6	4.37795768321063\\
-1.54	4.46619037293186\\
-1.48	4.55056633911157\\
-1.42	4.63102904939645\\
-1.36	4.70754942330675\\
-1.3	4.78012251050271\\
-1.26	4.82631869007797\\
-1.22	4.87077632799148\\
-1.18	4.91350759036314\\
-1.14	4.95452727091589\\
-1.1	4.99385240685308\\
-1.06	5.0315019155648\\
-1.02	5.06749625342489\\
-0.98	5.10185709737878\\
-0.940000000000001	5.13460704955473\\
-0.899999999999999	5.16576936473608\\
-0.859999999999999	5.19536770021589\\
-0.82	5.2234258873051\\
-0.780000000000001	5.24996772357167\\
-0.739999999999998	5.27501678475499\\
-0.699999999999999	5.2985962552026\\
-0.66	5.32072877562274\\
-0.620000000000001	5.34143630692417\\
-0.579999999999998	5.36074000891617\\
-0.539999999999999	5.37866013266755\\
-0.5	5.395215925362\\
-0.460000000000001	5.41042554654214\\
-0.419999999999998	5.42430599469486\\
-0.379999999999999	5.43687304320282\\
-0.34	5.44814118475865\\
-0.300000000000001	5.45812358341542\\
-0.259999999999998	5.46683203352514\\
-0.219999999999999	5.47427692489795\\
-0.18	5.48046721358731\\
-0.140000000000001	5.48541039779072\\
-0.0999999999999979	5.48911249842709\\
-0.0599999999999987	5.49157804402953\\
-0.0199999999999996	5.49281005966445\\
0.0199999999999996	5.49281005966445\\
0.0599999999999987	5.49157804402951\\
0.100000000000001	5.48911249842723\\
0.140000000000001	5.48541039779085\\
0.18	5.48046721358751\\
0.219999999999999	5.47427692489808\\
0.260000000000002	5.46683203352528\\
0.300000000000001	5.45812358341541\\
0.34	5.4481411847589\\
0.379999999999999	5.43687304320301\\
0.420000000000002	5.42430599469496\\
0.460000000000001	5.41042554654226\\
0.5	5.39521592536213\\
0.539999999999999	5.37866013266764\\
0.579999999999998	5.36074000891624\\
0.620000000000001	5.34143630692423\\
0.66	5.32072877562288\\
0.699999999999999	5.29859625520269\\
0.739999999999998	5.27501678475492\\
0.780000000000001	5.24996772357163\\
0.82	5.22342588730515\\
0.859999999999999	5.19536770021596\\
0.899999999999999	5.16576936473601\\
0.940000000000001	5.13460704955464\\
0.98	5.10185709737871\\
1.02	5.06749625342482\\
1.06	5.03150191556466\\
1.1	4.99385240685292\\
1.14	4.95452727091567\\
1.18	4.91350759036285\\
1.22	4.87077632799119\\
1.26	4.82631869007773\\
1.3	4.78012251050245\\
1.36	4.70754942330656\\
1.42	4.63102904939613\\
1.48	4.55056633911122\\
1.54	4.46619037293154\\
1.6	4.37795768321023\\
1.66	4.2859554131814\\
1.72	4.19030416457749\\
1.78	4.09116036736842\\
1.84	3.98871799499672\\
1.9	3.88320944921825\\
1.96	3.77490545314123\\
2.04	3.62668872696006\\
2.12	3.47487734533937\\
2.22	3.28146148058591\\
2.58	2.5810250905265\\
2.66	2.43067942823395\\
2.74	2.2841876522365\\
2.8	2.17727557397203\\
2.86	2.07318148212387\\
2.92	1.972117801639\\
2.98	1.87425751455668\\
3.04	1.77973525875051\\
3.1	1.68864920449522\\
3.16	1.60106353139579\\
3.22	1.51701133082762\\
3.28	1.43649777165701\\
3.34	1.35950338641566\\
3.4	1.28598735828687\\
3.46	1.21589071368474\\
3.52	1.14913934890288\\
3.58	1.08564684094449\\
3.64	1.02531701141377\\
3.7	0.968046227905194\\
3.76	0.913725439661636\\
3.82	0.862241953597234\\
3.88	0.813480963449329\\
3.94	0.767326849251379\\
4	0.723664266921052\\
4.06	0.682379048945055\\
4.12	0.643358937259425\\
4.18	0.606494168782532\\
4.24	0.571677932901231\\
4.3	0.538806718741832\\
4.36	0.507780568426018\\
4.42	0.478503250835914\\
4.48	0.450882368771527\\
4.54	0.424829410835372\\
4.6	0.400259757953695\\
4.68	0.369668553432611\\
4.76	0.341390260561653\\
4.84	0.315255585834269\\
4.92	0.291106442202928\\
5	0.268795447063543\\
5.08	0.248185388753523\\
5.16	0.229148679970862\\
5.24	0.211566810856418\\
5.34	0.191468578308214\\
5.44	0.173274965040239\\
5.54	0.156806660133498\\
5.64	0.141900882732934\\
5.76	0.125869568816697\\
5.88	0.111647229860445\\
6	0.0990302684158451\\
6.14	0.0860994234427928\\
6.28	0.0748558224039044\\
6.44	0.0637913714713179\\
6.6	0.0543616301869534\\
6.78	0.045408171525807\\
6.98	0.0371780358161615\\
7.2	0.0298366455524501\\
7.44	0.0234706514811336\\
7.72	0.0177389038679223\\
8.04	0.0128811648209179\\
8.42	0.00880896285348598\\
8.9	0.00545085044407401\\
9.52	0.00293225717656753\\
10.4	0.00121625051846763\\
11.88	0.000276864471000948\\
15.96	4.68463871428071e-06\\
19.98	1.68341461659338e-07\\
};
\addlegendentry{FOM}

\addplot [color=green, only marks, line width=1.5pt, mark size=2pt,mark=x, mark options={solid, green}]
  table[row sep=crcr]{%
-20	1.68286295121334e-07\\
-18	6.24352097133851e-07\\
-16	4.50133002516395e-06\\
-14	3.32322839362575e-05\\
-12	0.000245555219738236\\
-10	0.00181442107001928\\
-8	0.0134067649169403\\
-6	0.0990298164920063\\
-4	0.723662278281772\\
-2	3.70129801914601\\
0	5.49296123328359\\
2	3.70129801914503\\
4	0.723662278281516\\
6	0.0990298164920063\\
8	0.0134067649169403\\
10	0.00181442107001928\\
12	0.000245555219738236\\
14	3.32322839362575e-05\\
16	4.50133002516395e-06\\
18	6.24352097133851e-07\\
};
\addlegendentry{L-OpInf ROM $r=14$}

\addplot [color=blue, only marks, line width=1.5pt, mark size=3pt, mark=o, mark options={solid, blue}]
  table[row sep=crcr]{%
-20	1.68444657333566e-07\\
-17	1.66129839129781e-06\\
-14	3.32313660180716e-05\\
-11	0.00066746986169619\\
-8	0.0134064057748056\\
-5	0.268795167588301\\
-2	3.70130042514094\\
1	5.08487930516394\\
4	0.723666209213309\\
7	0.0364409333237816\\
10	0.00181437101737103\\
13	9.03322222320924e-05\\
16	4.50082704617216e-06\\
19	2.57769187328449e-07\\
};
\addlegendentry{Intrusive Lagrangian ROM $r=14$};

\end{axis}
\end{tikzpicture}%
\subcaption{$t=5$}
    \end{subfigure}
    \hspace{1.2cm}
    \begin{subfigure}{.25\textwidth}
           \setlength\fheight{3 cm}
           \setlength\fwidth{\textwidth}
\raisebox{-59mm}{
%
%
\definecolor{mycolor1}{rgb}{0.00000,0.44700,0.74100}%
\definecolor{mycolor2}{rgb}{0.85000,0.32500,0.09800}%
\definecolor{mycolor3}{rgb}{0.92900,0.69400,0.12500}%
\begin{tikzpicture}

\begin{axis}[%
width=1.251\fheight,
height=1\fheight,
at={(0\fheight,0\fheight)},
scale only axis,
xmin=-20,
xmax=20,
xlabel style={font=\color{white!15!black}},
xlabel={$x$},
ymin=-1,
ymax=7,
ylabel style={font=\color{white!15!black}},
axis background/.style={fill=white},
xmajorgrids,
ymajorgrids
]
\addplot [color=red, line width=1.5pt, forget plot]
  table[row sep=crcr]{%
-20	4.95076442064146e-07\\
-11.7	0.00098977459247962\\
-10.56	0.00313660699573504\\
-9.86	0.006301729347161\\
-9.34	0.0105809402300316\\
-8.94	0.0157690629393272\\
-8.6	0.0221411381289727\\
-8.32	0.0292847924034305\\
-8.08	0.0372192952869028\\
-7.86	0.0463699618999875\\
-7.66	0.0566288679792706\\
-7.48	0.067789839977717\\
-7.32	0.0795448970549977\\
-7.16	0.0933380945211866\\
-7.02	0.107354778467389\\
-6.88	0.123474737081462\\
-6.76	0.139202949362829\\
-6.64	0.156931423319964\\
-6.52	0.176912936508185\\
-6.42	0.195489165409814\\
-6.32	0.216009313869801\\
-6.22	0.238674384603804\\
-6.12	0.263705325269527\\
-6.04	0.285595746039874\\
-5.96	0.309290006918214\\
-5.88	0.334933120122503\\
-5.8	0.362680779900945\\
-5.72	0.392699884090934\\
-5.64	0.425169004213121\\
-5.56	0.46027877650199\\
-5.5	0.488465669079428\\
-5.44	0.518342079718558\\
-5.38	0.55000217397232\\
-5.32	0.583543990338878\\
-5.26	0.619069323241597\\
-5.2	0.656683537739525\\
-5.14	0.69649530272542\\
-5.08	0.738616227704231\\
-5.02	0.783160386541262\\
-4.96	0.830243709963032\\
-4.9	0.879983227183523\\
-4.84	0.932496135973793\\
-4.78	0.987898680004125\\
-4.72	1.04630481260183\\
-4.66	1.10782462751893\\
-4.6	1.17256254024388\\
-4.54	1.24061520823512\\
-4.48	1.31206918562051\\
-4.42	1.38699831780328\\
-4.36	1.46546089437675\\
-4.3	1.54749659494966\\
-4.24	1.63312328185933\\
-4.18	1.72233371586816\\
-4.12	1.81509229490827\\
-4.06	1.91133194028766\\
-4	2.01095127743827\\
-3.94	2.11381227660177\\
-3.86	2.25569230368123\\
-3.78	2.40249033550628\\
-3.7	2.55355691184319\\
-3.6	2.74721804206906\\
-3.48	2.98440812807997\\
-3.2	3.54018999656952\\
-3.1	3.73330588435303\\
-3.02	3.88379148594198\\
-2.94	4.02988811982096\\
-2.86	4.17095453019007\\
-2.8	4.27313658351383\\
-2.74	4.37202167308699\\
-2.68	4.4674747101725\\
-2.62	4.55939688263139\\
-2.56	4.64772322985273\\
-2.5	4.73241977659349\\
-2.44	4.81348038934564\\
-2.38	4.89092349989614\\
-2.32	4.96478881776265\\
-2.26	5.03513412878672\\
-2.2	5.1020322533437\\
-2.14	5.16556821577874\\
-2.08	5.22583665764855\\
-2.02	5.28293951151195\\
-1.96	5.336983939405\\
-1.9	5.38808053056983\\
-1.84	5.43634174611147\\
-1.78	5.48188059365328\\
-1.72	5.52480951230193\\
-1.66	5.56523944692417\\
-1.6	5.60327909052752\\
-1.54	5.63903427409696\\
-1.48	5.6726074843477\\
-1.42	5.70409749128569\\
-1.36	5.73359906908103\\
-1.3	5.76120279543965\\
-1.24	5.78699491632852\\
-1.18	5.81105726449529\\
-1.12	5.83346722172505\\
-1.04	5.86090218870212\\
-0.960000000000001	5.88568935117079\\
-0.879999999999999	5.90797556356546\\
-0.800000000000001	5.92789417857129\\
-0.719999999999999	5.94556537996906\\
-0.640000000000001	5.96109653978121\\
-0.559999999999999	5.97458258029611\\
-0.48	5.98610632591032\\
-0.399999999999999	5.99573883321797\\
-0.32	6.00353969052673\\
-0.239999999999998	6.00955728017454\\
-0.16	6.01382899871595\\
-0.0799999999999983	6.01638143140407\\
0	6.01723047845022\\
0.0799999999999983	6.01638143140425\\
0.16	6.01382899871598\\
0.239999999999998	6.00955728017475\\
0.32	6.00353969052689\\
0.399999999999999	5.99573883321828\\
0.48	5.98610632591075\\
0.559999999999999	5.97458258029637\\
0.640000000000001	5.9610965397816\\
0.719999999999999	5.94556537996947\\
0.800000000000001	5.92789417857156\\
0.879999999999999	5.90797556356573\\
0.960000000000001	5.88568935117108\\
1.04	5.86090218870243\\
1.12	5.83346722172528\\
1.18	5.81105726449561\\
1.24	5.78699491632888\\
1.3	5.76120279544011\\
1.36	5.73359906908152\\
1.42	5.70409749128626\\
1.48	5.6726074843481\\
1.54	5.63903427409738\\
1.6	5.60327909052787\\
1.66	5.56523944692441\\
1.72	5.52480951230216\\
1.78	5.4818805936535\\
1.84	5.43634174611168\\
1.9	5.38808053056995\\
1.96	5.33698393940533\\
2.02	5.28293951151228\\
2.08	5.22583665764893\\
2.14	5.16556821577906\\
2.2	5.10203225334399\\
2.26	5.03513412878702\\
2.32	4.96478881776297\\
2.38	4.89092349989652\\
2.44	4.81348038934622\\
2.5	4.73241977659401\\
2.56	4.64772322985329\\
2.62	4.55939688263174\\
2.68	4.46747471017311\\
2.74	4.37202167308751\\
2.8	4.2731365835144\\
2.86	4.17095453019047\\
2.92	4.06564868689687\\
3	3.92075421284115\\
3.08	3.77129922807154\\
3.16	3.61801275929436\\
3.26	3.42229741981591\\
3.42	3.10401341640172\\
3.58	2.78645957899755\\
3.68	2.59190326610764\\
3.76	2.43988478941545\\
3.84	2.2919540559703\\
3.92	2.14879028331208\\
3.98	2.04488600176282\\
4.04	1.94416890713494\\
4.1	1.84678966705834\\
4.16	1.75286164008942\\
4.22	1.66246304279553\\
4.28	1.57563963479376\\
4.34	1.49240775172936\\
4.4	1.41275753237742\\
4.46	1.33665620759403\\
4.52	1.26405134288754\\
4.58	1.19487395067042\\
4.64	1.12904141111244\\
4.7	1.06646016087605\\
4.76	1.00702812629393\\
4.82	0.950636891580647\\
4.88	0.897173603557086\\
4.94	0.846522622402084\\
5	0.798566933511424\\
5.06	0.753189339079423\\
5.12	0.710273449914308\\
5.18	0.669704498662451\\
5.24	0.631369995356881\\
5.3	0.59516024531866\\
5.36	0.560968748128118\\
5.42	0.528692494824984\\
5.48	0.498232178817009\\
5.54	0.469492334285018\\
5.62	0.433691169789292\\
5.7	0.400580237181266\\
5.78	0.369965835422981\\
5.86	0.341666409683778\\
5.94	0.3155121852086\\
6.02	0.291344714152334\\
6.1	0.269016364632016\\
6.2	0.243483862971601\\
6.3	0.220363961634451\\
6.4	0.199431523364765\\
6.5	0.180481688485262\\
6.6	0.163328239801835\\
6.72	0.144878312483996\\
6.84	0.128509641295214\\
6.96	0.113988448756832\\
7.1	0.0991061279855323\\
7.24	0.0861660027533908\\
7.4	0.0734325198924175\\
7.58	0.0613423501002117\\
7.76	0.0512431975517487\\
7.96	0.0419604331600567\\
8.18	0.0336805505180884\\
8.44	0.0259770501658174\\
8.74	0.0192529708855105\\
9.08	0.0137132751924085\\
9.48	0.00920257779818456\\
9.98	0.00559168540018362\\
10.64	0.00289577443449218\\
11.62	0.00107517614796038\\
13.34	0.00016997979433242\\
19.98	4.9517527855869e-07\\
};
\addplot [color=green, only marks, line width=1.5pt, mark size=2pt,mark=x, mark options={solid, green}]
  table[row sep=crcr]{%
-20	1.31869137476315e-06\\
-18	5.45097535109562e-07\\
-16	1.3965278121475e-05\\
-14	9.62118600114081e-05\\
-12	0.000710925543000229\\
-10	0.00525301825139479\\
-8	0.0387904512243473\\
-6	0.282022640774759\\
-4	2.0087790174779\\
-2	5.30333713398917\\
0	6.01599467895047\\
2	5.30333713404676\\
4	2.00877901756983\\
6	0.28202264087901\\
8	0.0387904512555295\\
10	0.00525301824959001\\
12	0.000710925542545482\\
14	9.62118599545647e-05\\
16	1.39652781108168e-05\\
18	5.45097517345994e-07\\
};
\addplot [color=blue, only marks, line width=1.5pt, mark size=3pt, mark=o, mark options={solid, blue}]
  table[row sep=crcr]{%
-20	1.28819936442426e-05\\
-17	6.2158554996472e-06\\
-14	9.98682332387091e-05\\
-11	0.00200582255914128\\
-8	0.0402803579571227\\
-5	0.794823818404947\\
-2	5.30014949894112\\
1	5.87336192071245\\
4	2.00785321377745\\
7	0.109367216768138\\
10	0.00545236709949037\\
13	0.000271462476387541\\
16	1.38645513239055e-05\\
19	1.19404164422576e-05\\
};
\end{axis}
\end{tikzpicture}%
}
\subcaption{$t=15$}
    \end{subfigure}
    \hspace{1cm}
        \begin{subfigure}{.25\textwidth}
           \setlength\fheight{3 cm}
           \setlength\fwidth{\textwidth}
\raisebox{-59mm}{
%
%
\begin{tikzpicture}

\begin{axis}[%
width=1.251\fheight,
height=1\fheight,
at={(0\fheight,0\fheight)},
scale only axis,
xmin=-20,
xmax=20,
xlabel style={font=\color{white!15!black}},
xlabel={$x$},
ymin=-1,
ymax=7,
ylabel style={font=\color{white!15!black}},
axis background/.style={fill=white},
xmajorgrids,
ymajorgrids
]
\addplot [color=red, line width=1.5pt]
  table[row sep=crcr]{%
-20	-2.95046549254607e-05\\
-12.1	0.00109635465581448\\
-10.96	0.00347664385176216\\
-10.26	0.00701924929420983\\
-9.74	0.0118157590338477\\
-9.34	0.0176318902192136\\
-9	0.0247746332948786\\
-8.72	0.0327813186939672\\
-8.46	0.0425153078438605\\
-8.24	0.0529765552567376\\
-8.04	0.0647038342994506\\
-7.86	0.0774613760251484\\
-7.7	0.0908970831953653\\
-7.54	0.106661232276856\\
-7.4	0.122679615471885\\
-7.26	0.141100129616316\\
-7.14	0.15907135757665\\
-7.02	0.17932604063305\\
-6.9	0.202151968244475\\
-6.8	0.223369628796885\\
-6.7	0.246804036772524\\
-6.6	0.272683219886545\\
-6.5	0.301257417254217\\
-6.42	0.326240585161074\\
-6.34	0.353275664057126\\
-6.26	0.382525814541154\\
-6.18	0.414165728738521\\
-6.1	0.448382052546773\\
-6.02	0.485373706618653\\
-5.94	0.525352065147001\\
-5.86	0.568540941372991\\
-5.8	0.603181015641066\\
-5.74	0.639863372967639\\
-5.68	0.678695319793491\\
-5.62	0.719787286325928\\
-5.56	0.763252403972128\\
-5.5	0.809205959978893\\
-5.44	0.857764710263723\\
-5.38	0.909046030184562\\
-5.32	0.963166882212533\\
-5.26	1.02024257939587\\
-5.2	1.08038532438895\\
-5.14	1.14370250604417\\
-5.08	1.21029473949531\\
-5.02	1.280253641718\\
-4.96	1.35365934317312\\
-4.9	1.43057774767581\\
-4.84	1.5110575673578\\
-4.78	1.59512717755682\\
-4.72	1.68279135742636\\
-4.66	1.77402800533314\\
-4.6	1.86878494250799\\
-4.54	1.96697694213446\\
-4.46	2.10302874465929\\
-4.38	2.24457683386483\\
-4.3	2.39110641632726\\
-4.22	2.541979767708\\
-4.12	2.73552158400829\\
-3.98	3.01263372906013\\
-3.72	3.5296333675223\\
-3.62	3.72345320790174\\
-3.54	3.87461226804599\\
-3.46	4.02147211651711\\
-3.38	4.1633829834223\\
-3.3	4.2998175161086\\
-3.24	4.39830211234463\\
-3.18	4.49335070210236\\
-3.12	4.58487318155712\\
-3.06	4.67281279226796\\
-3	4.75714318733288\\
-2.94	4.83786525312438\\
-2.88	4.91500382655954\\
-2.82	4.98860442443062\\
-2.76	5.05873007692749\\
-2.7	5.12545833390046\\
-2.64	5.188878491206\\
-2.58	5.24908906607176\\
-2.52	5.30619553525484\\
-2.46	5.36030833791837\\
-2.4	5.41154113608659\\
-2.34	5.46000931925149\\
-2.28	5.5058287355073\\
-2.22	5.54911462914059\\
-2.16	5.58998076363434\\
-2.1	5.62853870893336\\
-2.04	5.66489727260734\\
-1.98	5.69916205566869\\
-1.92	5.73143511533042\\
-1.84	5.77153690134386\\
-1.76	5.8084940133017\\
-1.68	5.8425154735795\\
-1.6	5.87379726195144\\
-1.52	5.90252238796507\\
-1.44	5.92886111795612\\
-1.36	5.95297131212241\\
-1.28	5.97499883614852\\
-1.2	5.99507801935636\\
-1.1	6.01762393421988\\
-1	6.03753510788979\\
-0.899999999999999	6.05500492046209\\
-0.800000000000001	6.07020384560887\\
-0.699999999999999	6.08328076621313\\
-0.599999999999998	6.09436417474406\\
-0.48	6.10518524212963\\
-0.359999999999999	6.11344575458923\\
-0.239999999999998	6.11926337103713\\
-0.120000000000001	6.12272105635682\\
0	6.12386815888964\\
0.120000000000001	6.12272105635695\\
0.239999999999998	6.11926337103716\\
0.359999999999999	6.11344575458959\\
0.48	6.10518524213015\\
0.600000000000001	6.09436417474458\\
0.699999999999999	6.08328076621369\\
0.800000000000001	6.07020384560927\\
0.899999999999999	6.05500492046253\\
1	6.03753510789033\\
1.1	6.01762393422051\\
1.2	5.99507801935697\\
1.28	5.97499883614884\\
1.36	5.95297131212258\\
1.44	5.92886111795617\\
1.52	5.90252238796531\\
1.6	5.87379726195151\\
1.68	5.84251547357971\\
1.76	5.80849401330213\\
1.84	5.77153690134439\\
1.92	5.73143511533099\\
1.98	5.69916205566956\\
2.04	5.66489727260863\\
2.1	5.62853870893486\\
2.16	5.58998076363588\\
2.22	5.54911462914228\\
2.28	5.50582873550932\\
2.34	5.46000931925335\\
2.4	5.41154113608876\\
2.46	5.36030833792044\\
2.52	5.3061955352568\\
2.58	5.24908906607387\\
2.64	5.18887849120811\\
2.7	5.12545833390282\\
2.76	5.05873007692998\\
2.82	4.98860442443317\\
2.88	4.91500382656242\\
2.94	4.8378652531276\\
3	4.75714318733626\\
3.06	4.67281279227139\\
3.12	4.58487318156061\\
3.18	4.49335070210606\\
3.24	4.39830211234849\\
3.3	4.29981751611247\\
3.38	4.16338298342643\\
3.46	4.02147211652121\\
3.54	3.87461226805012\\
3.62	3.72345320790603\\
3.72	3.52963336752669\\
3.86	3.25230902231913\\
4.12	2.73552158401267\\
4.22	2.54197976771234\\
4.3	2.39110641633145\\
4.38	2.24457683386882\\
4.46	2.10302874466329\\
4.54	1.96697694213826\\
4.6	1.86878494251174\\
4.66	1.77402800533676\\
4.72	1.68279135742996\\
4.78	1.59512717756026\\
4.84	1.51105756736116\\
4.9	1.4305777476791\\
4.96	1.35365934317629\\
5.02	1.28025364172107\\
5.08	1.21029473949831\\
5.14	1.14370250604699\\
5.2	1.08038532439169\\
5.26	1.02024257939852\\
5.32	0.963166882215216\\
5.38	0.909046030186964\\
5.44	0.857764710265965\\
5.5	0.809205959981\\
5.56	0.763252403974143\\
5.62	0.719787286327762\\
5.68	0.678695319795025\\
5.74	0.639863372968843\\
5.8	0.603181015642189\\
5.86	0.568540941374117\\
5.94	0.525352065148034\\
6.02	0.485373706619495\\
6.1	0.4483820525477\\
6.18	0.414165728739114\\
6.26	0.382525814541982\\
6.34	0.353275664057637\\
6.42	0.326240585161386\\
6.5	0.301257417254433\\
6.6	0.272683219886794\\
6.7	0.246804036772815\\
6.8	0.223369628797318\\
6.9	0.202151968245055\\
7.02	0.179326040633526\\
7.14	0.1590713575774\\
7.26	0.141100129617133\\
7.4	0.122679615472318\\
7.54	0.106661232276966\\
7.7	0.0908970831955855\\
7.86	0.0774613760254219\\
8.04	0.0647038342995607\\
8.24	0.0529765552565244\\
8.46	0.0425153078437361\\
8.7	0.0334435906445698\\
8.98	0.0252752176154587\\
9.3	0.0183517872497916\\
9.68	0.0125470910046452\\
10.14	0.00791606907420572\\
10.72	0.00442488197366586\\
11.54	0.00193709008027554\\
12.86	0.000500022031868497\\
15.96	1.51441727993529e-05\\
19.98	-2.94972762873158e-05\\
};
\addplot [color=green, only marks, line width=1.5pt, mark size=3pt,mark=x, mark options={solid, green}]
  table[row sep=crcr]{%
-20	3.34648671440618e-06\\
-18	-3.82681125543627e-06\\
-16	1.93976189883927e-05\\
-14	0.000134265977383308\\
-12	0.000992127646519236\\
-10	0.00733074739292405\\
-8	0.0541002121312708\\
-6	0.387205053530938\\
-4	2.73495419283641\\
-2	5.6244990770312\\
0	6.10058214222099\\
2	5.6244990789185\\
4	2.73495419263192\\
6	0.387205053620008\\
8	0.0541002121973513\\
10	0.00733074738836592\\
12	0.00099212764544987\\
14	0.000134265977258963\\
16	1.93976189706291e-05\\
18	-3.82681127319984e-06\\
};
\addplot [color=blue, only marks, line width=1.5pt, mark size=3pt,mark=o, mark options={solid, blue}]
  table[row sep=crcr]{%
-20	8.83907817019747e-05\\
-17	1.86293616408761e-05\\
-14	0.000160014725000934\\
-11	0.00321332535943952\\
-8	0.0645048733424609\\
-5	1.23324145496715\\
-2	5.65632302148465\\
1	6.02486712604291\\
4	2.86728721251054\\
7	0.17474625473735\\
10	0.00873458776877101\\
13	0.000434906422711379\\
16	2.33565579748074e-05\\
19	8.04125039763903e-05\\
};
\end{axis}
\end{tikzpicture}%
}
\subcaption{$t=25$}
    \end{subfigure}
    \caption{Sine-Gordon equation: Plots show the numerical approximation of the solution of \eqref{eq:sg} using low-dimensional ($r=14$) intrusive and nonintrusive Lagrangian ROMs at different $t$ values. The nonintrusive Lagrangian ROM captures the correct wave shape at $t=25$, which is 400\% outside the training time interval.}
\label{fig:SG_traj}
\end{figure}
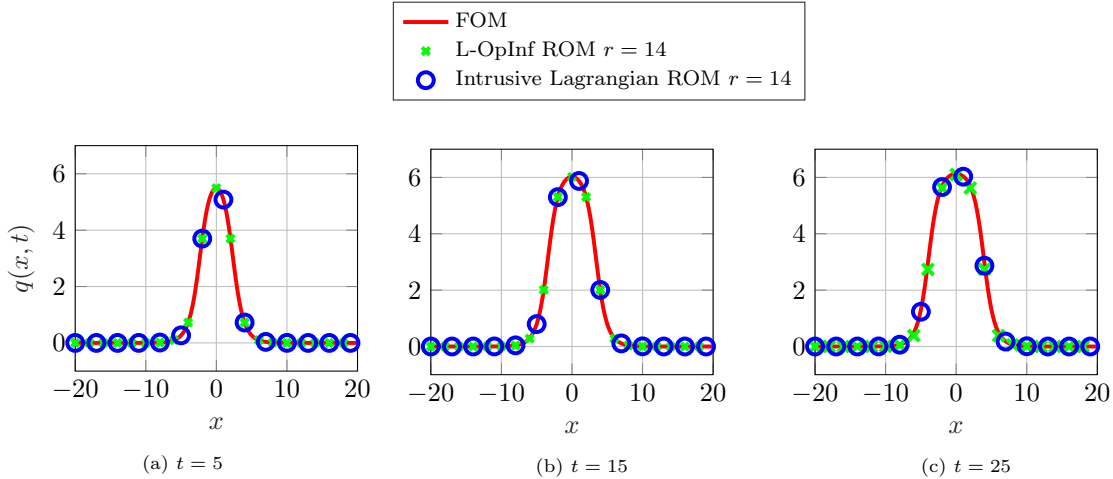
\begin{figure}[h]
\centering
\setlength\fheight{8.5 cm}
\setlength\fwidth{\textwidth}
\input{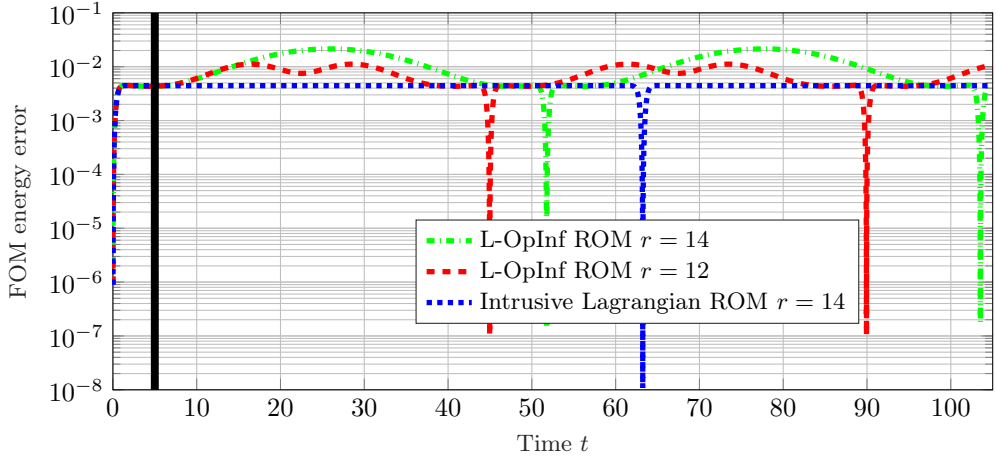}
 \caption{Sine-Gordon equation: The learned Lagrangian ROMs of reduced dimension $r=12$ and $r=14$ achieve bounded FOM energy error, which is remarkable because after $t=5$ the Lagrangian ROM simulations are purely predictive. The black line indicates end of training time interval at $t=5$.}
 \label{fig:SG_energy}
\end{figure}
Figure \ref{fig:SG_state} shows the relative state approximation error for intrusive and nonintrusive Lagrangian ROMs with increasing ROM order. The state error plots in Figure \ref{fig:SG_state_train} over the training time interval $[0, 5]$ show that the nonintrusive approach performs better than the intrusive Lagrangian ROM. For $r>10$, we observe the state error leveling off for the learned Lagrangian ROM where the accuracy does not improve with increasing reduced dimension in the training data regime. Figure~\ref{fig:SG_state_test} shows the relative state approximation error for the testing time interval $[5, 25]$. We can see from this figure that both intrusive and nonintrusive approaches yield Lagrangian ROMs with comparative accuracy up to $r=8$. For $r>8$, the intrusive Lagrangian ROM exhibits lower state error compared to the nonintrusive Lagrangian ROM in the testing data regime.
\par 
We have compared the approximate numerical solution using ROMs of size $r=14$ for both intrusive and nonintrusive approaches with the FOM solution in Figure \ref{fig:SG_traj}. Even though the reduced operators are learned from trajectory data in the training interval $[0,5]$, the nonintrusive Lagrangian ROM of reduced dimension $r=14$ captures the correct wave shape at $t=25$ which is $400\%$ past the training time interval. The ability of \textsf{L-OpInf} to provide accurate and stable predictions along with bounded FOM energy error  outside the training data for complex nonlinear wave phenomena is the key takeaway from this study. 
\par 
For the FOM energy error comparison, nonintrusive Lagrangian ROMs of size $r=12$ and $r=14$ are simulated until $t=105$ (which is $2,000\%$ past the training interval) to demonstrate the long-time stability of nonintrusive ROM simulations. The FOM energy error plots in Figure~\ref{fig:SG_energy} show that both intrusive and nonintrusive Lagrangian ROMs have similar bounded energy error behavior. Interestingly, the nonintrusive Lagrangian FOM energy error plots for both $r=12$ and $r=14$ change their qualitative behavior after leaving the training data regime but the error still remains bounded up to $t=105$. Despite the fact that nonintrusive reduced operators are learned purely from trajectory data up to $t=5$, the Lagrangian nature of our learned ROM ensures accurate prediction along with bounded energy error far outside the training data regime. 

\subsection{Soft-robotic fishtail model}
\label{sec:fishtail}
Autonomous underwater vehicles (AUVs) have become indispensable in a wide range of civilian and military applications. Due to their fast, agile and efficient underwater movement, AUVs mimicking a fish's swimming behavior have gained increasing interest for underwater surveillance and exploration. Soft-bodied robots designed for these applications utilize an actuation concept called fluid muscle, where the actuation is distributed over the whole fishtail via an array of fluid elastomers. The compliant nature of these soft-bodied robots allows them to achieve continuum motion and perform escape response maneuvers involving rapid body accelerations of very short duration.

\begin{figure}
\begin{subfigure}{.45\textwidth}
\begin{annotationimage}{width=\linewidth,height=4.5cm}{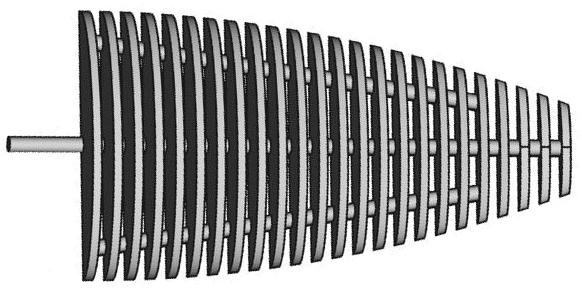}
       \draw[annotation above = {\textcolor{brown}{Main tubes} at 0.12}] to (0.08,0.5) ;
        \draw[annotation above = {\textcolor{olive}{Side tubes} at 0.5}] to (0.38,0.2) ;
        \draw[annotation above = { at 0.5}] to (0.6,0.7);
        \draw[annotation above = {\textcolor{magenta}{Chambers} at 0.8}] to (0.78,0.74) ;
        \draw[annotation above = { at 0.8}] to (0.9,0.68);
\end{annotationimage}
\caption{Fluid chamber system}
 \label{fig:fishtail_skeleton}
\end{subfigure}
\hspace{0.5 cm}
\begin{subfigure}{.45\textwidth}
\begin{annotationimage}{width=\linewidth,height=4.5cm}{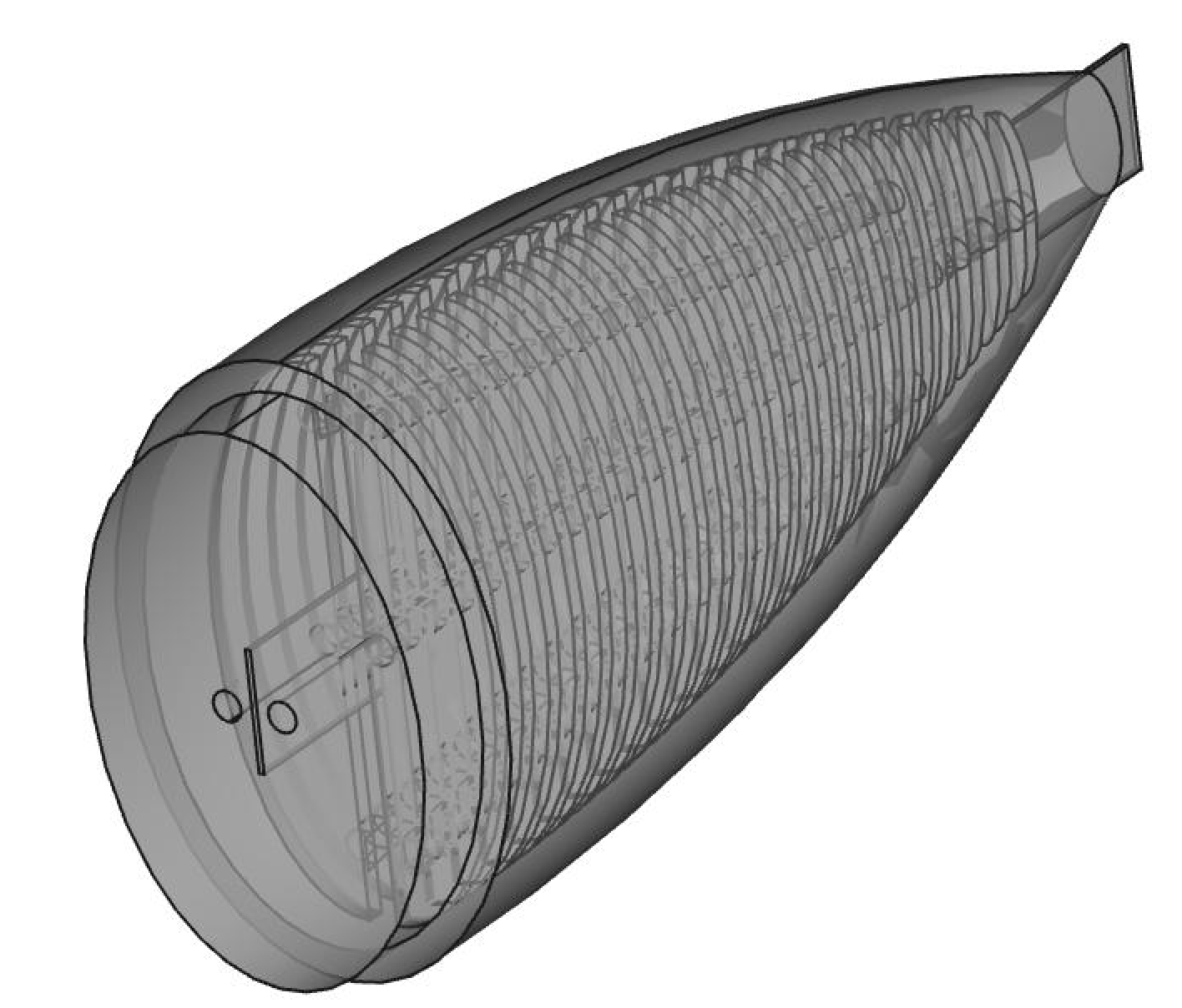}
\draw[thick,->] (0,0.2,0) -- (0.05,0.17,0.1) node[left=-0.07]{$z_1$};
\draw[thick,->] (0,0.2,0) -- (-0.03,0.1,-0.2) node[left=0.]{$z_2$};
\draw[thick,->] (0,0.2,0) -- (0.06000,0.2200,-0.0500) node[left=-0.07]{$z_3$};
       \draw[annotation above = {\textcolor{teal}{Carbon center beam} at 0.1}] to (0.19,0.3) ;
        \draw[annotation above = {\textcolor{orange}{Silicon hull} at 0.5}] to (0.6,0.83) ;
        \draw[annotation above = {\textcolor{violet}{POI} at 0.8}] to (0.98,0.84) ;
\end{annotationimage}
  \caption{Complete fishtail model}
  \label{fig:fishtail_CAD}
\end{subfigure}%
\caption{Soft-robotic fishtail: The model is based on a 3D model of the fishtail designed using the open-source software FreeCAD \cite{riegel2016freecad}. Original images shared by the authors of \cite{siebelts2018modeling}, copyright Elsevier.} 
\label{fig:fishtail_pics}
\end{figure}

The soft-robotic fishtail model considered here is based on the CAD model created in \cite{siebelts2018modeling}. More details about the design idea behind this model can be found in \cite{marchese2014autonomous}.
The fishtail design considered in this CAD model uses two separate fluid chambers to mimic the natural antagonistic muscle pair interaction. To minimize equipment in the experimental setup, both fluid chambers are pre-pressurized with $u_{\text{pre}}$ and a difference pressure $u(t)$ is used as the control input. This setup leads to an effective pressure of $u_{\text {pre}}+ u(t)$ in one chamber and $u_{\text {pre}}- u(t)$ in the other.
Figure \ref{fig:fishtail_skeleton} shows one of the two  identical fluid chambers where the main tube at the center along with the side tubes provide a single pressure supply to the chambers. The complete fishtail model shown in Figure \ref{fig:fishtail_CAD} is divided into two parts: (i) carbon center beam and (ii) silicone hull. The smooth silicone hull design is generated by three ellipses and a thick carbon beam is placed in the center of the silicone hull to increase stability of the soft-robotic fishtail against torsional moments.  The fluid chamber systems are placed on both sides of the carbon center beam.

In \cite{saak2019comparison}, various frequency-domain based second-order reduced-order modeling methods have been applied to the soft-robotic fishtail model. However, all the ROMs considered therein are intrusive in that they assume access to the FOM operators. In contrast, the proposed approach does not need access to the FOM operators as it learns the reduced operators nonintrusively directly from data. We only require solution of a linear least-squares problem with symmetric positive-definite constraints at the ROM level, whereas \cite{saak2019comparison} requires the solution of Lyapunov equations at the FOM level.

\subsubsection{PDE formulation}
\label{sec:fishtail_pde}
The governing PDE for both parts of the soft-robotic fishtail is derived using the theory of linear elastic materials; all equations in this subsection hold for both parts with different material properties. The first-order displacement tensor $\underline{q}(t,\bz)$ and the symmetric second-order strain tensor $\underline{\epsilon}(t,\bz)$ are related by the kinetic equation
\begin{equation*}
\underline{\epsilon}(t,\bz)=\frac{1}{2} \left( \nabla_{\bz} \underline{q}(t,\bz) + \nabla_{\bz}^{\top} \underline{q}(t,\bz) \right),
\end{equation*}
where $\bz=(z_1,z_2,z_3)$ is the three-dimensional spatial variable. Using the stress-strain relationship for isotropic elastic materials, the stress tensor $\underline{~\sigma}(t,\bz)$ is 
given by
\begin{equation*}
\underline{~\sigma}(t,\bz)=\lambda \trace{\underline{\epsilon}(t,\bz)} \underline{I} + 2\mu \underline{\epsilon}(t,\bz) ,
\end{equation*}
where $\trace{\cdot}$ is the trace operator, $\underline{I}$ is the identity tensor and  $\lambda$ and $\mu$ are the Lam\'{e} parameters. The governing PDE can be written in terms of the second time-derivative of the displacement and divergence of the stress tensor 
\begin{equation}
\rho \frac{\partial ^2 \underline{q}(t,\bz)}{\partial t^2} = \nabla_{\bz} \cdot \underline{~\sigma}(t,\bz),
\label{eq:fish_pde}
\end{equation}
where $\rho$ is the constant density of a material with a homogeneous mass distribution.
\subsubsection{FOM implementation}
\label{sec:fishtail_fom}
Space discretization of the governing PDE \eqref{eq:fish_pde} using tetrahedral finite elements yields
\begin{equation*}
\bM\ddot{\bq}(t)  + \bK\bq(t) = \bB u(t),
\end{equation*}
where $\bM,  \bK \in \real^{n \times n}$ are symmetric positive-definite matrices, $\bB \in \real^n$ is the input vector, and $u(t) \in \real $ in the scalar control input. Since the original mathematical model does not account for any dissipative effects, Rayleigh damping \cite{meirovitch2010fundamentals} is introduced in form of a damping matrix $\bC \in \real^{n \times n}$ to obtain realistic behavior of the soft structure, i.e.,
\begin{equation}
\bC=\alpha\bM + \beta \bK \in \real^{n \times n},
\end{equation}
where $\alpha=10^{-4}$ and $\beta=2 \times 10^{-4}$. The state vector $\bq$ for the FOM considered in this example has dimension $n=779,232$ which results from a computational grid of $259,744$ points and three spatial DOFs at each grid point. We consider a stationary initial condition, i.e. $\bq(0)=\dot{\bq}(0)=\bzero$, as the system is forced for $t>0$. 

\par

In order to achieve a fish-like movement, we require control over the displacement of the fish's fin. Thus, as an output, we are interested in the location ($0$~\si{\metre}, $0$~\si{\metre}, $0.21$~\si{\metre}) which is denoted as the point of interest (POI) and is equal to a single mesh point at the end of the center beam, see Figure \ref{fig:fishtail_CAD}. This then defines the output
\begin{equation*}
\by(t)=\bE \bq(t) \in \real^3,
\end{equation*}
which represents displacement of the rear tip of the carbon center beam in the three spatial directions.

\begin{remark}
The high-dimensional fishtail model reflects a typical gray-box setting where we have knowledge about the Lagrangian nature of the mechanical system, but details about the spatial discretization are unavailable. Unlike the beam example in Section~\ref{sec:beam}, the derivation of FOM system matrices for the fishtail model is labor-intensive as it requires a multitude of different software packages (having dependencies) and settings. Moreover, we do not have have access to the finite element code or meshing data used for generating these system matrices. The system matrices defining the Lagrangian FOM for this example are obtained from \cite{siebelts_dirk_2019_2558728}.
\end{remark}

The soft-robotic fishtail FOM is numerically integrated using the Newmark integrator with a fixed time step $\Delta t=0.001$~\si{\s}. Algorithm~\ref{alg:newmark} summarizes the Newmark integrator for the soft-robotic fishtail FOM with external control input. The Newmark integrator implementation for this FOM requires solving an $n$-dimensional linear system at every time step. We use the \texttt{SuiteSparse} package \cite{davis2014suitesparse} for an efficient FOM implementation in MATLAB. \texttt{SuiteSparse} is a suite of sparse matrix algorithms which exploits the sparsity in large matrices to achieve speedup in computations. Since the matrix $A_{\Delta t}$ in Algorithm~\ref{alg:newmark} is constant for all time steps, we have used the \texttt{Factorize} object in \texttt{SuiteSparse} for solving the large linear system. This object computes the Cholesky factorization of $A_{\Delta t}$ once and returns it as an object that can be reused for every linear solve. Numerical time integration of this large-scale FOM for total time $T=2$~\si{\s}  using $\Delta t=0.001$~\si{\s} requires approximately 110 minutes (MATLAB wall clock time) on a personal computer. Details about hardware and software used for these simulations can be found in Table~\ref{tab:computer}.

 \begin{algorithm}[h]
\caption{Newmark integrator for soft-robotic fishtail model with control input}
\begin{algorithmic}[1]
\Require Initial conditions $(\bq(0),\dot{\bq}(0))$, system matrices $\bM,\bK,\bC,\bB$, control input $u(t)$, simulation time $T$, and fixed time step $\Delta t$. 
\Ensure Discrete solution trajectories $\bQ$
    \State Compute total number of time steps $K=\frac{T}{\Delta t}$.
    \State Initialize $\ddot{\bq}_0:=\ddot{\bq}(0)$ by solving the following linear system of equations
    \begin{equation*}
        \bM \ddot{\bq}_0= \bB u(0) -\bC \dot{\bq}(0) - \bK \bq(0).
    \end{equation*}
    \While{$k \leq K$}
    \State Compute predicted `mean' values
\begin{equation*}
\breve{\bq}_{k+1}=\bq_k + \Delta t \dot{\bq}_k + \frac{\Delta t^2}{4}\ddot{\bq}_k, \quad \quad \dot{\breve{\bq}}_{k+1}=\dot{\bq}_k + \frac{\Delta t}{2}\ddot{\bq}_k.
\end{equation*}
\State Linear solve to obtain $\ddot{\bq}_{k+1}$
\begin{equation*}
\underbrace{\left( \bM + \frac{\Delta t}{2}\bC + \frac{\Delta t^2}{4}\bK \right)}_{\bA_{\Delta t}}\ddot{\bq}_{k+1} = \bB u_{k+1} - \bC \dot{\breve{\bq}}_{k+1} - \bK\breve{\bq}_{k+1}.
\end{equation*}
\State Explicit update equations to obtain $\bq_{k+1}$ and $\dot{\bq}_{k+1}$
\begin{equation*}
\bq_{k+1}= \breve{\bq}_{k+1} + \frac{\Delta t^2}{4}\ddot{\bq}_{k+1}, \quad \quad \dot{\bq}_{k+1}= \dot{\breve{\bq}}_{k+1} + \frac{\Delta t}{2}\ddot{\bq}_{k+1}.
\end{equation*}
\EndWhile
\State Assemble discrete trajectories to construct snapshot matrix
\begin{equation*}
    \bQ=[\bq_1, \cdots, \bq_{K}] \in \real^{n\times K}.
\end{equation*}
\end{algorithmic}
\label{alg:newmark}
\end{algorithm}

\begin{table}[h]
\caption{Soft-robotic fishtail implementation details}
\parbox{.45\linewidth}{
\centering
\begin{tabular}{|l|l|}
\hline
Processor & 2.3 GHz Intel Core i7\\
Cores & Quad-core \\
RAM & 32 GB 3733 MHz LPDDR4X \\
Operating system & macOS Catalina 10.15.7 \\
MATLAB  & 2020b\\
SuiteSparse &  5.10.1\\
CVX & 2.2 \\ 
\hline 
\end{tabular}
\label{tab:computer}
\caption*{a) Hardware and software specifications}
}
\hspace{1.6cm}
\parbox{.4\linewidth}{
\centering
\begin{tabular}{|c|l|}
\hline
Part & \quad \quad Parameter \\
\hline
Silicone hull & $\rho_1=1.07 \times 10^{-3}$~\si{\kilogram \per \metre\tothe{3}}    \\

 & $\lambda_1=2.03 \times 10^{5}$~\si{\newton \per \metre\tothe{2}}  \\

         & $\mu_1=8.45 \times 10^{3}$~\si{\newton \per \metre\tothe{2}} \\
\hline
Center beam & $\rho_2=1.40 \times 10^{3}$~\si{\kilogram \per \metre\tothe{3}} \\

 & $ \lambda_2=1.71 \times 10^{10}$~\si{\newton \per \metre\tothe{2}} \\

         & $\mu_2=1.14 \times 10^{10}$~\si{\newton \per \metre\tothe{2}}\\ 
\hline

\end{tabular}
\vspace{0.05cm}
\caption*{b) Material parameters}
}
\end{table}

The numerical studies presented in this section consider several different inputs (forcing) on the system to demonstrate that \textsf{L-OpInf} performs well for a range of different dynamics. In particular, we consider:
\begin{itemize}
    \item Step input $ u(t) =\left\{ \begin{matrix}
     0 & \text{if }t < 0.1 \\
      5000 & \text{if }t \geq 0.1
    \end{matrix}\right.\, $.
    \item Sinusoidal input $ u(t) =\left\{ \begin{matrix}
     0 & \text{if }t < 0.1 \\
      2500(~\sin(10\pi(t-1.75))+1) & \text{if }t \geq 0.1
    \end{matrix}\right.\,$.
    \item Ramp input $ u(t) =\left\{ \begin{matrix}
     50000t & \text{if }t < 0.1 \\
      5000 & \text{if }t \geq 0.1
    \end{matrix}\right.\, $.
    \item Sigmoid input $ u(t) =\left\{ \begin{matrix}
    \frac{5000\sqrt{17}}{4}\frac{(t/0.025)}{ \sqrt{1 + (t/0.025)^2}} & \text{if }t < 0.1 \\
      5000 & \text{if }t \geq 0.1
    \end{matrix}\right.\, $.
\end{itemize}

\subsubsection{Results}
\label{sec:fishtail_results}
The main motivation for ROM development for the soft-robotic fishtail is to make real time predictions on the limited hardware of an AUV. Therefore, we consider ROMs of size $r=2$ for comparison between the intrusive Lagrangian ROMs and nonintrusive learned ROMs. Compared to the approximate FOM run time of 110 minutes,  numerical time integration of
the learned Lagrangian ROM of size $r=2$ for total time $T = 2$~\si{\s} using $\Delta t=0.001$~\si{\s} requires approximately 0.0130 seconds
(MATLAB wall clock time averaged over 20 runs), which is a factor of 507,692x speedup. For illustration purposes, we show the second component ($z_2$ direction) of the output $\by(t)$ in the plots below.

\begin{figure}[h]
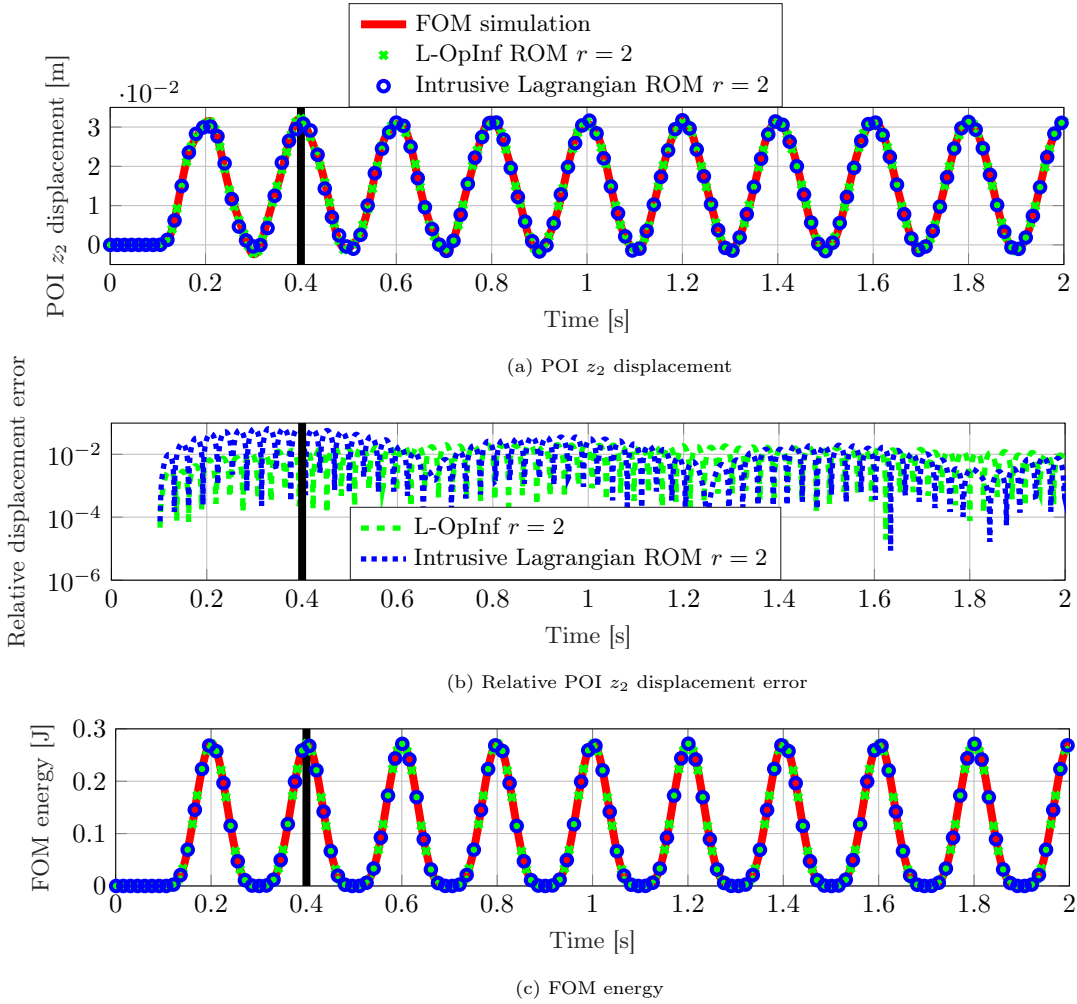

\captionsetup[subfigure]{oneside,margin={1cm,0 cm}}
\begin{subfigure}{\textwidth}
       \setlength\fheight{3.9 cm}
        \setlength\fwidth{\textwidth}
\hspace{5 mm}\input{Figures/fishtail/sine_POI.tex}
\caption{POI $z_2$ displacement}
\label{fig:fish_sine_poi}
    \end{subfigure}\\
    \vspace{-0.5cm}
\captionsetup[subfigure]{oneside,margin={1.2cm,0 cm}}

    \begin{subfigure}{\textwidth}
           \setlength\fheight{3.9 cm}
           \setlength\fwidth{\textwidth}
\input{Figures/fishtail/sine_POI_err.tex}
\caption{Relative POI $z_2$ displacement error}
\label{fig:fish_sine_poi_err}
    \end{subfigure} \\
   \captionsetup[subfigure]{oneside,margin={0.2cm,0 cm}}
\begin{subfigure}{\textwidth}
           \setlength\fheight{3.9 cm}
           \setlength\fwidth{\textwidth}
\hspace{3 mm}\input{Figures/fishtail/sine_energy.tex}
\caption{FOM energy}
\label{fig:fish_sine_poi_energy}
    \end{subfigure}
\caption{Soft-robotic fishtail: Both intrusive and nonintrusive Lagrangian ROMs capture the POI $z_2$ displacement accurately when the sinusoid input was used for training and testing. The relative displacement error in (b) shows that the \textsf{L-OpInf} ROM performs marginally better than the intrusive Lagrangian ROM in the training data regime $[0, 0.4]$~\si{\s}. Both ROMs track the change in FOM energy accurately far outside the training data. The black line indicates the end of  the training time interval.}
 \label{fig:fish_sine}
\end{figure}

\paragraph{Prediction outside of training time interval} 
We demonstrate the ability of the \textsf{L-OpInf} ROMs to predict model behavior and have accurate energy behavior far outside the training interval. Figure~\ref{fig:fish_sine} compares the performance of both intrusive and nonintrusive ROMs for the sinusoidal input case. For this case, snapshot data from the response of the soft-robotic fishtail model to the sinusoidal input over training interval $[0,0.4]$~\si{\s} is used to train the nonintrusive Lagrangian ROM of the form \eqref{eq:rom_learn}. The testing data consists of snapshot data from the response of the FOM to the sinusoidal input over the testing interval $[0.4, 2]$~\si{\s}. We have compared the POI $z_2$ displacement approximation from both ROMs with the FOM simulation in Figure~\ref{fig:fish_sine_poi}, and both approaches yield Lagrangian ROMs of comparative accuracy. The relative error comparison in Figure~\ref{fig:fish_sine_poi_err} shows that the learned ROM performs better than the intrusive ROM in the training data regime whereas after $t=0.6$~\si{\s} both ROMs exhibit relative error of approximately $10^{-2}$. We compare the total energy from FOM simulations with the system energy approximation from both methods in Figure~\ref{fig:fish_sine_poi_energy}. Even though the reduced operators are learned from training data $[0, 0.4]$~\si{\s}, the nonintrusive Lagrangian ROM tracks the time-varying energy accurately at $t=2$~\si{\s} which is $400\%$ outside the training time interval. 

\begin{figure}[h]
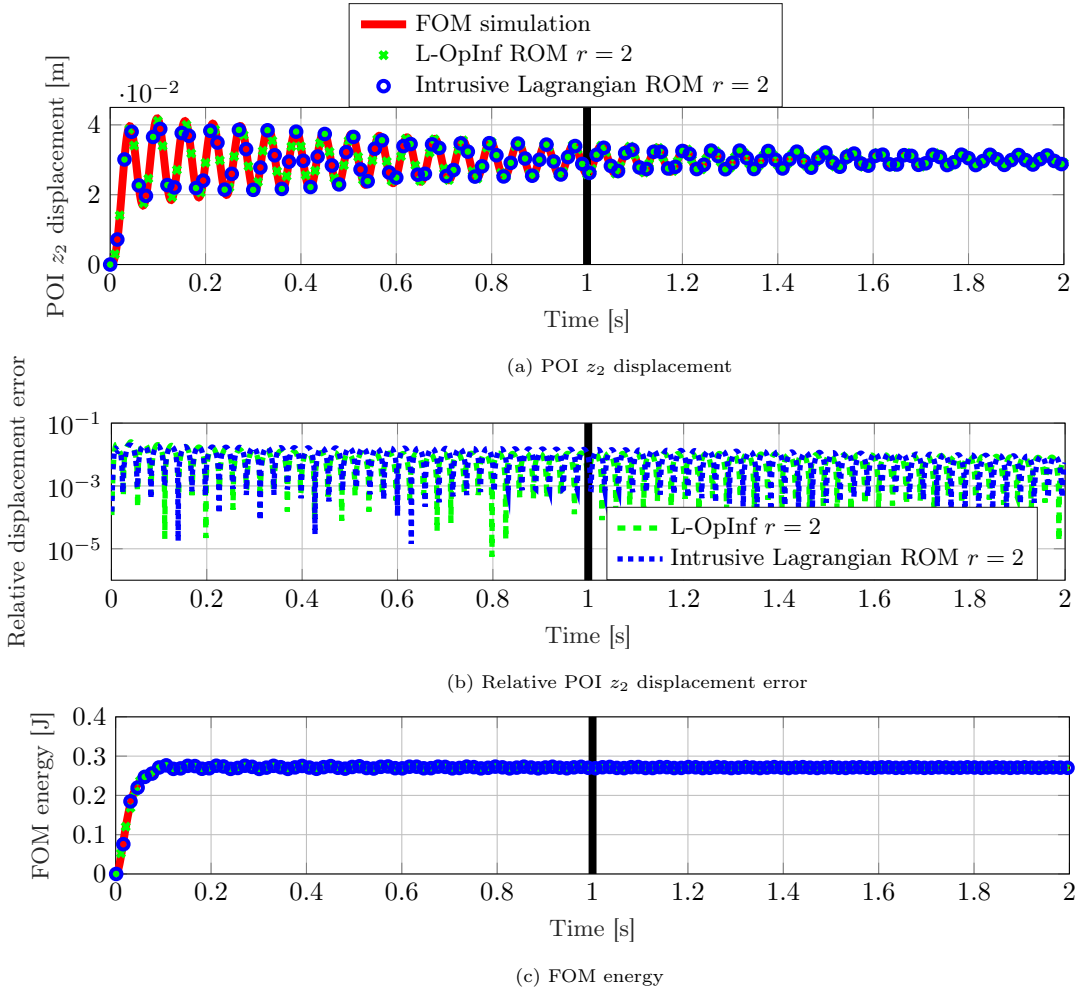

\captionsetup[subfigure]{oneside,margin={1cm,0 cm}}
\begin{subfigure}{\textwidth}
       \setlength\fheight{3.9 cm}
        \setlength\fwidth{\textwidth}
\hspace{5mm}\input{Figures/fishtail/sigmoid_POI.tex}
\caption{POI $z_2$ displacement}
\label{fig:fish_sigmoid_poi}
    \end{subfigure} \\
 \vspace{-0.5cm}
 \captionsetup[subfigure]{oneside,margin={1.2cm,0 cm}}

\begin{subfigure}{\textwidth}
           \setlength\fheight{3.9 cm}
           \setlength\fwidth{\textwidth}
\input{Figures/fishtail/sigmoid_POI_err.tex}
\caption{Relative POI $z_2$ displacement error}
\label{fig:fish_sigmoid_poi_err}
    \end{subfigure} \\
    \captionsetup[subfigure]{oneside,margin={0.2cm,0 cm}}

\begin{subfigure}{\textwidth}
           \setlength\fheight{3.9 cm}
           \setlength\fwidth{\textwidth}
\hspace{3mm}\input{Figures/fishtail/sigmoid_energy.tex}
\caption{FOM energy}
\label{fig:fish_sigmoid_energy}
    \end{subfigure}
\caption{Soft-robotic fishtail: When the sigmoid input was used for training and testing, the \textsf{L-OpInf} ROM with $r=2$ achieves similar output error performance to that of intrusive Lagrangian ROM with $r=2$ with relative error of approximately $10^{-2}$. Both ROMs also track the change in FOM energy accurately $100\%$ outside the training data regime. The black line indicates the end of the training time interval.}
 \label{fig:fish_sigmoid}
\end{figure}

In Figure~\ref{fig:fish_sigmoid}, we compare the numerical performance of intrusive and nonintrusive Lagrangian ROMs for the sigmoid input case. Both intrusive and nonintrusive Lagrangian ROMs are trained with snapshot data from the response of the soft-robotic fishtail model to the sigmoid input over training interval $[0,1]$~\si{\s}. The test data is generated from the response of the FOM to the sigmoid input over the testing interval $[1, 2]$~\si{\s}. The POI $z_2$ displacement plot in Figure~\ref{fig:fish_sigmoid_poi} shows that both in intrusive and learned ROM approximations agree with the FOM simulation. Both approaches exhibit a relative error of $10^{-2}$ for the entire simulation in Figure~\ref{fig:fish_sigmoid_poi_err}. Figure~\ref{fig:fish_sigmoid_energy} demonstrates that preserving the underlying Lagrangian structure leads to accurate FOM energy prediction for the sigmoid input case. 

\begin{figure}[h]
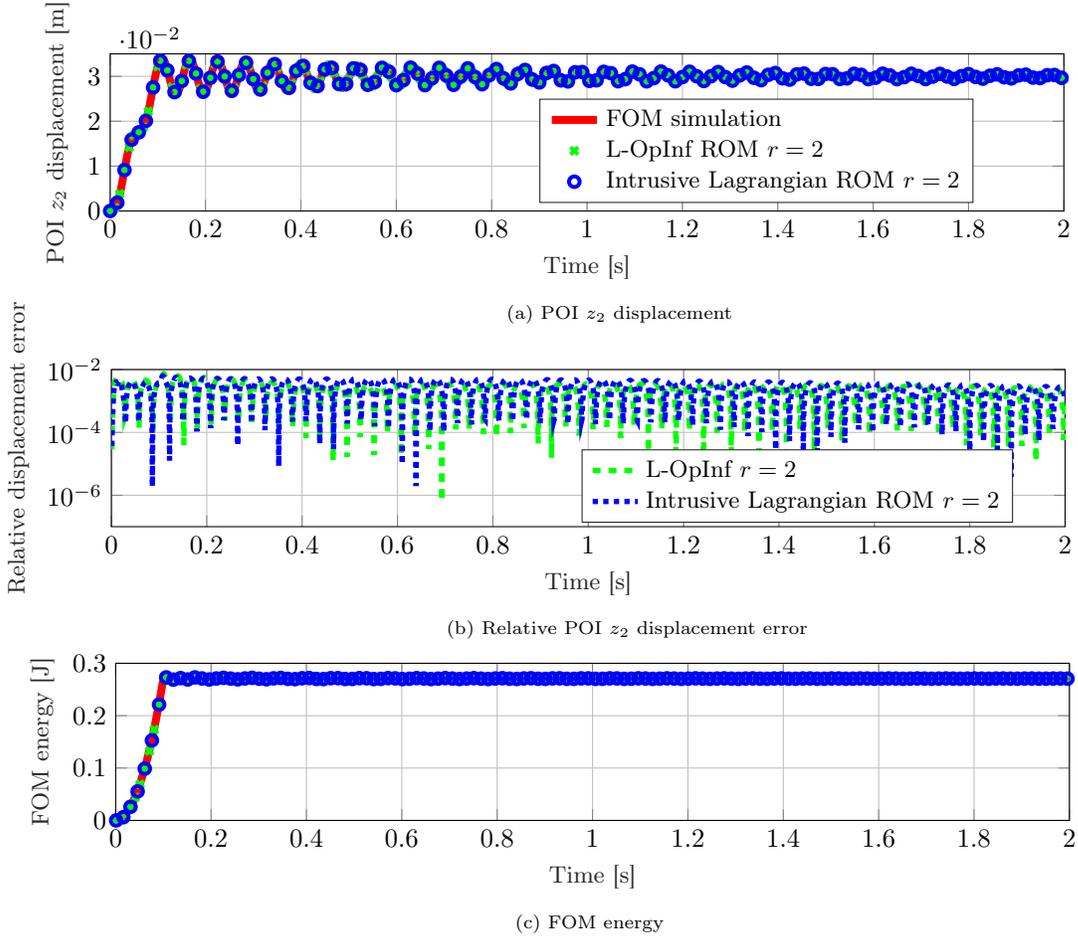

\captionsetup[subfigure]{oneside,margin={1cm,0 cm}}
\begin{subfigure}{\textwidth}
       \setlength\fheight{3.9 cm}
        \setlength\fwidth{\textwidth}
\hspace{5mm}\input{Figures/fishtail/ramp_unknown_sigmoid_POI.tex}
\caption{POI $z_2$ displacement}
\label{fig:fish_ramp_poi}
    \end{subfigure} \\
 \vspace{-0.5cm}
 \captionsetup[subfigure]{oneside,margin={1.2cm,0 cm}}
 
\begin{subfigure}{\textwidth}
           \setlength\fheight{3.9 cm}
           \setlength\fwidth{\textwidth}
\input{Figures/fishtail/ramp_unknown_sigmoid_POI_err.tex}
\caption{Relative POI $z_2$ displacement error}
\label{fig:fish_ramp_poi_err}
    \end{subfigure} \\
    \captionsetup[subfigure]{oneside,margin={0.2cm,0 cm}}
    
\begin{subfigure}{\textwidth}
           \setlength\fheight{3.9 cm}
           \setlength\fwidth{\textwidth}
\hspace{3mm}\input{Figures/fishtail/ramp_unknown_sigmoid_energy.tex}
\caption{FOM energy}
\label{fig:fish_ramp_energy}
    \end{subfigure}
\caption{Soft-robotic fishtail: Both intrusive and nonintrusive Lagrangian ROMs trained using the sigmoid input exhibit similar performance for the ramp input test case. The relative POI $z_2$ displacement error for both ROMs stays below $10^{-2}$ for the entire simulation.} 
\label{fig:fish_ramp}
\end{figure}

\paragraph{Prediction for unseen inputs} 
We demonstrate the ability of the learned ROM to generalize to new (unseen) inputs. Having a surrogate model that is robust to unknown inputs is desirable in control applications where we can not foresee what control input the system will be subjected to. The \textsf{L-OpInf} ROMs are learned with data from the sigmoid input.
First, we consider the ramp input as a testing input in Figure~\ref{fig:fish_ramp}. The POI $z_2$ displacement comparison in Figure~\ref{fig:fish_ramp_poi} shows that the learned ROM predicts the soft-robotic fishtail model behavior under the ramp input accurately. The relative error comparison in Figure~\ref{fig:fish_ramp_poi_err} shows that the relative error of both the intrusive and nonintrusive Lagrangian ROMs remains below $10^{-2}$ for the entire simulation. The FOM energy in Figure~\ref{fig:fish_ramp_energy} shows that the \textsf{L-OpInf} ROM captures the initial increase in the energy accurately before settling into a constant energy state.

\begin{figure}[h]
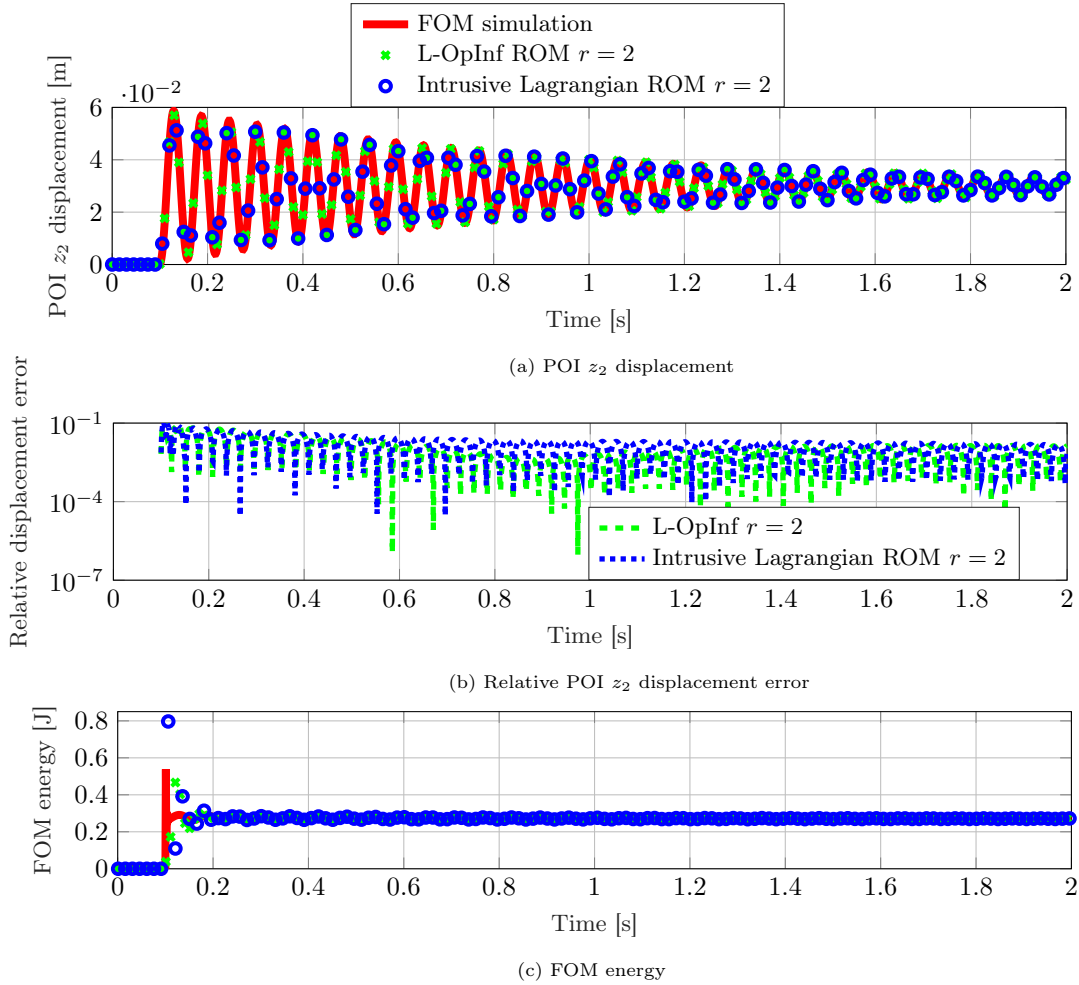

\captionsetup[subfigure]{oneside,margin={1cm,0 cm}}
\begin{subfigure}{\textwidth}
       \setlength\fheight{3.9 cm}
        \setlength\fwidth{\textwidth}
\hspace{5 mm}\input{Figures/fishtail/step_unknown_sigmoid_POI.tex}
\caption{POI $z_2$ displacement}
\label{fig:fish_step_poi}
    \end{subfigure} \\
 \vspace{-0.5cm}
 \captionsetup[subfigure]{oneside,margin={1.2cm,0 cm}}
 
\begin{subfigure}{\textwidth}
           \setlength\fheight{3.9 cm}
           \setlength\fwidth{\textwidth}
\input{Figures/fishtail/step_unknown_sigmoid_POI_err.tex}
\caption{Relative POI $z_2$ displacement error}
\label{fig:fish_step_poi_err}
    \end{subfigure} \\
     \captionsetup[subfigure]{oneside,margin={0.2cm,0 cm}}
\begin{subfigure}{\textwidth}
           \setlength\fheight{3.9 cm}
           \setlength\fwidth{\textwidth}
\hspace{3mm}\input{Figures/fishtail/step_unknown_sigmoid_energy.tex}
\caption{FOM energy}
\label{fig:fish_step_energy}
    \end{subfigure}
\caption{Soft-robotic fishtail: For the step input test case, both intrusive and nonintrusive ROMs trained using the sigmoid input case capture the POI $z_2$ displacement behavior accurately. Both Lagrangian ROMs exhibit similar displacement error with the \textsf{L-OpInf} ROM performing marginally better. In the FOM energy plot, both ROMs fail to capture the sudden jump in the FOM energy when the control input is turned on at $t=0.1$~\si{\s}.}
 \label{fig:fish_step}
\end{figure}

We have used the learned Lagrangian ROM learned from the snapshots of the sigmoid input case to study the step response of the soft-robotic fishtail model in Figure~\ref{fig:fish_step}. Unlike the other three input cases with smooth actuation, the step input case involves a non-smooth actuation at $t=0.1$~\si{\s}. The displacement plots in Figure~\ref{fig:fish_step_poi} shows how well the learned ROM predicts the dominant motion of the soft-robotic fishtail in the $z_2$ direction. The relative error plots show a sharp jump in the displacement error for both approaches at $t=0.1$~\si{\s} when the control input is applied. Relative errors for both methods gradually decline from $t=0.1$~\si{\s} to $t=0.6$~\si{\s} before settling to a relative error of approximately $10^{-2}$ with the learned ROM performing marginally better. Unlike the other three input cases, both intrusive and nonintrusive Lagrangian ROMs fail to capture the FOM energy in the transient phase in Figure~\ref{fig:fish_step_energy}. The FOM energy predicted by both reduced models shows oscillatory behavior before converging to the energy predicted by the FOM simulation.

\section{Conclusions} \label{sec:conclusions}
We presented a data-driven model reduction method,  \textit{Lagrangian operator inference} (\textsf{L-OpInf}), that learns Lagrangian ROMs directly from high-dimensional data via nonintrusive structure-preserving operator inference. The data we considered came from large-scale semi-discretized PDE models that are marched forward in time using a time integrator.~The presented \textsf{L-OpInf} method exploits the underlying geometric structure of Lagrangian systems and applies to mechanical systems with nonconservative forcing and nonlinear wave equations.~The inference of the reduced operators is based on a constrained optimization problem that ensures that the reduced models are Lagrangian and also respect the symmetric (positive definite) property of system matrices. For Lagrangian systems that do not fall into the class of Lagrangian FOMs considered in Section~\ref{sec:FOM}, our work may be used as a template to define a problem-specific constrained optimization problem that yields Lagrangian ROMs that preserve the structure of the high-dimensional Lagrangian system.
The proposed method only assumes a Lagrangian nature of the large-scale system and does not require access to FOM operators or require information about the spatial discretization used to derive the high-dimensional systems. This setting is common in many applications where the full model is either given as a black box or the source code is very complicated and understanding these implementation details is tedious and time consuming. In those settings, intrusive model reduction may be less appealing to novice users. 

\par 

The numerical experiments demonstrate the advantages of preserving the underlying geometric structure for both conservative and forced Lagrangian systems. The conservative Euler-Bernoulli beam example shows that the proposed method learns stable ROMs with bounded energy error, while facilitating accurate long-time predictions far outside the training time interval and demonstrating robustness to unseen initial conditions that are outside the training dataset. The nonlinear sine-Gordon wave equation example shows that the proposed method can learn accurate and stable ROMs of semi-discretized nonlinear wave equations purely from  the trajectory data. For the high-dimensional soft-robotic fishtail model, the learned Lagrangian ROMs track the change in system energy accurately in the presence of dissipation and time-dependent control input. Notably, the learned Lagrangian ROMs work well even for unknown control inputs. Moreover, the proposed method achieves significant reduction in state dimension, which makes the learned ROMs ideal for real-time control and state estimation. 

\par 

Future research directions motivated by this work are: combining \textsf{L-OpInf} with structure-preserving machine learning methods to learn Lagrangian ROMs of dynamical systems with unknown nonlinear operators; studying connection between reduced operators and predicted energy of the intrusive and nonintrusive ROMs; and applying the proposed method to noisy data coming from experiments instead of the simulated data setting considered in this work.
\section*{Acknowledgments} We thank Prof. Thomas Meurer and Dirk Wolfram for sharing the images in Figure \ref{fig:fishtail_pics}. This research was in part financially supported by the Ministry of Trade, Industry and Energy (MOTIE) and the Korea Institute for Advancement of Technology (KIAT) through the International Cooperative R\&D program (No.~P0019804, Digital twin based intelligent unmanned facility inspection solutions) and the U.S. Office of Naval Research under award number N00014-22-1-2624.

\bibliography{references}
\bibliographystyle{vancouver}

\end{document}